\documentclass[epsf,11pt]{article}

\usepackage{chngcntr}
\counterwithin*{equation}{section}

\usepackage{latexsym,amssymb,amsmath}
\usepackage[dvips]{graphicx}
\usepackage{psfrag}
\usepackage[usenames]{color}
\usepackage{amsmath}
\newcommand{\BYDEF}{\,\shortstack{\textup{\tiny def} \\ = }\,}

\newcommand {\R} {\rm I\!R}

\newcommand{\nexto}{\kern -0.54em}

\newcommand{\dZ}{{\cal Z \kern -0.7em Z}}
\newcommand{\dC}{{\rm\hbox{C \kern-0.8em\raise0.2ex\hbox{\vrule height5.4pt width0.7pt}}}}
\newcommand{\dQ}{{\rm\hbox{Q \kern-0.85em\raise0.25ex\hbox{\vrule height5.4pt width0.7pt}}}}

\newtheorem{Lemma}{Lemma}
\newtheorem{Theorem}{Theorem}
\newtheorem{Corollary}{Corollary}
\newtheorem{Proposition}{Proposition}

\newtheorem{Definition}{Definition}
\newtheorem{Remark}{Remark}
\newtheorem{Assumption}{Assumption}

\newcommand{\beq}{\begin{equation}}
\newcommand{\deq}{\end{equation}}
\newcommand{\baq}{\begin{eqnarray}}
\newcommand{\daq}{\end{eqnarray}}

\newcommand{\baqm}{\begin{eqnarray*}}
\newcommand{\daqm}{\end{eqnarray*}}

\newcommand\old[1]{}
\usepackage{epsfig}
\usepackage{latexsym}

\title{On  near optimal control  of  systems with slow observables}
\author{Vladimir Gaitsgory\thanks{Department of Mathematics, Macquarie University, Sydney 2109, Australia, vladimir.gaitsgory@mq.edu.au; The work of V. Gaitsgory was supported by the Australian Research Council Discovery Grants DP130104432}
\and
Sergey Rossomakhine
\thanks{Flinders Mathematical Sciences Laboratory, School of Computer Science, Engineering and Mathematics, Flinders University, GPO Box 2100, Adelaide SA 5001, Australia, serguei.rossomakhine@flinders.edu.au; The work of S. Rossomakhine was
supported by the Australian Research Council Discovery Grant
DP150100618}}
\begin{document}

\maketitle
\begin{abstract}
The paper deals with a problem of  control of a system characterized by the fact that the influence of controls on the dynamics of certain functions of state variables (called observables) is relatively weak and the rates of change of these observables are much slower than the rates of change of the state variables themselves.   The contributions of the paper are twofold. Firstly,  the averaged  system whose solutions   approximate  the trajectories of the slow observables is introduced, and
it is shown that the optimal value of the problem of optimal control
 with time discounting criterion considered on the solutions of the system with slow observables  (this problem is referred to as perturbed) converges to the optimal value of the corresponding problem
 of optimal control of the averaged system. Secondly,  a way how  an asymptotically optimal control of the  perturbed problem can be constructed on the basis of an optimal solution of the averaged problem is indicated,  sufficient and necessary optimality conditions for the averaged problem are stated, and  a way how a near optimal solution of the latter can be constructed numerically is outlined (the construction being illustrated with an example).
\end{abstract}

\bigskip

{\bf Key words.} Control of slow observables, Averaged system,  Average control generating families, Occupational measures, Numerical solution

\bigskip

{\bf AMS subject classifications.} 34E15, 34C29, 34A60, 93C70

\bigskip
\bigskip

\section{Introduction and preliminaries}\label{Introduction}

In this paper, we will study a problem of  control of a system characterized by the fact that
the influence of controls on the dynamics of certain functions of state variables
 (called {\it observables}) is relatively weak and the rates of change of these observables are much slower than the rates of change of the state variables themselves. Formally, a system with slow observables can be viewed as one embedded into a family of systems depending on a small parameter. When the parameter is set to zero, the observables  remain constant. With non-zero but small positive values of the parameter, the observables may move slowly, and their dynamics can be influenced  by controls.

 In more detail, we will be considering the control system
\begin{equation}\label{e-CSO-1}
\frac{dy_{\epsilon}(\tau)}{d\tau}= f(u(\tau), y_{\epsilon}(\tau)) + \epsilon g(u(\tau), y_{\epsilon}(\tau)) , \ \ \ \ \ \ y_{\epsilon}(0)=y_0,
\end{equation}
where  $\epsilon > 0$ is a small parameter; the vector functions $\ f(u,y):
     U\times \mathbb{R}^m \to \mathbb{R}^m$  and
  $\ g(u,y): U\times \mathbb{R}^m \to \mathbb{R}^m$
are continuous  and satisfy Lipschitz conditions in $y$ (uniformly in $u\in U$);  the controls $u(\cdot) $ are  measurable functions of time  satisfying  the inclusion $\ u(t) \in U $
($U$ being a given  compact subset of a finite dimensional space). The system (\ref{e-CSO-1}) will be referred to as the {\it perturbed system}.

Along with the perturbed system (\ref{e-CSO-1}), we will be considering the {\it reduced system}
\begin{equation}\label{e-CSO-2}
\frac{dy(\tau)}{d\tau}= f(u(\tau), y(\tau)),
\end{equation}
 (obtained from (\ref{e-CSO-1}) by taking $\epsilon = 0$), and everywhere in what follows, it is assumed that

{\it The reduced system (\ref{e-CSO-2})
 has $k$ constants of motion or, more precisely, there exists a  vector function
$F(y):\mathbb{R}^m \to \mathbb{R}^k$ such that
\begin{equation}\label{e-CSO-3}
F'(y)f(u,y)= 0 \ \ \ \ \ \forall u\in U, \ \ \forall y\in Y,
\end{equation}
where $F'(\cdot)$ is the Jacobian matrix and $Y$ is a compact subset of $\R^m$ in a neighbourhood of which the partial derivatives of the components of $F(\cdot) $ exist and satisfy Lipschitz conditions.}

From (\ref{e-CSO-3}) it follows that any solution $y(\tau)$ of the reduced system (\ref{e-CSO-2}) that has the following property: if $F(y(0))=z$, then $F(y(\tau))=z\ \  \forall \ \tau\geq 0 $. That is,
the set
\begin{equation}\label{e-CSO-4}
Y_z\BYDEF \{y\in \R^m \ : \ F(y)=z\}
\end{equation}
 is  invariant with respect to the reduced  system. (In what follows, it will be assumed that $Y_z\subset Y$ for $z$ from a sufficiently large area.)  This  invariance property does not hold true for the solutions of the perturbed
 system (\ref{e-CSO-1}). In fact, let $y_{\epsilon} (\tau)$ be the solution of (\ref{e-CSO-1}) obtained with  a control $u(\cdot)$ such that it is contained in $Y$ for all $\tau \geq 0$.  Then, taking the derivative of
 $z_{\epsilon}(\tau)\BYDEF F(y_{\epsilon}(\tau))$, one arrives at the  equation
\begin{equation}\label{e-CSO-0}
\frac{dz_{\epsilon}(\tau)}{d\tau}= \epsilon h(u(\tau), y_{\epsilon}(\tau)), \ \ \ \ \ z_{\epsilon}(0)= z_0\BYDEF F(y_0),
\end{equation}
where
\begin{equation}\label{e-CSO-0-1}
h(u,y)\BYDEF F'(y)g(u,y).
\end{equation}
(Below, we will be interested only in those solutions of (\ref{e-CSO-1}) that are contained in $Y$ for  sufficiently small $\epsilon$.) The components of the vector  $z_{\epsilon}(\tau)$ are called  {\it observables}.
Note that the observables are changing their values slowly but these  changes can be significant if the time period is long enough.

Systems with slow observables were studied mostly in the uncontrolled setting (see, e.g., \cite{Arnold}, \cite{Art-Kevrekidis},  \cite{FroyGottHam}, \cite{Givon}, \cite{Fima}, \cite{Grigorius}, \cite{Sanders}) with the only exception being problems of control of singularly perturbed Markov chains  (see \cite{Yin} and references therein). The latter  present an important  special case  of problems that are dealt with in this paper. It  is characterized by that $f(u,y)=A(u)y$ and $g(u,y)=B(u)y$, where $A(u), B(u)$ are matrix functions having some special properties (see Remark \ref{Remark-SP-Markov}  below).

Note that
 the system consisting of (\ref{e-CSO-1}) and (\ref{e-CSO-0}) (for convenience,  it  will be referred to as  the {\it augmented perturbed system}) is a special case of the system
\begin{equation}\label{e-CSO-0-101}
\frac{dy_{\epsilon}(\tau)}{d\tau}= f_1(\epsilon , u(\tau), y_{\epsilon}(\tau)), \ \ \ \ \ \  \frac{dz_{\epsilon}(\tau)}{d\tau}= \epsilon f_2(\epsilon , u(\tau), y_{\epsilon}(\tau), z_{\epsilon}(\tau)),
\end{equation}
in which $f_1(\cdot) $ and $f_2(\cdot)$ are smooth in $\epsilon$.
Being rewritten in the \lq\lq shrinked" time scale $t=\epsilon\tau $, such systems are commonly called weakly coupled singularly perturbed ({\bf SP}) systems.
Problems of control of SP  systems (including weakly coupled SP systems)
have received a great deal of attention (see \cite{Alv}, \cite{Art2}, \cite{Art3}, \cite{Art1}, \cite{AG}, \cite{Ben}, \cite{Bor-Gai-1}, \cite{Dmitriev}, \cite{Dont-Don}, \cite{Fil4}, \cite{Gai1}, \cite{Gai8},
\cite{Kab2}, \cite{Kok1}, \cite{Kus3}, \cite{Nai}, \cite{Reo1}, \cite{plot}, \cite{Marc-New}, \cite{QW}, \cite{Serea}, \cite{Vel}, \cite{Nai-2} and references therein).
However, to the best of our knowledge, in all  the earlier works, except  the recent paper \cite{Marc-New},
 the SP control systems were considered under the  assumption that the set of occupational measures generated by the controls and the corresponding solutions of the reduced system
 \begin{equation}\label{e-CSO-0-101-1}
\frac{dy(\tau)}{d\tau}= f_1(0 , u(\tau), y(\tau))
\end{equation}
is independent of the initial conditions  from a sufficiently large set (the validity of this being guaranteed by imposing certain  stability or controllability type conditions).
In \cite{Marc-New}, a weakly coupled SP system similar to (\ref{e-CSO-0-101}) was considered under the assumption that the reduced  system
(\ref{e-CSO-0-101-1})
satisfies the non-expansivity condition introduced in \cite{Marc-J}, the latter does not exclude the situation when  the set of occupational measures generated by the controls and the solutions of  (\ref{e-CSO-0-101-1}) depends  on the initial conditions.

In our case, the set of occupational measures generated by the controls and the solutions of the reduced system (\ref{e-CSO-2}) does depend on the initial conditions. In fact, due to the fact that the set $Y_z$ is invariant with respect to the system, any  measure generated by a control $u(\cdot)$ and the corresponding solution $y(\cdot)$ of (\ref{e-CSO-2}) will be supported on $U\times Y_z$ if $y(0)\in Y_z$. That is, the sets of occupational measures  corresponding to the initial conditions in $Y_{z_1}$ and in $Y_{z_2}$ ($z_1\neq z_2$) are different. Also, the system (\ref{e-CSO-2}) may not satisfy the non-expansivity condition of \cite{Marc-J}. Thus, the available techniques are generally not applicable to the class of problems that we are considering.

The contributions of the paper are twofold. Firstly,  the averaged  system whose solutions approximate  the trajectories of the slow observables is introduced, and
it is shown that the optimal value of the  optimal control problem
\begin{equation}\label{Vy-perturbed-streched}
\displaystyle{
  \epsilon \inf _{u(\cdot)  }\int_0 ^ { + \infty }e^{-(C\epsilon \tau)}r(u(\tau), y_{\epsilon}(\tau))d\tau} \stackrel{\rm def}{=} R^*(\epsilon, y_0) ,
\end{equation}
where $C >0$ is a discount rate  and  $inf$ is sought over the controls and the corresponding solutions of the system (\ref{e-CSO-1}), converges to the optimal value of an optimal control problem considered
 on the solutions  of the averaged system. (These two problems will be referred to as {\it perturbed} and {\it averaged} problems respectively.)  Note that the function $r(u,y): U\times \R^m \mapsto \R$ is assumed to be continuous in $(u,y)$ and
 satisfying Lipschitz conditions in $y$ (uniformly with respect to $u\in U$) and that the multiplier $\epsilon$ before the integral in (\ref{Vy-perturbed-streched}) is introduced for the values of the objective function to remain bounded with $\epsilon$ tending to zero.

 Secondly,
  a way how  an optimal or near optimal solution of the averaged problem can be found is described and a way how the latter can be used for construction of asymptotically optimal/near optimal solution of the perturbed  problem is indicated.
  Note that the perturbed problem may  become very ill-conditioned with small values of the parameter $\epsilon$ which makes it  difficult (and sometimes even impossible) to  tackle  with the  existing  optimisation techniques.
  In the framework of the approach that we are proposing, a near optimal control
of the perturbed problem is constructed on the basis of an optimal solution of the averaged problem that is independent of $\epsilon$.
Thus, the difficulties caused by the ill-conditioning are eliminated, and, in fact, the approximation  of the optimal control
 improves as $\epsilon$ decreases towards $0$.

 The paper consists of seven sections and is organized as follows. Section \ref{Introduction} is this introduction. In Section \ref{Sec-Averaged system}, we introduce the averaged system and the averaged optimal control problem, and we establish that, under the validity of certain assumptions, these can be used for an approximation of the perturbed system and the perturbed optimal control problem. In Section \ref{S-SC-2}, the general assumptions used in Section \ref{Sec-Averaged system} are shown to be valid if certain stability or controllability conditions are satisfied on the invariant sets of the reduced system.
 In Section \ref{S-SC1}, the assumptions of Section \ref{Sec-Averaged system} are verified to be true for
  a special case when the reduced system is uncontrolled.
 In Section \ref{Sec-ACG}, we give a  definition of  an average control generating (ACG) family and establish
 sufficient and necessary conditions for an ACG family to be optimal in the averaged problem, and we also outline an algorithm for finding near optimal ACG families (based on results of \cite{GR-long}). In Section \ref{Sec-near-opt-construction}, we discuss a way how asymptotically optimal/near optimal controls of the perturbed problem can be constructed on the basis of optimal/near optimal ACG families, and we illustrate the construction with a numerical example. In Section \ref{Sec-Proofs-I}, we give proofs of results stated in Sections \ref{Sec-Averaged system} and \ref{Sec-near-opt-construction}.


Let us conclude this section with some notations and definitions.
Given a compact metric space $W$, $\mathcal{B}(W)$ will stand for
the $\sigma$-algebra of its Borel subsets and $\mathcal{P}(W)$
will denote the set of probability measures defined on
$\mathcal{B}(W)$. The set $\mathcal{P}(W)$ will always be treated
as a compact metric space with a metric $\rho$,
 which is  consistent
 with its weak$^*$  topology. That is, a sequence
$\gamma^k \in \mathcal{P}(W), k =1,2,... ,$ converges to $\gamma
\in \mathcal{P}(W)$ in this metric if and only if
$$
 \lim_{k\rightarrow \infty}\int_{W} \phi(w) \gamma^k (dw) \ = \
 \int_{W} \phi(w) \gamma (dw),
$$
for any continuous $\phi(w): W \rightarrow \mathbb{R}^1$.

Using this metric $\rho$,
one can define the Hausdorff metric
$\rho_H$ on the set of subsets of $\mathcal{P}(W)$ as follows: $\forall \Gamma_i \subset \mathcal{P}(W) \ , \ i=1,2 \ ,$
\vspace{-0.15cm}
\begin{equation}\label{e:intro-3}
\rho_H(\Gamma_1, \Gamma_2) \stackrel{def}{=} \max \{\sup_{\gamma \in
\Gamma_1} \rho(\gamma,\Gamma_2), \sup_{\gamma \in \Gamma_2}
\rho(\gamma,\Gamma_1)\}, \
\end{equation}
where $\rho(\gamma, \Gamma_i) \BYDEF
\inf_{\gamma' \in \Gamma_i} \rho(\gamma,\gamma') \ .$
Note that, although, by some abuse of terminology,  we refer to
$\rho_H(\cdot,\cdot)$ as to a metric on the set of subsets of
${\cal P} (Y \times U)$, it is, in fact, a semi metric on this set
(note that $\rho_H(\Gamma_1, \Gamma_2)=0$ implies  $\Gamma_1
= \Gamma_2$ if  $\Gamma_1$ and $\Gamma_2$ are closed but the equality may not be true if at least one of these sets is not closed).

 Given a measurable function $w(\cdot): [0,\infty) \rightarrow W$, the {\it occupational measure} generated by
 this function on the interval $[0,S] $ is the probability measure $\gamma^{w(\cdot),S}\in \mathcal{P}(W)$ defined by the equation
 \begin{equation}\label{e:occup-S}
\gamma^{w(\cdot),S} (Q) \BYDEF \frac{1}{S} \int_0^{S}  1_Q(w(t))dt ,
 \ \ \forall Q \in \mathcal{B}(W),
\end{equation}
where $1_Q(\cdot) $ is the indicator function. The occupational measure  generated by
 this function on the interval $[0,\infty) $  is the probability measure $\gamma^{w(\cdot)}\in \mathcal{P}(W)$ defined as the limit (assumed to exist)
 \begin{equation}\label{e:occup-S-infty}
\gamma^{w(\cdot)} (Q) \BYDEF \lim_{S\rightarrow\infty}\frac{1}{S} \int_0^{S}  1_Q(w(t))dt ,
 \ \ \forall Q \in \mathcal{B}(W).
\end{equation}
Note that  (\ref{e:occup-S}) is equivalent to that
\begin{equation}\label{e:oms-0-1}
\int_{W} q(w) \gamma^{w(\cdot),S}(dw) = \frac{1}{S} \int _0 ^ S
q (w(t)) dt
\end{equation}
for any $q(\cdot)\in C(W)$, and
(\ref{e:occup-S-infty}) is equivalent to that
\begin{equation}\label{e:oms-0-1-infy}
\int_{W} q(w) \gamma^{w(\cdot)}(dw) =  \lim_{S\rightarrow\infty}\frac{1}{S} \int _0 ^ S q (w(t)) dt
\end{equation}
 for any $q(\cdot)\in C(W)$. Note that, if $w(\cdot)$ is $T$-periodic, then (\ref{e:oms-0-1-infy}) becomes
 \begin{equation}\label{e:oms-0-1-infy-per}
\int_{W} q(w) \gamma^{w(\cdot)}(dw) =  \frac{1}{T} \int _0 ^ T q (w(t)) dt.
\end{equation}

\section{Averaged system and averaged optimal control problem}\label{Sec-Averaged system}
It has been understood that, in dealing with uncontrolled systems with slow observables
($f(u,y)=f(y) $, $\ g(u,y)=g(y)$),
 the trajectories of the slow observables can be approximated by solutions of the averaged system (see, e.g., Theorem 6.5 in \cite{Art-Kevrekidis})
\begin{equation}\label{e:Art-2-1}
\frac{dz(\tau)}{d\tau}=  \epsilon\tilde h (\mu_{z(\tau)}), 
\end{equation}
where $\ \tilde h(\mu_z)\BYDEF \int_{ Y_z} F'(y)g(y)\mu_z(dy)\ $
and where $\mu_z(dy)$ is
 the occupational measure generated (induced) by  the solutions of the reduced system that are contained in $Y_z$.

In the presence of controls, there are multiple occupational measures that are generated by the controls and the corresponding
solutions of the reduced system staying in $Y_z$. The closed convex hull of the set of such measures
is presentable in the form (see Theorem 2.1 in \cite{Gai8})
\begin{equation}\label{e-CSO-7}
 W(z) \BYDEF \{\mu\in \mathcal{P}(U\times Y) \ : \ \int_{U\times Y_z}\nabla \phi(y)^Tf(u,y)\mu(du,dy)=0 \ \ \forall \phi(\cdot)\in C^1, \ \ \ \ \ supp (\mu)\subset U\times Y_z\}.\footnote{Given a compact set $X$  and $\mu\in \mathcal{P}(X)$, $supp (\mu)$ stands for the support of $\mu$, that is, the  \lq\lq minimal" closed subset of $X$ the $\mu$-measure of which is equal to $1$.    }
\end{equation}
 It is, therefore, natural to consider the following  averaged system
 \begin{equation}\label{e-CSO-17-11}
\frac{dz(\tau)}{d\tau} = \epsilon \tilde h(\mu(\tau)), \ \ \ \ \ z(0)=z_0,
\end{equation}
where
\begin{equation}\label{e-CSO-16}
\tilde h(\mu)\BYDEF \int_{U\times Y}F'(y)g(u,y)\mu(du,dy)
\end{equation}
and where
\begin{equation}\label{e-CSO-18-11}
\mu(\tau)\in W(z(\tau)).
\end{equation}
 In contrast to (\ref{e:Art-2-1}), the system (\ref{e-CSO-17-11}) is controlled, with the role of controls being played by the measure-valued functions $\mu(\tau) $ satisfying the inclusion (\ref{e-CSO-18-11}). Note that the map $W(z): Z\leadsto \mathcal{P}(U\times Y) $ defined in (\ref{e-CSO-7})  is convex and weak$^*$ compact-valued.

 By changing the time scale $t = \epsilon \tau $, one can rewrite (\ref{e-CSO-17-11}), (\ref{e-CSO-18-11}) in the form
 \begin{equation}\label{e-CSO-17}
\frac{dz(t)}{dt} = \tilde h(\mu(t)), \ \ \ \ \ z(0)=z_0,
\end{equation}
\vspace{-.2in}
\begin{equation}\label{e-CSO-18}
\mu(t)\in W(z(t)).
\end{equation}
In what follows, it is the system (\ref{e-CSO-17}) that will be referred to as the {\it averaged system}. A pair $(\mu(\cdot), z(\cdot)) $ that satisfies
(\ref{e-CSO-17}) and (\ref{e-CSO-18}) will be called a {\it solution of the averaged system}.

Consider the following optimal control problem
\begin{equation}\label{Vy-ave}
\displaystyle{
  \inf _{(\mu(\cdot),z(\cdot))  }\int_0 ^ { + \infty }e^{-C t}\tilde r(\mu(t))dt} \stackrel{\rm def}{=} \tilde R^*(y_0) ,
\end{equation}
where
\begin{equation}\label{Vy-ave-r}
 \tilde{r}(\mu) \BYDEF  \int_{U\times Y} r(u,y)
\mu(du,dy)
\end{equation}
and   $inf$ is sought over all solutions of the averaged system. This will be referred to as the
{\it averaged optimal control problem}.

 Note that in the shrinked time scale $t=\tau\epsilon$,  the problem
(\ref{Vy-perturbed-streched}) is rewritten in the form
\begin{equation}\label{Vy-perturbed}
\displaystyle{
   \inf _{u(\cdot)  }\int_0 ^ { + \infty }e^{-C t}r(u(t), y_{\epsilon}(t))dt} \stackrel{\rm def}{=} R^*(\epsilon, y_0),
\end{equation}
where $inf$ is over the controls and the corresponding solutions of the system
\begin{equation}\label{e-CSO}
\frac{dy_{\epsilon}(t)}{dt}= \frac{1}{\epsilon}f(u(t),y_{\epsilon}(t)) + g(u(t), y_{\epsilon}(t)), \ \ \ \ \ y(0)=y_0.
\end{equation}
This is the shrinked time scale  version  of the perturbed system (\ref{e-CSO-1}).
Note that in this time scale, the system (\ref{e-CSO-0}) describing the dynamics of the observables takes the form
\begin{equation}\label{e-CSO-0-2}
\frac{dz_{\epsilon}(t)}{dt}=  h(u(t), y_{\epsilon}(t)), \ \ \ \ \ z_{\epsilon}(0)=z_0.
\end{equation}

In this section,  we will present results establishing that  trajectories of the observables are approximated by solutions of the averaged system and that the optimal value of the  problem (\ref{Vy-perturbed}) is approximated by the optimal value of the averaged optimal control problem
(\ref{Vy-ave-r}). These results are valid under certain assumptions that are introduced below.

Let $Z$ be a compact subset of $\R^k $ such that $\ Y_z\subset Y $ for any $\ z\in Z$. In stating our results, we will be assuming (without mentioning it) that the multivalued map $Y_{z}:Z\leadsto Y$ is Lipschitz continuous on $Z$.
That is,
\begin{equation}\label{eq-Prop-Special-Cas-3-2}
d_H (Y_{z'}, Y_{z''})\leq L_0||z' - z''|| \ \ \ \ \forall z',z''\in Z,\ \ \ \ \          L_0= {\rm const}.
\end{equation}
(Here and in what follows, $d_H(\cdot,\cdot) $ stands for the Hausdorff metric defined on compact subsets of a finite-dimensional space.)
Note that the  Lipschitz continuity of $Y_{z}$  can be  verified to be true
 if there exists a $m\times (m-k)$ matrix  $D $ (which may depend on $y$) such that the inverse of  the matrix
$\hat D(y)\BYDEF (F'(y)^T, D) $
exists and is uniformly bounded: $ || \hat D(y)^{-1} ||\leq c ={\rm const}\ \forall y\in Y$ (see examples in Sections \ref{S-SC-2} and \ref{S-SC1}).

\begin{Assumption}\label{Continuity-W}
The multivalued map $W(z): Z\leadsto \mathcal{P}(U\times Y) $ is continuous and the multivalued map
$V(z): Z\leadsto \R^{k+1} $,
\begin{equation}\label{e-CSO-21}
V(z)\BYDEF \cup_{\mu\in W(z)}\{(\tilde h(\mu),\tilde r(\mu)\},
\end{equation}
is Lipschitz continuous.  That is,
 \begin{equation}\label{e-CSO-21-1}
d_H(V(z'),V(z''))\leq L ||z'-z''|| \ \ \ \ \forall z',z''\in Z,\ \ \ \ \          L= {\rm const}.
\end{equation}
\end{Assumption}
Sufficient conditions for Assumption \ref{Continuity-W} (and for Assumption \ref{Ass-existence-LOMS} introduced below) to be satisfied  are discussed in Sections \ref{S-SC-2} and \ref{S-SC1}.

The following proposition is important for understanding of the approximation results stated below.
\begin{Proposition}\label{Prop-inv-measure-set-2}
  If the multi-valued map $W(\cdot): Z\leadsto \mathcal{P}(U\times Y)$ is continuous, then the occupational measure $\mu^{u(\cdot),y(\cdot),S} $ generated
  on the interval $[0,S]$ by an arbitrary control $u(\cdot)$ and the solution $y(\cdot)$ of the system (\ref{e-CSO-2}) obtained with this control and an initial condition $y(0)\in Y_z$ satisfies the inequality
  \begin{equation}\label{e-CSO-11}
\rho(\mu^{u(\cdot),y(\cdot),S}, W(z))\leq \beta(S), \ \ \ \ \ \ \ \ \ \ \lim_{S\rightarrow\infty}\beta(S) = 0,
\end{equation}
the estimate being uniform with respect to $y(0)\in Y_z$ and $z\in Z$.
\end{Proposition}
{\it Proof.}
The proof of the proposition is given at the end of this section.
$\ \Box$

Denote by $\mathcal{M}_z (S, y(0))$ the set of occupational measures generated by all pairs $(u(\cdot), y(\cdot))$ on an interval $[0,S]$, where $u(\cdot)$ is a control and $y(\cdot)$ is the corresponding solution of
the reduced system (\ref{e-CSO-2}) such that $y(0)\in Y_z$. Then the estimate (\ref{e-CSO-11}) of Proposition \ref{Prop-inv-measure-set-2} can be rewritten in the form
\begin{equation}\label{e-CSO-further-1}
sup_{\mu\in \mathcal{M}_z (S, y(0))}\rho(\mu, W(z))\leq \beta(S) \ \ \ \ \ \ \forall y(0)\in Y_z, \ \ \ \ \ \forall z\in Z.
\end{equation}
\begin{Assumption}\label{Ass-existence-LOMS}
The following estimate is valid
\begin{equation}\label{e-CSO-further-2}
max_{\mu\in W(z))}\rho(\mu, \mathcal{M}_z (S, y(0)))\leq \beta(S) \ \ \ \ \ \forall y(0)\in Y_z, \ \ \ \ \ \forall z\in Z, \ \ \ \ \
\lim_{S\rightarrow\infty}\beta(S) = 0.
\end{equation}
\end{Assumption}
Note that, by the definition of the Hausdorff metric (see (\ref{e:intro-3})),   the validity of (\ref{e-CSO-further-1}) and (\ref{e-CSO-further-2})
is equivalent to
\begin{equation}\label{e-CSO-further-3}
\rho_H(\mathcal{M}_z (S, y(0)), W(z))\leq \beta(S) \ \ \ \ \ \forall y(0)\in Y_z, \ \  \forall z\in Z,
\end{equation}
where  $ \lim_{S\rightarrow\infty}\beta(S) = 0 $. The validity of  (\ref{e-CSO-further-3}) and, hence, the validity of Assumption \ref{Ass-existence-LOMS}  can be verified
with the help of the following statement (see instances of application of the latter in Sections \ref{S-SC-2} and \ref{S-SC1}).
\begin{Proposition}\label{Prop-Nec-Suf-LOMS}
The estimate (\ref{e-CSO-further-3}) is valid if and only if,
for any pair $(u^1(\cdot), y^1(\cdot))$, where
$u^1(\cdot)$ is an arbitrary control and $y^1(\cdot) $ is the solution of (\ref{e-CSO-2}) obtained with this control and $y^1(0)=y_0^1\in Y_z$ and for any  $y_0^2\in Y_z$, there exists a control
$u^2(\cdot)$ such that, for an arbitrary Lipschitz continuous function $q(u,y)$,
\begin{equation}\label{e-CSO-further-5}
|\frac{1}{S}\int_0^S q(u^1(\tau), y^1(\tau))d\tau - \frac{1}{S}\int_0^S q(u^2(\tau), y^2(\tau))d\tau |\leq \beta_q(S), \ \ \ \ \lim_{S\rightarrow\infty}\beta_q(S) = 0
\end{equation}
where
 $y^2(\cdot) $ is the solution of the system (\ref{e-CSO-2}) obtained with the control $u^2(\cdot) $ and with the initial condition $y^2(0)=y_0^2$.
\end{Proposition}
{\it Proof.}
The proof follows from Theorem 2.1(iii) and Proposition 4.1 in \cite{Gai8}.
$\ \Box$

The main results of this section are  the following two theorems.

\begin{Theorem}\label{Prop-Averaging-1}
Let Assumption  \ref{Continuity-W} be satisfied and let the averaged system  be viable in $Z$ (see \cite{aubin0}).
Then, corresponding to any pair  $(u_{\epsilon}(t), y_{\epsilon}(t))$, where $u_{\epsilon}(t)$ is a control and  $y_{\epsilon}(t)$ is the solution
of (\ref{e-CSO}) obtained with this control  such that
\begin{equation}\label{e-CSO-19}
y_{\epsilon}(t)\in Y, \ \ \ \ \ \ \ \ \ F(y_{\epsilon}(t))\in \check Z\subset int Z \ \ \ \ \ \forall t\in [0,\infty)
\end{equation}
($\check Z$ being a compact subset of the interior of  $Z$),
  there exists a solution $(\mu(t), z(t)) $ of the averaged system (\ref{e-CSO-17}), (\ref{e-CSO-18}) such that $\ z(t)\in Z \ \ \forall t\in [0,\infty)\ $ and such that
\begin{equation}\label{e-CSO-22}
\max_{t\in [0,T]}|| F(y_{\epsilon}(t)) - z(t)||\leq \beta(\epsilon, T),
\end{equation}
\begin{equation}\label{e:intro-0-3-8}
 |\int_0^T e^{-Ct}r(u_{\epsilon}(t),y_{\epsilon}(t)dt - \int_0^Te^{-Ct} \tilde r(\mu(t))dt| \leq \beta(\epsilon , T),
\end{equation}
where $ \ \lim_{\epsilon\rightarrow 0}\beta(\epsilon, T) = 0$.

\end{Theorem}

\medskip

\begin{Theorem}\label{Prop-converse-convergence}
Let Assumptions  \ref{Continuity-W} and \ref{Ass-existence-LOMS} be valid, and let the augmented perturbed system (\ref{e-CSO}), (\ref{e-CSO-0-2}) be viable in $Y\times Z$. Then, corresponding to any solution $(\mu(t), z(t))$ of the averaged system (\ref{e-CSO-17}), (\ref{e-CSO-18}) such that
\begin{equation}\label{e-CSO-20}
z(t) \in \check Z\subset int Z,\ \ \ \  \ \ \   Y_{z(t)}\subset \check Y\subset int Y\ \ \ \ \ \ \ \forall t\in [0,\infty)
\end{equation}
($\check Z $ and $\check Y $ being compact subsets of the interiors of $Y$ and $Z$ respectively),
there exists a control $u_{\epsilon}(t) $  such that the corresponding solution $y_{\epsilon}(t)$ of the system  (\ref{e-CSO}) satisfies the inclusions
\begin{equation}\label{e-CSO-20-101}
(y_{\epsilon}(t), F(y_{\epsilon}(t)))\in Y\times Z \ \ \ \forall \ t\in [0, \infty),
\end{equation}
and  the estimates (\ref{e-CSO-22}) and (\ref{e:intro-0-3-8}) are valid.
\end{Theorem}
{\it Proof.}
The proofs  of the theorems are given in Section \ref{Sec-Proofs-I}.
$\ \Box$

\begin{Corollary}\label{Corollary-inequality-values}
Let Assumption  \ref{Continuity-W} be satisfied and  the averaged system  be viable in $Z$. Let also there exist a control   $u_{\epsilon}(\tau) $ such that the corresponding solution $y_{\epsilon}(\tau) $ of the system (\ref{e-CSO})
 satisfies
 the inclusions (\ref{e-CSO-19}) and also
\begin{equation}\label{e-Cor-1-101}
\displaystyle{
   \liminf_{\epsilon\rightarrow 0}\int_0 ^ { + \infty }e^{-C t}r(u_{\epsilon}(t)(t), y_{\epsilon}(t))dt} = \liminf_{\epsilon\rightarrow 0} R^*(\epsilon, y_0).
\end{equation}
Then the optimal values of the problems (\ref{Vy-perturbed}) and (\ref{Vy-ave}) are related by the inequality
\begin{equation}\label{e-opt-val-ineq}
\liminf_{\epsilon\rightarrow 0}R^*(\epsilon, y_0)\geq \tilde R^*(y_0).
\end{equation}
\end{Corollary}
{\it Proof.}
Let us choose a  function $T(\epsilon)$ in such a way that
\begin{equation}\label{e:SP-aug-5-1-5-1-1}
\lim_{\epsilon\rightarrow 0}T(\epsilon) = \infty \ , \ \ \ \ \ \lim_{\epsilon\rightarrow 0}\beta(\epsilon , T(\epsilon)) = 0
\end{equation}
(the fact that such $T(\epsilon)$ exists is readily verifiable). By Theorem \ref{Prop-Averaging-1}, there exists a solution  $(\mu(t), z(t)) $ of the averaged system (\ref{e-CSO-17}) such that the estimate (\ref{e:intro-0-3-8}) is valid. Then
\begin{equation}\label{e:intro-0-3-8-new-1}
 |\int_0^{\infty}e^{-Ct} r(u_{\epsilon}(t),y_{\epsilon}(t))dt - \int_0^{\infty}e^{-Ct}\tilde r(\mu(t))dt| \leq \frac{2}{C}\max_{(u,y)\in U\times Y}|r(u,y)|e^{-CT(\epsilon)}+ \beta(\epsilon , T(\epsilon))\BYDEF \bar{\beta}(\epsilon),
\end{equation}
where  $\lim_{\epsilon\rightarrow 0}\bar{\beta}(\epsilon)=0$ (due to (\ref{e:SP-aug-5-1-5-1-1})).
Consequently,
$$
 \int_0^{\infty}e^{-Ct} r(u_{\epsilon}(t),y_{\epsilon}(t))dt\geq \int_0^{\infty}e^{-Ct}\tilde r(\mu(t))dt -\bar{\beta}(\epsilon)\ \geq\ \tilde R^*(y_0) - \bar{\beta}(\epsilon).
 $$
By (\ref{e-Cor-1-101}), the  latter implies (\ref{e-opt-val-ineq}).
$\ \Box$

\begin{Corollary}\label{Corollary-equailty-values}
Let Assumptions \ref{Continuity-W} and \ref{Ass-existence-LOMS} be satisfied.  Let the averaged system be viable in $Z$ and the augmented perturbed system be viable in $Y\times Z$. Let also  there exist a control   $u_{\epsilon}(\tau) $ such that the corresponding solution $y_{\epsilon}(\tau) $ of (\ref{e-CSO})  satisfies (\ref{e-CSO-19}),
(\ref{e-Cor-1-101}) and there exist an optimal solution $(\mu(t), z(t))$ of the averaged problem that satisfies (\ref{e-CSO-20}). Then the limit of the optimal value of the problem (\ref{Vy-perturbed}) exists and is equal to the optimal value of the averaged problem (\ref{Vy-ave}):
\begin{equation}\label{e-opt-val-ineq-1}
\lim_{\epsilon\rightarrow 0}R^*(\epsilon, y_0)=\tilde R^*(y_0).
\end{equation}
\end{Corollary}
{\it Proof.}
Similarly to the proof of Corollary \ref{Corollary-inequality-values}, it is established that from Theorem \ref{Prop-converse-convergence} it follows that
\begin{equation}\label{e-opt-val-ineq-2}
\limsup_{\epsilon\rightarrow 0}R^*(\epsilon, y_0)\leq \tilde R^*(y_0).
\end{equation}
This along with (\ref{e-opt-val-ineq}) prove (\ref{e-opt-val-ineq-1}).
$\ \Box$

\begin{Remark}
{\rm The statements of the Theorems \ref{Prop-Averaging-1} and \ref{Prop-converse-convergence} are simplified if one assumes  that $Z$ is invariant with respect to the solutions of the averaged system and that   $Y\times Z$ is invariant with respect to the solutions of the augmented perturbed system. In this case, corresponding to any pair  $(u_{\epsilon}(t), y_{\epsilon}(t))$, where $u_{\epsilon}(t)$ is an arbitrary control and $y_{\epsilon}(t)$ is the solution
of (\ref{e-CSO}) obtained with this control,  there exists a solution $(\mu(t), z(t)) $ of the averaged system such that (\ref{e-CSO-22}), (\ref{e:intro-0-3-8}) are satisfied (provided that Assumption \ref{Continuity-W} is valid). Conversely, corresponding to any solution $(\mu(t), z(t)) $  of the averaged system, there exists a control $u_{\epsilon}(t) $ such that (\ref{e-CSO-22}), (\ref{e:intro-0-3-8}) are satisfied
with $y_{\epsilon}(t)$ being the solution
of (\ref{e-CSO}) obtained with the control $u_{\epsilon}(t)$ (provided that Assumptions \ref{Continuity-W} and \ref{Ass-existence-LOMS} are valid). Also
(\ref{e-opt-val-ineq-1}) is valid  in this case as well.}
\end{Remark}

{\it Proof of Proposition \ref{Prop-inv-measure-set-2}.}
Assume the statement of the proposition is not true. Then there exist a number $\alpha > 0$ and sequences $S_i, \ \ z_i\in Z, \ \ u^i(\cdot), \ \ y^i(\cdot)$ \ ($y^i(0)
\in Y_{z_i}$),\ \ $i=1,2,..., $ such that
 \begin{equation}\label{e-CSO-12}
\rho(\mu^{u^i(\cdot),y^i(\cdot),S_i}, W(z_i))\geq \alpha, \ \ \ \ \ \ \ \ \ \ \lim_{i\rightarrow\infty}S_i = \infty .
\end{equation}
Due to compactness of $Z$ and due to compactness of $\mathcal{P}(U\times Y)$, one may assume (without loss of generality) that there exist limits
 \begin{equation}\label{e-CSO-13}
\lim_{i\rightarrow\infty}z_i=z\in Z, \ \ \ \ \ \ \ \ \ \ \ \ \ \ \ \ \   \lim_{i\rightarrow\infty}\mu^{u^i(\cdot),y^i(\cdot),S_i}=\mu\in \mathcal{P}(U\times Y).
\end{equation}
Also, $\lim_{i\rightarrow\infty}W(z_i)=W(z) $ (due to assumed continuity of $W(\cdot)$). Hence, from (\ref{e-CSO-12}) it follows that
\begin{equation}\label{e-CSO-14}
\rho(\mu, W(z))\geq \alpha .
\end{equation}
By the definition of the occupational measure (see (\ref{e:oms-0-1})), for any $\phi(\cdot)\in C^1$,
$$
\lim_{i\rightarrow\infty}\int_{U\times Y}\nabla \phi(y)^Tf(u,y)\mu^{u^i(\cdot),y^i(\cdot),S_i}(du,dy)= \lim_{i\rightarrow\infty}
\frac{1}{S_i}\int_0^{S_i}\phi(y^i(\tau)f(u^i(\tau),y^i(\tau))
$$
\vspace{-.2in}
$$
=  \lim_{i\rightarrow\infty}\frac{1}{S_i}[\phi(y^i(S_i))- \phi(y^i(0))] = 0.
$$
Consequently, due to the second equality in (\ref{e-CSO-13}),
\begin{equation}\label{e-CSO-15}
\int_{U\times Y}\nabla \phi(y)^Tf(u,y)\mu(du,dy)= 0.
\end{equation}
 From the fact that $Y_z$ is upper semicontinuous it follows that, for any $\kappa > 0 $,
 $\
 Y_{z_i}\subset Y_z + \kappa \bar B
 $
if $i$ is large enough ($\bar B$ being the closed unit ball in $\R^m) $. Consequently,
$$ 1 = \mu^{u^i(\cdot),y^i(\cdot),S_i}(U\times Y_{z_i})= \mu^{u^i(\cdot),y^i(\cdot),S_i}(U\times (Y_z + \kappa \bar B)), $$
and
$$
\mu (U\times (Y_z + \kappa \bar B)=  \lim_{i\rightarrow\infty}\mu^{u^i(\cdot),y^i(\cdot),S_i}(U\times (Y_z + \kappa \bar B)) = 1.
$$
Since $\kappa$ is arbitrary small, it follows that
$$
\mu (U\times Y_z) = 1 \ \ \Rightarrow \ \ supp(\mu)\subset U\times Y_z.
$$
The latter along with (\ref{e-CSO-15}) imply that $\mu\in W(z) $, which contradicts (\ref{e-CSO-14}).
$\ \Box$

\section{Verification of Assumptions \ref{Continuity-W} and \ref{Ass-existence-LOMS}}\label{S-SC-2}
In this section we will discuss some stability and controllability type conditions that ensure the validity of Assumptions \ref{Continuity-W} and \ref{Ass-existence-LOMS}.
\begin{Proposition}\label{Prop-1-stability}
Let there exist a positive definite matrix $Q$ and a positive number $\alpha$ such that
\begin{equation}\label{eq-Prop-Special-Cas-3-1}
\left( f(u,y') - f(u,y'') \right)^T Q (y'-y'')\leq - \alpha ||y'-y''||^2\ \ \ \ \ \forall u\in U, \ \ \forall y',y''\in Y_z, \ \ z\in Z.
\end{equation}
Then Assumptions \ref{Continuity-W} and \ref{Ass-existence-LOMS} are satisfied.

\end{Proposition}

{\it Proof.}
Let us first establish that  Assumption \ref{Ass-existence-LOMS} is satisfied. By Proposition \ref{Prop-Nec-Suf-LOMS}, to prove the validity of the latter, it is sufficient to show that, given a pair $(u^1(\cdot), y^1(\cdot))$, where
$u^1(\cdot)$ is an arbitrary control and $y^1(\cdot) $ is a solution of (\ref{e-CSO-2}) obtained with this control and  an initial condition $y^1(0)=y_0^1\in Y_z$ and given  $y_0^2\in Y_z$, there exists a control
$u^2(\cdot)$ such that, for an arbitrary Lipschitz continuous function $q(u,y)$, the estimate (\ref{e-CSO-further-5}) is valid (with $y^2(\cdot) $ being the solution of the system (\ref{e-CSO-2}) obtained with the control $u^2(\cdot) $ and with the initial condition $y^2(0)=y_0^2$).

Take $u^2(\cdot) $ to be equal to $u^1(\cdot)$. Then, by (\ref{eq-Prop-Special-Cas-3-1}),
$$
\frac{d }{d\tau}(y^1(\tau) - y^2(\tau))^T Q (y^1(\tau) - y^2(\tau))
= 2\left( f(u^1(\tau),y^1(\tau)) - f(u^1(\tau),y^2(\tau)) \right)^T Q (y^1(\tau)-y^2(\tau))
$$
\vspace{-.2in}
$$
\leq - 2\alpha ||y^1(\tau)-y^2(\tau)||^2\leq -\frac{2\alpha}{\lambda_{max}}(y^1(\tau) - y^2(\tau))^T Q (y^1(\tau) - y^2(\tau)),
$$
where $\lambda_{max} $ is the maximal eigenvalue of $Q$. It follows that
$$
(y^1(\tau) - y^2(\tau))^T Q (y^1(\tau) - y^2(\tau))\leq e^{- \frac{2\alpha}{\lambda_{max}}\tau}
(y_0^1 - y_0^2)^T Q (y_0^1 - y_0^2)\ \ \ \ \forall \ \tau\geq 0,
$$
the latter implying
\begin{equation}\label{eq-Prop-Special-Cas-3-3}
||y^1(\tau) - y^2(\tau)||\leq \beta_1 e^{-\beta_2\tau}||y_0^1 - y_0^2|| \ \ \ \forall \tau\geq 0,
\end{equation}
where $\beta_1 , \beta_2 $ are  constants ($\beta_1 > 1, \ \beta_2 > 0 $). Using (\ref{eq-Prop-Special-Cas-3-3}) and Lipschitz continuity of $q(\cdot)$, one obtains
$$
|\frac{1}{S}\int_0^S q(u^1(\tau), y^1(\tau))d\tau - \frac{1}{S}\int_0^S q(u^1(\tau), y^2(\tau))d\tau |
\leq L_q\frac{1}{S}\left( \int_0^S ||y^1(\tau) - y^2(\tau)|| d\tau\right)
$$
\vspace{-.2in}
$$
 \leq L_q\beta_1 ||y_0^1 - y_0^2||\left(\frac{1}{S}\int_0^S e^{-\beta_2\tau}d\tau \right)
\leq \frac{1}{S} \left(L_q\frac{\beta_1}{\beta_2}\max_{y',y''\in Y}||y'-y''||\right),
$$
where $L_q $ is a Lipschitz constant of $q(\cdot) $. Thus, (\ref{e-CSO-further-5}) is satisfied with $\beta_q(S) = O(\frac{1}{S}) $, and the validity of Assumption \ref{Ass-existence-LOMS} is established.

To prove the validity of Assumption \ref{Continuity-W}, it is sufficient to show that, for any $z', \ z'' \in Z$,
\begin{equation}\label{eq-Prop-Special-Cas-3-5}
\min_{\mu''\in W(z'')}|\int_{U\times Y_{z'}}q(u,y)\mu'(du,dy) - \int_{U\times Y_{z''}}q(u,y)\mu''(du,dy)|\leq \hat{\kappa}_q(||z'-z'' ||) \ \ \  \ \forall \ \mu'\in W(z'),
\end{equation}
where $\lim_{\theta\rightarrow 0}\hat{\kappa}_q(\theta)=0$ for  any continuous $q(u,y)$, and that
\begin{equation}\label{eq-Prop-Special-Cas-3-6}
\min_{\mu''\in W(z'')}|\int_{U\times Y_{z'}}q(u,y)\mu'(du,dy) - \int_{U\times Y_{z''}}q(u,y)\mu''(du,dy)|\leq \hat{L}_q||z'-z'' || \ \ \  \ \forall \ \mu'\in W(z'),
\end{equation}
where $\hat{L}_q$  is a positive constant if $q(u,y)$ satisfies Lipschitz conditions in $y$ (uniformly in $u\in U$).

Let $\mu'\in W(z') $. By (\ref{e-CSO-further-2}), there exist a control $u'(\cdot) $ and the corresponding solution $y'(\cdot) $ of the system (\ref{e-CSO-2}) satisfying an initial condition $y'(0)=y_0'\in Y_{z'} $ such that
\begin{equation}\label{eq-Prop-Special-Cas-3-7}
|\int_{U\times Y_{z'}}q(u,y)\mu'(du,dy) - \frac{1}{S}\int_0^S q(u'(\tau),y'(\tau))d\tau|\leq \beta_q(S) \ \ \  \ {\rm with} \ \  \ \lim_{S\rightarrow\infty}\beta_q(S) = 0.
\end{equation}
Let $y_0'' \BYDEF {\rm argmin}_{y\in Y_{z''}}\{||y- y_0'||\}$. By (\ref{eq-Prop-Special-Cas-3-2}),
\begin{equation}\label{eq-Prop-Special-Cas-3-8}
||y_0' - y_0''||\leq L_0 ||z'-z''||.
\end{equation}
Denote by $y''(\cdot)$ the solution of the system (\ref{e-CSO-2}) obtained with the control $u'(\cdot) $ that satisfies the initial condition
$y''(0) = y_0''$ (note that $y''(\tau)\in Y_{z''}\ \ \forall \tau \geq 0$). It will be shown below that
\begin{equation}\label{eq-Prop-Special-Cas-3-9}
||y'(\tau) - y''(\tau)||\leq \hat{L}||z'-z''|| \ \ \ \ \forall \ \ \tau\geq 0,
\end{equation}
where $\hat{L}$ is a positive constant. Using this inequality, one can obtain
\begin{equation}\label{eq-Prop-Special-Cas-3-10}
| \frac{1}{S}\int_0^S q(u'(\tau),y'(\tau))d\tau - \frac{1}{S}\int_0^S q(u'(\tau),y''(\tau))d\tau|\leq \kappa_q(\hat{L}||z'-z''||),
\end{equation}
where
$$ \kappa_q(\theta)\BYDEF \max_{u\in U}\max_{y',y''\in Y}\{|q(u,y')-q(u,y'')|\ | \ ||y'-y''||\leq \theta \}.$$ (Note that $\lim_{\theta\rightarrow 0}\kappa_q(\theta)=0$.)
  Denote by $\mu''_S \in \mathcal{M}_{z''} (S, y_0'')$, the occupational measure generated by the pair $(u'(\cdot), y''(\cdot)) $ on the interval $[0,S]$. From the  validity of (\ref{e-CSO-further-3}) it follows that there exists $\bar{\mu}''\in W(z'')$ such that
$$
\rho(\mu''_S, \bar{\mu}'')\leq \beta(S), \ \ \ \ \ \lim_{S\rightarrow\infty}\beta(S) = 0,
$$
which implies
\begin{equation}\label{eq-Prop-Special-Cas-3-11}
| \frac{1}{S}\int_0^S q(u'(\tau),y''(\tau))d\tau- \int_{U\times Y_{z''}}q(u,y)\bar{\mu}''(du,dy) |\leq
\bar{\beta}_q(S), \ \ \ \ \ \lim_{S\rightarrow\infty}\bar{\beta}_q(S) = 0.
\end{equation}
The latter along with (\ref{eq-Prop-Special-Cas-3-7}) and (\ref{eq-Prop-Special-Cas-3-10})  imply
$$
|\int_{U\times Y_{z'}}q(u,y)\mu'(du,dy) - \int_{U\times Y_{z''}}q(u,y)\bar{\mu}''(du,dy)|\leq \kappa_q(\hat{L}||z'-z'' ||) +\beta_q(S)
+ \bar{\beta}_q(S),
$$
which lead to
$$
\min_{\mu''\in W(z'')}|\int_{U\times Y_{z'}}q(u,y)\mu'(du,dy) - \int_{U\times Y_{z''}}q(u,y)\mu''(du,dy)|\leq \kappa_q(\hat{L}||z'-z'' ||) +\beta_q(S)
+ \bar{\beta}_q(S).
$$
Passing to the limit as $S\rightarrow\infty$ in the expression above, one obtains (\ref{eq-Prop-Special-Cas-3-5}) with $\hat{\kappa}_q(\theta) = \kappa_q(\hat{L}\theta) $. It also leads to
(\ref{eq-Prop-Special-Cas-3-6}) with $\hat{L}_q = \hat{L} L_q $ for a function $q(\cdot)$ that is Lipschitz continuous in $y$ with a constant $L_q$ (since $\kappa_q(\theta) \leq L_q\theta$ in this case).

To finalize the proof of the proposition, one now needs to establish the validity of
(\ref{eq-Prop-Special-Cas-3-9}). Let $\sigma \in (0,1) $ be such that
\begin{equation}\label{eq-Prop-Special-Cas-3-12}
\beta_1 e^{-\beta_2\sigma}= \frac{1}{2}
\end{equation}
(that is, $\sigma= \frac{1}{\beta_2}\ln 2\beta_1 $) and let $\tau_l\BYDEF l\sigma, \ l=0,1,... \ $.
Define $y_l'' $ by the equation $y_l'' \BYDEF {\rm argmin}_{y\in Y_{z''}}\{||y- y'(\tau_l)||\}$. Similarly to (\ref{eq-Prop-Special-Cas-3-8}), from (\ref{eq-Prop-Special-Cas-3-2}) it follows that
\begin{equation}\label{eq-Prop-Special-Cas-3-8-1}
||y'(\tau_l) - y_l''||\leq L_0 ||z'-z''||, \ \ l=0,1,... \ .
\end{equation}
Let $y_l''(\tau) $ stand for the solution of the system (\ref{e-CSO-2}) considered on the interval $[\tau_l, \tau_{l+1}] $ that is
obtained with  the control
$u'(\tau) $ and that satisfies the initial condition $y_l''(\tau_l)= y_l'' $ (note that $y_l''(\tau)\in Y_{z''} $). Using (\ref{eq-Prop-Special-Cas-3-8-1}) and Gronwall- Bellman lemma, one can readily verify that
\begin{equation}\label{eq-Prop-Special-Cas-3-8-2}
\max_{\tau\in [\tau_l,\tau_{l+1}]}||y'(\tau) - y_l''(\tau)||\leq  \hat{L}_0 ||z'-z''||, \ \ l=0,1,... \ ,
\end{equation}
where $\hat{L}_0 $ is an appropriately chosen constant. Also, similarly to (\ref{eq-Prop-Special-Cas-3-3}),
one can obtain that
$$
||y''(\tau) - y_l''(\tau)||\leq  \beta_1 e^{-\beta_2\tau}||y''(\tau_l) - y_l''||\ \ \ \ \ \forall \
\tau\in [\tau_l, \tau_{l+1}],
$$
which, by  (\ref{eq-Prop-Special-Cas-3-12}), leads to
$$
||y''(\tau_{l+1}) - y_l''(\tau_{l+1})||\leq  \beta_1 e^{-\beta_2\sigma}||y''(\tau_l) - y_l''||
\leq \frac{1}{2} ||y''(\tau_l) - y_l''||
$$
Consequently, one can use (\ref{eq-Prop-Special-Cas-3-8-2}) to obtain
$$
||y'(\tau_{l+1}) - y''(\tau_{l+1})||\leq ||y'(\tau_{l+1}) - y_l''(\tau_{l+1})   ||
+ || y_l''(\tau_{l+1}) - y''(\tau_{l+1})|| \leq \hat{L}_0 ||z'-z''|| + \frac{1}{2} ||y''(\tau_l) - y_l''||
$$
\vspace{-.2in}
$$
 \leq
\hat{L}_0 ||z'-z''|| + \frac{1}{2} [\ ||y''(\tau_l) - y'(\tau_l)|| + ||y'(\tau_l) - y_l''||\ ]\leq
\frac{3}{2}\hat{L}_0||z'-z''|| + \frac{1}{2} ||y'(\tau_l) - y''(\tau_l) ||, \ \ l=0,1,... \ .
$$
As can be readily verified, the estimates above imply (see also (\ref{eq-Prop-Special-Cas-3-8})) that,
for any $l=0,1,...$,
\begin{equation}\label{eq-Prop-Special-Cas-3-8-3}
||y'(\tau_{l}) - y''(\tau_{l})||\leq 3 \hat{L}_0||z'-z''|| + \frac{1}{2^l} ||y'_0 - y''_0 ||
\leq 3 \hat{L}_0||z'-z''|| + L_0 ||z'-z''|| \leq  (3 \hat{L}_0 + L_0) ||z'-z''||.
\end{equation}
By definition, for any $\tau\in [\tau_{l}, \tau_{l+1}] $,
$$
y'(\tau) = y'(\tau_{l}) + \int_{\tau_{l}}^\tau f(u'(s),y'(s))ds, \ \ \ \ \  \ \ y'(\tau) = y''(\tau_{l}) + \int_{\tau_{l}}^\tau f(u'(s),y''(s))ds.
$$
Hence, by (\ref{eq-Prop-Special-Cas-3-8-3}),
$$
||y'(\tau)- y''(\tau)||\leq (3 \hat{L}_0 + L_0) ||z'-z''|| + L_f \int_{\tau_{l}}^\tau
|| y'(s)- y''(s) ||ds\ \ \ \ \ \ \forall \tau \in [\tau_{l}, \tau_{l+1}],
$$
where $L_f$ is a Lipschitz constant of $f(u,y)$ (in $y$). Applying
Gronwall- Bellman lemma, one obtains
$$
\max_{\tau\in [\tau_l,\tau_{l+1}]}||y'(\tau) - y''(\tau)||\leq  \hat{L} ||z'-z''||, \ \ l=0,1,... \ ,
$$
where $\hat{L}$ is an appropriately chosen constant. This proves (\ref{eq-Prop-Special-Cas-3-9})
$\ \Box$

Note that  (\ref{eq-Prop-Special-Cas-3-1}) is a stability (dissipativity)   type condition. A similar condition was introduced in \cite{Dont} and was  used extensively in dealing with averaging of singularly perturbed control systems with decomposed fast and slow dynamics (see, e.g., \cite{Gai-Leiz}, \cite{GQ}, \cite{GR-long} and references therein). Emphasize also that (\ref{eq-Prop-Special-Cas-3-1})  is assumed to be satisfied for $y', y'' $ in the invariant set $Y_z$  $\ (\forall z\in Z$) and it is  not assumed to be true for all $y', y'' $ in $Y$.

Consider a special case in which
$ \ F(y) = V^T y$, with $rank(V)=k $. That is,
\begin{equation}\label{eq-Prop-Special-Cas-7}
Y_z = \{y \ | \ V^T y=z \}.
\end{equation}
Let us verify the inequality (\ref{eq-Prop-Special-Cas-3-2})  in this case. Let
$y'\in Y_{z'} $. That is, $\ V^T y'=z' $. Let $V_c $ be a $m\times (m-k)$ matrix such that the $m\times m$
matrix $\hat V\BYDEF (V, V_c) $ has rank $m$. Define $y''$ as the solution of the system
$$
 V^T y'' = z'', \ \ \ \ \ \ \ \ V_c^T y'' = V_c^T y'.
$$
Then $y''\in Y_{z''} $ and
$$
||y' - y''|| = || (\hat V^T)^{-1}\ ((z'-z'')^T, 0^T)^T||\leq ||(\hat V^T)^{-1}|| \ ||z'-z''||.
$$
This proves (\ref{eq-Prop-Special-Cas-3-2}) with $L_0=||(\hat V^T)^{-1}|| $. Note that the set $Y_z$ defined by the right-hand-side of (\ref{eq-Prop-Special-Cas-7}) is not compact. However, it can be shown that the solutions of the reduced system are contained in a compact set, say $\check Y $, for any $z\in Z$ (to show this, one may use an argument based on the validity of (\ref{eq-Prop-Special-Cas-3-3}) that is similar to one used in the proof of Theorem 3.1(ii) in \cite{Gai1})). In this case, for the results of Section \ref{Sec-Averaged system} to be valid,  $Y$ should be chosen large enough to contain $\check Y $ in its interior and $Y_z $ should be re-defined as the intersection of $Y$ and the set defined by the right-hand-side of (\ref{eq-Prop-Special-Cas-7}).

The condition (\ref{eq-Prop-Special-Cas-3-1}) is satisfied  if there exist a positive definite matrix $Q$ and a positive number $\alpha$ such that the Jacobian matrix $f'_y(u,y) $ satisfies the inequality
\begin{equation}\label{eq-Prop-Special-Cas-8}
v^T\left[\ f'_y(u,y)^T Q +Q f'_y(u,y)\ \right]v\leq - \alpha ||v||^2\ \ \ \   \forall u\in U, \ \ \ \forall y\in Y_z, \ \ \ \ \ \forall v\in \{v\ | \ V^Tv =0\} \ \ \ {\rm and} \ \ \  \forall z\in Z.
\end{equation}
To further illustrate the latter condition, let us consider an even more special case when
\begin{equation}\label{eq-Prop-Special-Cas-9}
f(u,y) = A(u)y,
\end{equation}
where $A(u) $ is a square matrix function of $u$ and the columns of $V $ are eigenvectors of the matrix $A^T(u) $ corresponding to its eigenvalue $0$ that has the multiplicity $k$. These  eigenvectors are assumed to be independent of $u$ as, for example, in the case when, for any $u\in U$, the matrix $A^T(u) $ is equal to the difference of a block-diagonal stochastic matrix and the identity matrix (see Remark \ref{Remark-SP-Markov} below).

Let $V_c $ be as above (that is, it is such that the rank of the matrix  $ \hat V = (V,V_c)\ $ is equal to $m$) and  let the $m\times (m-k) $ matrix $W_c$ be defined by the equation
$$
V^T W_c = 0, \ \ \ \ \ \ \ V_c^T W_c = I_c \ \ \ \ \ \ \ \Leftrightarrow \ \ \ \ \ \ \
W_c =  (\hat V^T)^{-1} (0, I_c)^T,
$$
where $ I_c$ is the $(m-k)\times (m-k) $ identity matrix. It can be shown (although we do not do it in the paper) that the inequality (\ref{eq-Prop-Special-Cas-8}) is satisfied if there exist a $(m-k)\times (m-k) $ positive definite matrix $Q_c $ and a positive number $\alpha_c $ such that
\begin{equation}\label{eq-Prop-Special-Cas-10}
\zeta^T [\ A_c^T(u) Q_c + Q_c A_c(u)\ ] \zeta \leq - \alpha_c || \zeta||^2 \ \ \  \forall u\in U \ \ \ \
{\rm and } \ \ \ \ \forall \zeta\in \R^{m-k},
\end{equation}
where $A_c(u) = V^T_c A(u) W_c $. Note that $Q_c $ and $\alpha_c $ satisfying (\ref{eq-Prop-Special-Cas-10}) exist if $A_c(u)$ is independent of $u$ and all its eigenvalues have negative real parts or if $A_c(u)$ is symmetric and  its eigenvalues are less than some negative constant for all $u\in U$.

\begin{Remark}\label{Remark-SP-Markov}
{\rm Systems of the form
\begin{equation}\label{e-Markov-Chains}
\frac{dy_{\epsilon}(\tau)}{d\tau}= \left[ A(u(\tau)) + \epsilon B(u(\tau))\right] y_{\epsilon}(\tau),
\end{equation}
(this being a special case of (\ref{e-CSO-1}), with $f(u,y) $ as in (\ref{eq-Prop-Special-Cas-9}) and $g(u,y) =B(u)y$) have been studied extensively in the context of singularly perturbed controlled Markov chains (see, e.g., \cite{Yin}
and references therein). In these studies, $A^T(u)$ was assumed to be
equal to the difference of a block-diagonal stochastic matrix and the identity matrix. Each block of $A^T(u)$ was assumed to be irreducible thus ensuring  that the zero eigenvalue of $A^T(u)$ has the multiplicity $k$ (provided that  the number of blocks is $k$) and that, corresponding to this eigenvalue, there are $k$ orthogonal eigenvectors  that consist of $1$s and $0 $s. In this case, the state variables are the probabilities of the states of the Markov chain and the observables represent the sums of the probabilities of the states that enter the blocks. Note that, with an appropriate change of variables, the system (\ref{e-Markov-Chains}) can be decomposed into fast and slow subsystems. However, such change of variables would make the nature of the problem much less transparent and it can also complicate the process of
finding a numerical solution; see, e.g., \cite{FAA02}, \cite{PG88} and  \cite{Yin}. }

\end{Remark}

The validity of Assumption  \ref{Ass-existence-LOMS} can also be verified under a controllability condition as stated in the proposition below.

\begin{Proposition}\label{Prop-1-controllability}
Let for any $z\in Z$ and any $y^1,y^2\in Y_z $, there exist a control $u(\cdot)$ that steers the system
(\ref{e-CSO-2}) from $y^1$ to $y^2$ in finite time $ T(y^1, y^2,z)\leq T_0 $ ($T_0$ being some positive constant). Then Assumption \ref{Ass-existence-LOMS} is satisfied.
\end{Proposition}
{\it Proof.}
Let $u^1(\cdot)$ be an arbitrary control and let $y^1(\cdot) $ be a solution of (\ref{e-CSO-2}) obtained with this control and  an initial condition $y^1(0)=y_0^1\in Y_z$. Also, let  $y_0^2\in Y_z$ and let $u_{1,2}(\tau)$ be the control that steers the system (\ref{e-CSO-2}) from the point $y_0^1$ to the point $y_0^2$ in time $T(y^1_0, y^2_0,z) $. Define the control $u^2(\tau) $ as being equal to $u_{1,2}(\tau)$ for
$\tau\in [0, T(y^1_0, y^2_0,z)]$ and as being equal to $u^1(\tau -T(y^1_0, y^2_0,z)) $ for $\tau > T(y^1_0, y^2_0,z)$. Denote by $y^2(\cdot) $ the solution of (\ref{e-CSO-2}) obtained with the control $u^2(\tau) $ and with the initial condition $y^2(0)= y^2_0 $. Note that, by construction, $y^2 (\tau) =y^1(\tau -T(y^1_0, y^2_0,z)) \ \forall \tau  > T(y^1_0, y^2_0,z) $. Hence, for any continuous $q(u,y)$, the following relationships are valid
$$
|\frac{1}{S}\int_0^S q(u^1(\tau), y^1(\tau))d\tau - \frac{1}{S}\int_0^S q(u^2(\tau), y^2(\tau))d\tau |
\leq \frac{1}{S}\int^S_{S-T(y^1_0, y^2_0,z)} |q(u^1(\tau), y^1(\tau))|d\tau
$$
$$
+ \frac{1}{S}\int_0^{T(y^1_0, y^2_0,z)} |q(u^2(\tau), y^2(\tau))|d\tau +  |\frac{1}{S}\int_0^{S-T(y^1_0, y^2_0,z)} q(u^1(\tau), y^1(\tau))d\tau - \frac{1}{S}\int_{T(y^1_0, y^2_0,z)}^S q(u^2(\tau), y^2(\tau))d\tau |
$$
$$
\leq \frac{2M_q}{S}, \ \ \ \ \ {\rm where} \ \ \ \ \ \ \ \ \ \ M_q\BYDEF\max_{(u,y)\in U\times Y}|q(u,y)|.
$$
The latter prove  (\ref{e-CSO-further-5}) with $\beta_q(S)= \frac{2M_q}{S}$, and, thus, the validity of Assumption \ref{Ass-existence-LOMS} follows from
Proposition \ref{Prop-Nec-Suf-LOMS}.

$\ \Box$

Let us give another sufficient condition for Assumption \ref{Continuity-W} to be satisfied provided that Assumption \ref{Ass-existence-LOMS} is valid
\begin{Proposition}\label{Prop-cont-Marc}
Let  Assumption \ref{Ass-existence-LOMS} be valid.  Then  Assumption \ref{Continuity-W} will be satisfied
 if
there exists a positive definite matrix $Q$  such that
\begin{equation}\label{eq-Prop-Special-Cas-3-1-1}
\left( f(u,y') - f(u,y'') \right)^T Q (y'-y'')\leq 0 \ \ \ \ \forall u\in U, \ \ \forall y',y''\in Y.
\end{equation}
\end{Proposition}

{\it Proof.}
The proof of the validity  of Assumption \ref{Continuity-W} follows exactly the same lines as in the proof of Proposition \ref {Prop-1-stability}, with the only difference being the way the estimate (\ref{eq-Prop-Special-Cas-3-9}) is established, where $y'(\tau), y''(\tau)$ are defined as in the aforementioned proof. To prove this estimate, note that, by
(\ref{eq-Prop-Special-Cas-3-1-1}),
$$
\frac{d }{d\tau}(y'(\tau) - y''(\tau))^T Q (y'(\tau) - y''(\tau))
= 2\left( f(u'(\tau),y'(\tau)) - f(u'(\tau),y''(\tau)) \right)^T Q (y'(\tau)-y''(\tau))\ \leq \ 0.
$$
Hence,
$$
(y'(\tau) - y''(\tau))^T Q (y'(\tau) - y''(\tau))\leq (y'_0 - y''_0)^T Q (y'_0 - y''_0)\ \ \ \ \forall \ \tau \geq 0,
$$
which leads to
\begin{equation}\label{e-CSO-further-5-extra}
\lambda_{min}||y'(\tau) - y''(\tau)||^2\leq \lambda_{max}||y'_0 - y''_0||^2\leq \lambda_{max}(L_0 ||z'-z''||)^2
\end{equation}
where the latter inequality follows from (\ref{eq-Prop-Special-Cas-3-8}) ($\lambda_{min} $ and $\lambda_{max}$ being the minimal and the maximum eigenvalues of $Q$).
The estimates (\ref{e-CSO-further-5-extra}) imply (\ref{eq-Prop-Special-Cas-3-9}) with $\hat L = L_0\sqrt{\frac{\lambda_{max}}{\lambda_{min}} }  $.
This completes the proof of the proposition.
$\ \Box$

Note that (\ref{eq-Prop-Special-Cas-3-1-1}) is  the \lq\lq non-expansivity" condition similar to one introduced in \cite{Marc-J}, and note that it is assumed to be satisfied for all vectors $y', y'' \in Y $ (and not only for ones
belonging to the same invariant set).
Let us conclude this section with an example.

{\it Example 1.} Assume that the perturbed system is of the form
\begin{equation}\label{ex-1-1}
\frac{dy_1(\tau)}{d\tau}= u(\tau)y_2(\tau) + \epsilon g_1(u(\tau), y(\tau)), \ \ \ \ \ \ \  \ \ \ \ \ \ \frac{dy_2(\tau)}{d\tau}= -u(\tau)y_1(\tau) + \epsilon g_2(u(\tau), y(\tau)),
\end{equation}
where $y(\tau)= (y_1(\tau), y_2(\tau))\in \R^2  $, \ $u(\tau)\in [-1,1]\in \R^1$ \ and \ $g_i(u, y): [-1,1]\times \R^2\rightarrow \R^1, \ i=1,2$ are continuous in $(u,y)$ and Lipschitz continuous in $y$ (uniformly in $u$). It can be readily verified that the reduced system (that is, the system obtained from (\ref{ex-1-1}) with $\epsilon = 0$) satisfies
(\ref{e-CSO-3}) with $F(y)=y_1^2 + y_2^2$.
 Hence,
\begin{equation}\label{ex-1-2}
Y_z = \{(y_1,y_2) \ | \ y_1^2 + y_2^2 = z\} \ \ \ \ \ \ \ \ \forall \ z\geq 0.
\end{equation}
Let $Z=[a,b]$, where $0<a<b$ are some constants. It is easy to check that (\ref{eq-Prop-Special-Cas-3-2}) is valid in this case. Indeed, for any $z',z''\in Z$ and any $y'\in Y_{z'} $ (note that $ ||y'|| = \sqrt{z'} $), take $y''= \frac{y' \sqrt{z''}}{\sqrt{z'}}$. Then $||y''||=\sqrt{z''}$  and  $y'' \in Y_{z''}  $. Also,
$$
||y'-y''|| = \left|\sqrt{z'}- \sqrt{z''}\right|\leq  \frac{\left|z' - z''\right|}{2\min\{\sqrt{z'}, \sqrt{z''} \}}\leq \frac{\left|z' - z''\right|}{2\sqrt{a}}.
 $$
This proves (\ref{eq-Prop-Special-Cas-3-2}) with $L_0= \frac{1}{2\sqrt{a}}$. It can also be readily seen that, with the use of the control $u(\tau)=1  $, a solution of the reduced system can reach any point in the set $Y_z$ from any other point of this set within a time interval that is less or equal than $2\pi $.
 Thus, the controllability condition of Proposition \ref{Prop-1-controllability} is satisfied and, hence, Assumption \ref{Ass-existence-LOMS} is valid. One can also verify that the inequality (\ref{Ass-existence-LOMS}) is satisfied with $Q=I$. Indeed,
$$
\left( f(u,y') - f(u,y'') \right)^T (y'-y'')= u\ [(y_2'- y_2''), -(y_1'- y_1'')]\ [(y_1'- y_1''), (y_2'- y_2'')]^T = 0 \ \ \ \forall y',y''.
$$
Thus, by Proposition \ref{Prop-cont-Marc}, Assumption \ref{Continuity-W} is valid as well.

\section{A  case of uncontrolled reduced system}\label{S-SC1}

Assume that $f(u,y)=f(y) $. That is, the reduced system (\ref{e-CSO-2})  has the form
\begin{equation}\label{e-CSO-2-SC-I}
\frac{dy(\tau)}{d\tau}= f( y(\tau))
\end{equation}
and is uncontrolled.
Let us introduce the following assumption about the system (\ref{e-CSO-2-SC-I}).
\begin{Assumption}\label{def-1}
For any $z\in Z$, the system (\ref{e-CSO-2-SC-I})  has a $T_z$-periodic solution $y_z(\tau)\in Y_z$ (with $T_z$ being uniformly bounded for $z\in Z$)
and  any other solution of this system with the initial condition  in $Y_z$ is a time shift of $y_z(\tau) $.
\end{Assumption}

If Assumption \ref{def-1} is satisfied, then, for any continuous function $q(y)$ and for any solution $y_{y_0}(\tau)$ of the system (\ref{e-CSO-2-SC-I}) satisfying the initial condition $\ y(0)=y_0\in Y_z $,
 \begin{equation}\label{e-CSO-5-11-SC-1}
 \lim_{S\rightarrow\infty}\int_0^Sq(y_{y_0}(\tau))d\tau=\frac{1}{T_z}\int_0^{T_z}q(y_z(\tau))d\tau,
\end{equation}
with the convergence being uniform with respect to $y_0\in Y_z$ and with respect to $z\in Z$. It implies that
\begin{equation}\label{e-CSO-5}
\lim_{S\rightarrow\infty}\rho(\nu^{y_{y_0}(\cdot),S}, \nu_z) = 0,
\end{equation}
where $\nu^{y_{y_0}(\cdot),S}\in\mathcal{P}(Y_z) $ is the occupational measure generated by $y_{y_0}(\tau)$ on the interval $[0,S]$ and $\nu_z\in\mathcal{P}(Y_z) $ is the occupational measure generated by the periodic solution $y_z(\tau) $, the latter being defined by the equation
\begin{equation}\label{e-CSO-5-11}
 \int_{Y_z}q(y)\nu_z(dy)=\frac{1}{T_z}\int_0^{T_z}q(y_z(\tau))d\tau \ \ \ \ \forall q(\cdot)\in C
\end{equation}
 Note  that $\nu_z$ is the invariant measure of the system (\ref{e-CSO-2-SC-I}) that is supported on $Y_z$ (see, e.g., \cite{Art-Kevrekidis}). Note also that from the fact that the convergence in (\ref{e-CSO-5-11-SC-1}) is uniform it follows that the convergence in (\ref{e-CSO-5}) is uniform.
That is,
\begin{equation}\label{e-CSO-6}
\rho(\nu^{y_{y_0}(\cdot),S}, \nu_z)\leq \beta(S) \ \ \ \forall \ y_0\in Y_z, \ \forall \ z\in Z ,\ \ \ \ \
\lim_{S\rightarrow\infty} \beta(S)= 0.
 \end{equation}

 \begin{Proposition}\label{Prop-simple}
Let the reduced system has the form (\ref{e-CSO-2-SC-I}) and let Assumption \ref{def-1} be satisfied . Then:

(i) The set $W(z)$  is presentable  in  the form
  \begin{equation}\label{e-CSO-8}
W(z)= \{\mu\in \mathcal{P}(U\times Y) \ : \ \mu(du,dy)=\sigma(du|y)\nu_z(dy), \ \ \ \sigma(du|\cdot): Y\to \mathcal{P}(U)\},
\end{equation}
 where $\nu_z$ is defined by (\ref{e-CSO-5-11}) and $\ \sigma(du|\cdot): Y\to \mathcal{P}(U) $ are   Borel measurable.

 (ii)  
  Assumption
\ref{Ass-existence-LOMS} is satisfied.

 (iii) If, in addition,
\begin{equation}\label{e-CSO-2-SC-I-101}
\frac{|T_{z'}- T_{z''}|}{\max\{T_{z'}, T_{z''}\}}\leq \tilde L||z'-z''|| \ \ \ \forall \ z', z'' \in Z \ \  \ \ \ \ {\rm for\ \ some} \ \ \tilde L={\rm const},
\end{equation}
  then Assumption \ref{Continuity-W} is satisfied as well.

\end{Proposition}

{\it Proof of the statements (i) and (ii).}
 By Theorem 2.1(iii) in \cite{Gai8}, from (\ref{e-CSO-6}) it follows that the set
\begin{equation}\label{e-CSO-6-1-0}
\mathcal{V}(z)\BYDEF \{\nu\in \mathcal{P}( Y)\ : \ \int_{ Y_z}\nabla \phi(y)^Tf(y)\nu(dy)=0 \ \ \forall \phi(\cdot)\in C^1, \ \ \ \ \ supp(\nu)\subset Y_z\}
\end{equation}
is a singleton and that
\begin{equation}\label{e-CSO-6-1}
\mathcal{V}(z)=\{\nu_z\}.
\end{equation}
Note that an arbitrary measure $\mu\in \mathcal{P}(U\times Y) $
can be “disintegrated” as follows
\begin{equation}\label{e-CSO-6-1-1}
\mu(du, dy) = \sigma(du|y)\nu(dy),
\end{equation}
where $\ \nu(dy)= \mu(U, dy) $ and $\sigma(du|y)\in \mathcal{P}(U) $ is such that the integral $\ \int_U q(u,y)\sigma(du|y) $ is a Borel measurable function on $Y$ for any continuous $q(u, y)$, and
\begin{equation}\label{e-CSO-6-1-2}
\int_{U\times Y}q(u,y)\mu(du, dy)= \int_{y\in Y} \left(\int_Uq(u,y)\sigma(du|y)\right)\nu(dy).
\end{equation}
Note that the Borel measurability of $\ \int_U q(u,y)\sigma(du|y) $ for any continuous $q(u,y)$ implies the Borel measurability of
$\ \sigma(du|\cdot): Y\to \mathcal{P}(U) $.

{\it Proof of (i)}. Assume that
 $\mu$ belongs to the right-hand-side of (\ref{e-CSO-8}). That is,
 $\mu(du,dy) =\sigma(du|y)\nu_z(dy)  $. Then $\mu\in W(z) $
(due to (\ref{e-CSO-6-1-0}) and (\ref{e-CSO-6-1})).
 Conversely,  if $\mu\in W(z) $, then using the disintegration (\ref{e-CSO-6-1-1}), one can verify that $\nu\in \mathcal{V}(z) $.
 Since the latter consists only of one element which is $\nu_z $, this proves that $\mu$ belongs to the right hand side of (\ref{e-CSO-8}). Thus, the validity  of (\ref{e-CSO-8}) is established and (i) is proved.

 {\it Proof of (ii)}. The statement (ii) follows from Proposition \ref{Prop-1-controllability} since
  any point in the set $Y_z$ can be reached by a solution of (\ref{e-CSO-2-SC-I}) from any other point of this set within a time interval that is less or equal than $T_z$.
 %
%
$\ \Box$

The proof of the statement (iii) of Proposition \ref{Prop-simple} is based on the following lemma.

\begin{Lemma}\label{Lemma-Cont-Uncontrolled-Case}
Let, the integral $\int_{Y_z}q(y)\nu_z(dy) $ be continuous in $z$ for any
 continuous  $q(y)$ and let this integral be Lipschitz continuous in $z$ with a Lipschitz  constant that can be chosen the same for all Lipschitz continuous  $q(y)$ that have a common Lipschitz constant. Then  Assumption \ref{Continuity-W} is satisfied.
\end{Lemma}
{\it Proof.}
Assume that the integral $\ \frac{1}{T_z}\int_0^{T_z}q(y_z(\tau))d\tau\ $ is continuous for any continuous $q(\cdot)$.
Due to the fact that the map $W(\cdot) $ is convex and compact valued, to prove the continuity of $W(\cdot)$, it is sufficient to show that, if $z_i\to z$ ($z_i\in Z$), then, for any continuous $q(u,y)$,
\begin{equation}\label{e-CSO-9}
\lim_{z_i\rightarrow z}\min_{\mu\in W(z_i)} \int_{U\times Y_{z_i}}q(u,y)\mu(du,dy) = \min_{\mu\in W(z)} \int_{U\times Y_{z}}q(u,y)\mu(du,dy).
\end{equation}
By (\ref{e-CSO-8}),
$$
\min_{\mu\in W(z_i)} \int_{U\times Y_{z_i}}q(u,y)\mu(du,dy) = \int_{ Y_{z_i}}q^*(y)\nu_{z_i}(dy)
$$
$$
 \min_{\mu\in W(z)} \int_{U\times Y_{z}}q(u,y)\mu(du,dy)= \int_{ Y_{z}}q^*(y)\nu_z(dy)
$$
where $\ q^*(y)\BYDEF \min_{u\in U}q(u,y)$. The function $q^*(y)$ is continuous. Hence,
$$
\lim_{z_i\rightarrow z}\int_{ Y_{z_i}}q^*(y)\nu_{z_i}(dy) = \int_{Y_z}q^*(y)\nu_z(dy).
$$
The latter proves (\ref{e-CSO-9}) and, thus,  the continuity of $W(z)$ is established.

 Let $L_1$ be a Lipschitz constant  of the vector function $X(u,y)\BYDEF(h(u,y),r(u,y))$ in $y$. That is,
\begin{equation}\label{e-SC-2-0}
||X(u,y') - X(u,y'')||\leq L_1||y'-y''||.
\end{equation}
Assume  that the integral $\ \frac{1}{T_z}\int_0^{T_z}q(y_z(\tau))d\tau\ $ is Lipschits continuous with a constant $L$ for any Lipschitz continuous $q(\cdot)$ that have a Lipschitz constant $L_1$.

The map $V(z): Z\leadsto \R^{k+1} $ defined in (\ref{e-CSO-21}) is convex and compact valued. Hence to prove (\ref{e-CSO-21-1}),
it is sufficient to prove that, for any $v\in \R^{k+1}$ such that  $ \ ||v||=1$,
\begin{equation}\label{e-SC-2}
|\min_{\mu\in W(z')} \int_{U\times Y_{z'}}v^TX(u,y)(u,y)\mu(du,dy) - \min_{\mu\in W(z'')} \int_{U\times Y_{z''}}v^TX(u,y)(u,y)\mu(du,dy)\ |
\leq L ||z'-z''||
\end{equation}
(see, e.g., Lemma P.2.9, p. 207 in \cite{Gai0}). By (\ref{e-CSO-8}),
$$
\min_{\mu\in W(z')} \int_{U\times Y_{z'}}v^TX(u,y)\mu(du,dy)=
\int_{Y_{z'}}q_v^*(y)\nu_{z'}(dy),
$$
$$
\min_{\mu\in W(z'')} \int_{U\times Y_{z''}}v^TX(u,y)\mu(du,dy)=
\int_{Y_{z''}}q_v^*(y)\nu_{z''}(dy),
$$
where $q_v^*(y)\BYDEF \min_{u\in U}\{v^TX(u,y)\}$. Note that from (\ref{e-SC-2-0}) it follows that
$\ |v^TX(u,y')-v^TX(u,y'')|\leq L_1||y'-y''|| \ $, which, in turn, implies $|q_v^*(y') - q_v^*(y'')|\leq L_1||y'-y''|| $.
Thus
$$\ \left|\int_{Y_{z'}}q_v^*(y)\nu_{z'}(dy)- \int_{Y_{z''}}q_v^*(y)\nu_{z''}(dy)\right|\leq L ||z'-z''||. $$
This proves (\ref{e-SC-2}), and, thus, the proof of the lemma is completed.
$\ \Box$

{\it Proof of the statement (iii) of Proposition \ref{Prop-simple}}. To prove the statement, it is sufficient to verify that the conditions of Lemma \ref{Lemma-Cont-Uncontrolled-Case} are satisfied.

Let $z', z''\in Z $. Let $y_{z'}(\tau) $ be a solution of (\ref{e-CSO-2-SC-I}) such that  $y_{z'}(0)\in Y_{z'} $ and let
$y_{z''}(\tau) $ be the solution of (\ref{e-CSO-2-SC-I}) that satisfies the initial condition
$y_{z''}(0)={\rm argmin}_{y\in Y_{z''}}\{||y-y_{z'}(0)|| \} $. From (\ref{eq-Prop-Special-Cas-3-2}) it follows that
$||y_{z''}(0)-y_{z'}(0)||\leq  L_0||z''-z'' ||   $. Using this estimate and applying Gronwall-Bellman Lemma,
one can deduce that
\begin{equation}\label{e-SC-2-21}
||y_{z''}(\tau)-y_{z'}(\tau)||\leq  \hat L_0||z''-z'' || \ \ \ \ \ \ \ \ \forall \tau\in [0, \max\{T_{z'}, T_{z''}\}],
\end{equation}
where $\hat L_0$ is an appropriately chosen constant.

Assume (without loss of generality) that $T_{z''}\geq T_{z'} $.  For an arbitrary continuous function
$q(y) $, one can write down the following inequalities
$$
\left|\frac{1}{T_{z'}}\int_0^{T_{z'}}q(y_{z'}(\tau))d\tau\ - \ \frac{1}{T_{z''}}\int_0^{T_{z''}}q(y_{z''}(\tau))d\tau\right|\ \leq \ \frac{1}{T_{z'}}\int_0^{T_{z'}}\left|q(y_{z'}(\tau))- q(y_{z''}(\tau))\right| d\tau
$$
\vspace{-.2in}
$$
+ \ \left|\frac{1}{T_{z''}} -  \frac{1}{T_{z'}}  \right| \int_0^{T_{z'}}\left| q(y_{z''}(\tau))\right| d\tau
\ + \ \frac{1}{T_{z''}}\int_{T_{z'}}^{T_{z''}}|q(y_{z''}(\tau))|d\tau\ \leq \
\kappa_q (\hat L_0||z''-z'' ||) + 2 M_q \frac{|{T_{z'}}-T_{z''}|}{T_{z''}}
$$
\vspace{-.2in}
$$
\leq \ \kappa_q (\hat L_0||z''-z'' ||) + 2 M_q \tilde L||z'-z''||
$$
where (\ref{e-SC-2-21}) and, subsequently, (\ref{e-CSO-2-SC-I-101}) have been taken into account and  where $\kappa_q(\cdot) $ is the modulus of continuity of $q(\cdot)$ on $Y$ ($\lim_{\theta\rightarrow 0}\kappa_q(\theta)=0$) and $M_q = \max_{y\in Y}|q(y)| $.
By (\ref{e-CSO-5-11}), the latter imply
\begin{equation}\label{e-SC-2-21-1}
\left|\int_{Y_{z'}}q(y)\nu_{z'}(dy) - \int_{Y_{z''}}q(y)\nu_{z''}(dy) \right|\ \leq\
\kappa_q (\hat L_0||z''-z'' ||) + 2 M_q \tilde L||z'-z''||,
\end{equation}
which proves that the integral $\int_{Y_{z}}q(y)\nu_{z}(dy) $ is continuous in $z$. Also, it proves that
that this integral is Lipschitz continuous with the constant $L_1\hat L_0 + 2 M_q \tilde L  $ for any
Lipschitz continuous $q(\cdot) $ that has a Lipschitz constant $L_1 $ (since
$\ \kappa_q (\hat L_0||z''-z'' ||)\leq L_1\hat L_0 ||z''-z'' ||\ $ in this case). This completes the proof of the statement.
$\ \Box$

Let us consider two examples in which the conditions of Proposition \ref{Prop-simple} are satisfied.

{\it Example 1 (continued).} Consider   system (\ref{ex-1-1}) with $u(\tau)\equiv 1 $ (an uncontrolled version of the system). As one can see, this system satisfies Assumption \ref{def-1} with $T_z =2\pi$ for all $z>0$. We have already checked that (\ref{eq-Prop-Special-Cas-3-2}) is satisfied on $Z=[a,b]$ ($0<a<b$). Also, one can see that (\ref{e-CSO-2-SC-I-101}) is valid with $\tilde L=0$.

{\it Example 2 (perturbed Lotka-Volterra equations)}.  Consider the system
\begin{equation}\label{e:ex-1}
\frac{dy_1(\tau)}{d\tau}= -y_1(\tau)+y_1(\tau)y_2(\tau)-\epsilon u(\tau)y_1(\tau) , \ \ \ \ \ \ \  \ \ \ \ \ \ \frac{dy_2(\tau)}{d\tau}= y_2(\tau)-y_1(\tau)y_2(\tau),
\end{equation}
where $y=(y_1,y_2)\in \R^2$ and controls $u(\cdot)$ are scalar functions of time satisfying the inclusion $\ u(\tau)\in [0,1]$. A  system similar to (\ref{e:ex-1}) but with no controls and with the presence of  small random perturbations was considered in \cite{Fima}.
The reduced system (\ref{e-CSO-2-SC-I}) is the Lotka-Volterra system
\begin{equation}\label{e:ex-1-LV}
\frac{dy_1(\tau)}{d\tau}= -y_1(\tau)+y_1(\tau)y_2(\tau), \ \ \ \ \ \ \  \ \ \ \ \ \ \frac{dy_2(\tau)}{d\tau}= y_2(\tau)-y_1(\tau)y_2(\tau),
\end{equation}
and it is known that (\ref{e-CSO-3}) is satisfied with
\begin{equation}\label{e:ex-1-LV-101}
F(y)=\ln y_2(t)-y_2(t)+\ln y_1(t)-y_1(t).
\end{equation}
Note that it can be shown that the system (\ref{e:ex-1-LV}) does not satisfy the non-expansivity condition (\ref{eq-Prop-Special-Cas-3-1-1}).

For any $z<-2$, the set $Y_z\BYDEF \{y\in \R^2 \ : \ F(y)=z\}$ is a closed curve, the solutions of the  reduced system being periodic and moving along this curve. Two  such curves (for $z=-3$ and $z=-2.05$) are depicted in Fig. 1 in Section \ref{Sec-near-opt-construction} below. For $z=-2$, the curve  is degenerated into  the equilibrium  point $(1,1)$. Note that   $T_{z''}< T_{z'}$ if $z''>z' $, where $T_z$ is  the
period of the solutions of (\ref{e:ex-1-LV}) that are  contained in $Y_z$ (see, e.g., Theorem 2 in \cite{Sze}).

Let
\begin{equation}\label{e:ex-1-LV-101-1}
Z=[-3, -2.05], \ \ \ \ \ \ Y=[a_1, b_1]\times [a_2, b_2],
 \end{equation}
 where $0< a_i< b_i, \ i=1,2, $ are  such that $Y_z\subset Y \ \forall z\in [-3, -2.05] $.  Assumption \ref{def-1} is obviously satisfied in this case, and it can  be readily verified that (\ref{eq-Prop-Special-Cas-3-2}) is valid (this follows from  the fact that $\min_{y\in Y_z}||F'(y)||\geq \alpha >0 \ \ \forall z\in Z$). Since  $T_z$ is a monotone decreasing, $T_{z}\geq T_{-2.05}>0\ \ \forall z\in Z$. Hence, to verify (\ref{e-CSO-2-SC-I-101}), it is sufficient to establish that $T_z$ is Lipschitz continuous on $Z$. The proof of this is, however,  quite technical, and we do not present it in the paper.

In Section \ref{Sec-near-opt-construction}, we will return to  this example in the context of numerical construction of asymptotically near optimal solution of an optimal control problem considered on the trajectories of the perturbed system (\ref{e:ex-1}). In the remainder of this section we will establish statements similar to (i) and (ii) of Proposition \ref{Prop-simple} under a weaker assumption.


\begin{Assumption}\label{Ass-almost-per}
There exists a solution $\bar y(t)\in Y_z$ of the reduced system (\ref{e-CSO-2-SC-I}) such that the following two conditions are satisfied:

(i) The solution $\bar y(\cdot)$  generates an occupational measure $\nu_z$. That is, for any  continuous function $q(y)$,
\begin{equation}\label{e-CSO-9-new-0-1}
\left|\frac{1}{S}\int_0^S q(\bar y(t))dt - \int_{Y_z}q(y)\nu_z(dy) \right| \leq \kappa_q(S), \ \ \ \ {\rm where} \ \ \ \ \lim_{S\rightarrow\infty}\kappa_q(S) = 0.
\end{equation}

(ii) There exist positive functions $\alpha_1(S)$ and $\alpha_2(S)$,
\begin{equation}\label{e-CSO-9-new-0}
\lim_{S\rightarrow\infty}\alpha_i(S) = 0, \ \ \ \ \ i=1,2,
\end{equation}
such that any  other solution $y(t)\in Y_z$ of the reduced system (\ref{e-CSO-2}) satisfies the inequality
\begin{equation}\label{e-CSO-9-new}
\max_{\tau\in [0,S]}||y(\tau+\theta_S) - \bar y(\tau)|| \leq \alpha_1 (S),
\end{equation}
for some $\theta_S $ such that
\begin{equation}\label{e-CSO-9-new-1}
\frac{\theta_S}{S}\in [0,\alpha_2(S)].
\end{equation}
\end{Assumption}

Note that Assumption \ref{Ass-almost-per} is satisfied if Assumption \ref{def-1} is satisfied,
 in which case one may take $\alpha_1(S)= 0$ and $\alpha_2(S)=\frac{T_z}{S} $,  where $T_z$ is the period of the solutions staying in $Y_z$.

\begin{Proposition}\label{Prop-Special-Case}
Let  Assumption \ref{Ass-almost-per} be satisfied. Then  the statements (i) and (ii)  of Proposition \ref{Prop-simple} are valid.
\end{Proposition}
{\it Proof.}
From Assumption \ref{Ass-almost-per}(ii) it follows that, for any Lipschitz continuous function $q(y)$ and  any solution $y(t)$ of the system
(\ref{e-CSO-2-SC-I}),
$$
\left|\frac{1}{S}\int_0^S q(y(\tau))d\tau - \frac{1}{S}\int_0^S q(\bar y(\tau))d\tau \right|\leq
\frac{1}{S}\int_0^{S-\theta_S} |q(y(\tau+\theta_S))dt -  q(\bar y(\tau))|d\tau + \frac{2M_q\theta_s}{S}
$$
\vspace{-.2in}
\begin{equation}\label{eq-Prop-Special-Cas-1}
\leq L_q \alpha_1(S) + \frac{2M_q\theta_s}{S} \ \leq\ L_q \alpha_1(S) + 2M_q \alpha_1(S)\BYDEF \alpha_q(S),
\end{equation}
where $L_q$ is a Lipschitz constant of $q(\cdot) $ and $\ M_q\BYDEF \max_{y\in Y}|q(y)|$.
Then
$$
\left|\frac{1}{S} \int_0^S q(y(\tau))d\tau\ - \
\int_{Y_z}q(y)\nu_z(dy) \right|\leq \alpha_q(S) + \kappa_q(S),
$$
where $\kappa_q(S) $ is as in (\ref{e-CSO-9-new-0-1}).
The latter implies that the occupational measure generated by any solution staying in $Y_z$
 converges to $\nu_z$. That is, (\ref{e-CSO-6}) is valid. From this point, the proofs of (i)  follows exactly the same lines as those of Proposition \ref{Prop-simple}(i).

 To prove (ii),
 let us  show that given an arbitrary pair $(u^1(\cdot), y^1(\cdot))$, where $y^1(\cdot) $ is a solution of (\ref{e-CSO-2-SC-I}) with $y^1(0)\in Y_z$ and given another solution $y^2(\cdot) $ of (\ref{e-CSO-2-SC-I}) with $y^2(0)\in Y_z$, there exists a control
$u^2(\cdot)$ such that, for any Lipschitz continuous function $q(u,y)$, the estimate (\ref{e-CSO-further-5}) is valid (see Proposition \ref{Prop-Nec-Suf-LOMS}). To this end, define, first, $\bar u(\tau) $ to be equal to
 $u^1(\tau+\theta_S)$ for $\tau\in [0,S-\theta_S]$ and to be equal to
an arbitrary $u\in U$ on the interval $(S-\theta_S, S]$. Then one can obtain
$$
\left|\frac{1}{S}\int_0^S q(u^1(\tau), y^1(\tau))d\tau - \frac{1}{S}\int_0^S q(\bar u(\tau), \bar y(\tau))d\tau \right|
$$
\vspace{-.2in}
$$
\leq
\frac{1}{S}\int_0^{S-\theta_S} |q(u^1(\tau+\theta_S),y^1(\tau+\theta_S))dt -  q(u^1(\tau+\theta_S),\bar y(\tau))|d\tau
+\ \frac{2M_q\theta_s}{S}\
$$
\vspace{-.2in}
\begin{equation}\label{eq-Prop-Special-Cas-2}
 \leq L_q \ \alpha_1(S) + \frac{2M_q\theta_s}{S} \ \leq\ L_q \alpha_1(S) + 2M_q \alpha_2(S)\BYDEF \alpha_q(S),
\end{equation}
where $L_q$ is a Lipschitz constant of $q(\cdot)$ and $\ M_q\BYDEF \max_{(u,y)\in U\times Y}|q(u,y)|$.
Define now $u^2(\tau)$ to be equal to an arbitrary $u\in U$ for $\tau\in [0,\theta_S)$ and to be equal to  $\bar u (\tau -\theta_S)$ for
$\tau\in [\theta_S ,S]$. Then,
 similarly to  (\ref{eq-Prop-Special-Cas-1}), one obtains
$$
\left|\frac{1}{S}\int_0^S q(u^2(\tau), y^2(\tau))d\tau - \frac{1}{S}\int_0^S q(\bar u(\tau), \bar y(\tau))d\tau \right|
$$
\vspace{-.2in}
$$
\leq
\frac{1}{S}\int_0^{S-\theta_S} |q(u^2(\tau+\theta_S),y^2(\tau+\theta_S))dt -  q(\bar u(\tau),\bar y(\tau))|d\tau +\ \frac{2M_q\theta_s}{S}
$$
\vspace{-.2in}
$$
=
\frac{1}{S}\int_0^{S-\theta_S} |q(\bar u(\tau),y^2(\tau+\theta_S))dt -  q(\bar u(\tau),\bar y(\tau))|d\tau +\ \frac{2M_q\theta_s}{S}
$$
\vspace{-.2in}
\begin{equation}\label{eq-Prop-Special-Cas-3}
\ \leq L_q \ \alpha_1(S) + \frac{2M_q\theta_s}{S} \ \leq\ L_q \alpha_1(S) + 2M_q \alpha_2(S)= \alpha_q(S),
\end{equation}
The validity of (\ref{eq-Prop-Special-Cas-2}) and (\ref{eq-Prop-Special-Cas-3}) implies (\ref{e-CSO-further-5}) with
$\beta_q(S)= 2\alpha_q(S)$.
$\ \Box$


\section{Optimal and near optimal average control generating families}\label{Sec-ACG}
In this section, we  discuss a way of characterization and  construction of optimal or near optimal  average control generating (ACG) families (a concept introduced in  \cite{GM-2015} and \cite{GR-long} for a different class of problems)
and , in the next section, we will demonstrate
how these can be used for the construction of asymptotically optimal or asymptotically near optimal solutions
 of the perturbed problem.

Let $u_z(\tau) $ be a parameterized by $z$ family of controls and let $y_z(\tau)\in Y_z $ be a family of solutions of the reduced system (\ref{e-CSO-2}) obtained with these controls
 such that, for any $z\in Z$,  the pair $(u_z(\cdot),y_z(\cdot))$ generates the occupational measure $\mu_z(du,dy)$
(see (\ref{e:oms-0-1-infy})),
 with the integral
$\int_{U\times Y}q(u,y)\mu_z(du,dy) $ being a measurable function of $z$, and
\begin{equation}\label{e-opt-OM-1-0}
 | S^{-1}\int_0^S q(u_{z}(\tau),y_{z}(\tau))d\tau - \int_{U\times Y}q(u,y)\mu_z(du,dy))|\leq \alpha_q(S)\ \ \forall z\in Z, \ \  \ \ \ \ \ \ \
 \lim_{S\rightarrow\infty}\alpha_q(S)
\end{equation}
 for any
continuous $q(u,y) $.
Note that the estimate (\ref{e-opt-OM-1-0}) is valid if $y_z(\cdot)$ and $u_z(\cdot)$ are $T_z$-periodic,
with $T_z $ being uniformly bounded on $ Z$. Note also that (similarly to the proof of Proposition \ref{Prop-inv-measure-set-2})
it can be verified that
\begin{equation}\label{e-opt-OM-1-0-11}
 \mu_z\in W(z) \ \ \ \ \forall \ z\in Z.
\end{equation}
\begin{Definition}\label{Def-ACG}
The family $(u_z(\cdot),y_z(\cdot))$ will be called {\it average control generating} (ACG) if the system
\begin{equation}\label{e-opt-OM-1}
z'(t)=\tilde h (\mu_{z(t)}), \ \ \  \  z(0) = z_0 ,
\end{equation}
where
\begin{equation}\label{e-opt-OM-1-101}
\tilde h (\mu_{z})= \int_{U\times Y}h(u,y)\mu_z(du,dy),
\end{equation}
has a unique solution  $z(t)\in  Z\ \forall t\in [0,\infty)$ (thus, the pair $(\mu_{z(t)}, z(t)) $ is a solution of the averaged
system (\ref{e-CSO-17}), (\ref{e-CSO-18})).
 \end{Definition}
   \begin{Definition}\label{Def-ACG-opt}
 An ACG family  $(u_z(\cdot),y_z(\cdot))$   is called  optimal  if the
 pair $(\mu_{z(t)}, z(t)) $, where $z(t)$ is the solution of (\ref{e-opt-OM-1}),
 is an optimal solution of the averaged problem (\ref{Vy-ave}). That is,
\begin{equation}\label{e-opt-OM-2-101}
\int_0 ^ { + \infty }e^{- C t}\tilde{r}(\mu_{z(t)})dt
 = \tilde{R}^*(z_0).
 \end{equation}
An ACG family  $(u_z(\cdot),y_z(\cdot))$   is called   {\it $\alpha $-near optimal} ($\alpha >0$)
 if
 \begin{equation}\label{e-opt-OM-2-101-1}
\int_0 ^ { + \infty }e^{- C t}\tilde{r}(\mu_{z(t)})dt
 \leq \tilde{R}^*(z_0) + \alpha.
 \end{equation}
 \end{Definition}

Let   $\Phi $ stand for the graph of $W(\cdot) $
(see (\ref{e-CSO-7})):
\begin{equation}\label{e:graph-w}
\Phi\BYDEF \{(\mu , z) \ : \ \mu\in W(z), \ \ z\in Z \} \ \subset \mathcal{P}(U\times Y)\times Z \ .
 \end{equation}
and let  the set $\tilde{\mathcal{D}}(z_0)$ be defined  by the equation
\begin{equation}\label{D-SP-new}
\tilde{\mathcal{D}}(z_0)\BYDEF \{ p \in {\cal P} (\Phi): \;
 \int_{\Phi }(\nabla\psi(z)^T \tilde{h}(\mu) + C (\psi (z_0) - \psi (z) )) p(d\mu,dz) = 0 \ \ \ \forall \psi(\cdot) \in C^1(\R^k)\}
  \end{equation}
  In what follows, it is always assumed that
  \begin{equation}\label{D-SP-new-101}
\tilde{\mathcal{D}}(z_0)\neq\emptyset .
  \end{equation}
    Note that, the latter is satisfied if   (\ref{e-CSO-20}) is valid (since the discounted occupational measure generated by any solution
  $(\mu(\cdot),z(\cdot))$ of the averaged system is contained in $\tilde{\mathcal{D}}(z_0)$; see  Lemma 2.1 in \cite{GQ-1}).

  Let us consider the following optimization problem.
\begin{equation}\label{e-ave-LP-opt-di}
\min_{p\in \tilde{\mathcal{D}}(z_0)}\int_{\Phi}\tilde r(\mu)p(d\mu , dz)\BYDEF \tilde{a}^*(z_0).
\end{equation}
Note  that this problem belongs to the class of infinite dimensional linear programming (IDLP) problems (since both the objective function and the constraints are linear in the \lq\lq decision variable" $p$). The
  optimal solution of the problem (\ref{e-ave-LP-opt-di}) exists (since, as can be readily verified,
 $\tilde{\mathcal{D}}(z_0)$ is weak$^*$ compact).

The IDLP problem (\ref{e-ave-LP-opt-di}) is closely related to the averaged optimal control problem
(\ref{Vy-ave}). In particular, the optimal values of these problem are related by the inequality (see Proposition 2.2 in \cite{GQ} and Lemma 2.1 in \cite{GQ-1}))
\begin{equation}\label{D-SP-new-2-inequality}
C \tilde{R}^*(z_0)\ \geq \ \tilde{a}^*(z_0),
  \end{equation}
which,
under certain  conditions,
becomes the equality
  \begin{equation}\label{D-SP-new-2}
C \tilde{R}^*(z_0)\ = \ \tilde{a}^*(z_0).
  \end{equation}
In fact, the equality  (\ref{D-SP-new-2}) is valid if, by allowing   the use of relaxed controls, we would not improve the optimal value of the averaged problem  (see Theorem 2.2 and  Corollary 2.3 in \cite{GQ-1}). It seems likely that the averaged problem does satisfy this property but we do not verify it in the paper.

  Consider the \lq\lq maximin" problem
\begin{equation}\label{e:DUAL-AVE}
\sup_{\zeta(\cdot)\in C^1(\R^k)} \min_{z\in Z}\{ \tilde H(\nabla \zeta(z), z) + C(\zeta(z_0) - \zeta(z))\}= \tilde{a}^*(z_0) ,
\end{equation}
where $sup $ is sought over all continuously differentiable functions $\zeta(\cdot):\R^k\rightarrow \R^1 $ and
$\tilde H(p,z)$ stands for the
 Hamiltonian of the averaged
system,
\begin{equation}\label{e:H-tilde}
\tilde H(p,z)\BYDEF \min_{\mu\in W(z)} \{\tilde r(\mu)+p^T \tilde h (\mu)\}.
\end{equation}
The problem (\ref{e:DUAL-AVE}) is dual with respect to the IDLP problem (\ref{e-ave-LP-opt-di}), the duality relationships
between the two problems include the equality of their optimal values
 (see Theorem 3.1 in \cite{GQ}).
 For convenience, (\ref{e:DUAL-AVE})  will be referred to as the  {\it  averaged dual problem}.

The right-hand-side in (\ref{e:H-tilde}) can be rewritten in the form of the problem (see (\ref{e-CSO-16}) and (\ref{Vy-ave-r}))
\begin{equation}\label{e:H-tilde-10-1}
 \min_{\mu\in W(z)} \{\int_{U\times Y}[r(u,y)+p^T  h (u,y )]\mu(du,dy)\} = \tilde H (p,z),
\end{equation}
which is of the IDLP class as well (since its objective functions and the constraints defining $W(z)$ are linear
in $\mu$; see (\ref{e-CSO-7})). The problem dual with respect to (\ref{e:H-tilde-10-1}) is the maximin problem
\begin{equation}\label{e:dec-fast-4}
  \sup_{\eta(\cdot)\in C^1(\R^m)} \min_{(u,y)\in U\times Y_z}\{r(u,y)+p^T  h (u,y) + \nabla \eta (y)^T  f(u,y)
\} =\tilde H (p,z),
\end{equation}
where $sup $ is sought over all continuously differentiable functions $\eta(\cdot): \R^m\rightarrow \R^1 $.
The  optimal values of the problems (\ref{e:H-tilde-10-1}) and (\ref{e:dec-fast-4}) are equal,
  this being one of the duality relationships between these two problems
 (see Theorem 4.1 in \cite{FinGaiLeb}). The problem (\ref{e:dec-fast-4}) will be referred to as the {\it associated dual problem}.

 \begin{Assumption}\label{Ass-Existence-Duals}\
 
 (i) A solution of the averaged dual problem exists. That is, there exists $\zeta^*(\cdot)\in C^1(\R^k) $ such that
  \begin{equation}\label{e:DUAL-AVE-0-1}
\min_{z\in Z}\{\tilde H(\nabla\zeta^*(z),z) + C (\zeta^*(z_0) - \zeta^*(z))\} =\tilde{a}^*(z_0).
\end{equation}
(ii)  A solution of the associated dual problem considered with $p= \nabla \zeta^*(z)$ (where $\zeta^*(z)$ is a solution of the averaged dual problem) exists for any $z\in Z$. That is, there exists $\eta_z^*(\cdot)\in C^1(\R^m) $  such that
 \begin{equation}\label{e:dec-fast-4-0-1}
\min_{(u,y)\in U\times Y_z}\{r(u,y)+\nabla \zeta^*(z)^T  h (u,y) + \nabla \eta_z^* (y)^T  f(u,y)\}= \tilde H (\nabla \zeta^*(z),z)
\ \ \ \forall u\in U\ \ \forall z\in Z.
\end{equation}
\end{Assumption}

Note that from that (\ref{e:H-tilde}) and (\ref{e:DUAL-AVE-0-1}) it follows that
\begin{equation}\label{e:dec-fast-4-0-12}
\min_{(\mu ,z) \in \Phi }\{\tilde r(\mu )+ \nabla \zeta^*(z)^T \tilde h (\mu ) + C (\zeta(z_0) - \zeta(z))\} = \tilde{a}^*(z_0)  ,
\end{equation}
and
 from (\ref{e:DUAL-AVE-0-1}) and (\ref{e:dec-fast-4-0-1}) it follows that
\begin{equation}\label{e:dec-fast-4-0-10}
\min_{z\in Z}\min_{(u,y)\in U\times Y_z} \{r(u,y)+\nabla \zeta^*(z)^T  h (u,y) + \nabla \eta_z^* (y)^T  f(u,y)+ C (\zeta^*(z_0) - \zeta^*(z))\} =  \tilde{a}^*(z_0).
\end{equation}

The result stated below gives  a sufficient  condition for an ACG family $(u_z(\cdot),y_z(\cdot))$ to be optimal. This  also
becomes a necessary condition of optimality  provided that the ACG family is periodic, that is,
\begin{equation}\label{e:H-tilde-10-3-1001}
(u_z(\tau + T_z),y_z(\tau + T_z)) = (u_z(\tau ),y_z(\tau )) \ \ \ \ \forall \ \tau \geq 0
\end{equation}
for some $\ T_z > 0 $.

\begin{Proposition}\label{Prop-necessary-opt-cond} Let Assumption \ref{Ass-Existence-Duals} be satisfied.
  Then, for an ACG family $(u_z(\cdot),y_z(\cdot))$ to be optimal and for the equality (\ref{D-SP-new-2}) to be valid it is sufficient and, provided that (\ref{e:H-tilde-10-3-1001}) is satisfied, also necessary that
  there exist sets $P_t\subset [0,\infty) $, $A\subset [0,\infty) $ such that
  \begin{equation}\label{e:H-tilde-10-3-1002-1-1}
 meas\{[0,\infty)\setminus P_t\} = 0 \ \ \ \forall t\in A \ \ \ \ \ {\rm and} \ \ \ \ \  meas\{[0,\infty)\setminus A\} = 0,
\end{equation}
($meas\{\cdot\} $ stands for the Lebesgue measure of the corresponding set) and such that
$$
r(u_{z(t)}(\tau),y_{z(t)}(\tau))+\nabla \zeta^*(z(t))^T  h (u_{z(t)}(\tau),y_{z(t)}(\tau))
+ \nabla \eta_{z(t)}^* (y_{z(t)}(\tau))^T  f(u_{z(t)}(\tau),y_{z(t)}(\tau))
$$
\vspace{-.2in}
\begin{equation}\label{e:dec-fast-4-0-9}
 +\  C (\zeta^*(z_0) - \zeta^*(z(t)))=\tilde{a}^*(z_0)  \ \ \ \ \ \forall \ \tau\in P_t, \ \ \forall \ t\in A,
\end{equation}
 where $z(t) $ is  the  solution of (\ref{e-opt-OM-1}) and $\mu_{z(t)}$ in (\ref{e-opt-OM-1})
is the occupational measure generated by $(u_{z(t)}(\cdot),y_{z(t)}(\cdot))$.
 \end{Proposition}
 {\it Proof.}
 We assume first that (\ref{e:dec-fast-4-0-9}) is satisfied and show that the ACG family $(u_z(\cdot),y_z(\cdot))$ is optimal and that
 (\ref{D-SP-new-2}) is valid. By integrating (\ref{e:dec-fast-4-0-9}) over $\tau\in [0,S]$, dividing the result
 by $S $ and passing to the limit with $S\rightarrow\infty $ (and also having in mind (\ref{e-CSO-16}), (\ref{Vy-ave-r})
 and (\ref{e-opt-OM-1-0})), one obtains
 \begin{equation}\label{e:dec-fast-4-0-8}
 \tilde r(\mu_{z(t)}) + \nabla \zeta^*(z(t))^T  \tilde h(\mu_{z(t)})
 +\  C (\zeta^*(z_0) - \zeta^*(z(t)))=\tilde{a}^*(z_0) \ \ \ \forall t\in A,
\end{equation}
 where it has been taken into account that
 $$
 \lim_{S\rightarrow\infty}\frac{1}{S}\int_0^{S}\nabla \eta_{z(t)}^* (u_{z(t)}(\tau), y_{z(t)}(\tau))^T  f(y_{z(t)}(\tau))d\tau
 = \lim_{S\rightarrow\infty}\frac{1}{S}\left( \eta_{z(t)}^* (y_{z(t)}(S))- \eta_{z(t)}^* (y_{z(t)}(0))\right) = 0.
 $$
 By multiplying (\ref{e:dec-fast-4-0-8}) by $e^{-Ct} $ and integrating the result over $t$ from $0$ to $\infty$,
  one  obtains
 \begin{equation}\label{e:dec-fast-4-0-7}
 \int_0^{\infty}e^{-Ct}\tilde r(\mu_{z(t)})dt = \frac{1}{C}\tilde{a}^*(z_0)
 \end{equation}
 where it has been taken into account that
\begin{equation}\label{e:dec-fast-4-0-6}
\int_0^{\infty}e^{-Ct}\left( \nabla \zeta^*(z(t))^T  \tilde h(\mu_{z(t)})
 +\  C (\zeta^*(z_0) - \zeta^*(z(t)))\right) dt = \int_0^{\infty}\frac{d}{dt}\left(e^{-Ct}\zeta^*(z(t)) \right)
 + \zeta^*(z_0) = 0.
\end{equation}
 Due to (\ref{D-SP-new-2-inequality}), the validity of (\ref{e:dec-fast-4-0-7}) implies both the optimality of
 the ACG under consideration and the validity of (\ref{D-SP-new-2}).

 Assume now the ACG family is optimal and (\ref{D-SP-new-2}) is valid. Then (\ref{e:dec-fast-4-0-7}) is true, which along
 with (\ref{e:dec-fast-4-0-6}) implies that
 \begin{equation}\label{e:dec-fast-4-0-5}
\int_0^{\infty}e^{-Ct}\left(\tilde r(\mu_{z(t)}) + \nabla \zeta^*(z(t))^T  \tilde h(\mu_{z(t)})
 +\  C (\zeta^*(z_0) - \zeta^*(z(t)))-\tilde{a}^*(z_0)\right)dt = 0.
\end{equation}
Due to (\ref{e:dec-fast-4-0-12}), from (\ref{e:dec-fast-4-0-5}) it follows that
 \begin{equation}\label{e:dec-fast-4-0-51}
\tilde r(\mu_{z(t)}) + \nabla \zeta^*(z(t))^T  \tilde h(\mu_{z(t)})
 +\  C (\zeta^*(z_0) - \zeta^*(z(t)))=\tilde{a}^*(z_0) \ \ \ \forall t\in A,
\end{equation}
 where $A\subset [0,\infty)$ is such that $meas\{[0,\infty)\setminus A\} = 0$. Having  in mind the fact that, under the periodicity condition
(\ref{e:H-tilde-10-3-1001}),
$$
\tilde r(\mu_{z(t)})=\frac{1}{T_{z(t)}}\int_0^{T_{z(t)}}r(u_{z(t)}(\tau),y_{z(t)}(\tau))d\tau, \ \ \ \
\tilde h(\mu_{z(t)})=\frac{1}{T_{z(t)}}\int_0^{T_{z(t)}}h(u_{z(t)}(\tau),y_{z(t)}(\tau))d\tau
$$
and
$$
\int_0^{T_{z(t)}}\nabla \eta_{z(t)}^* (y_{z(t)}(\tau))^T  f(u_{z(t)}(\tau),y_{z(t)}(\tau))d\tau \
=\  \eta_{z(t)}^* (y_{z(t)}(T_{z(t)}))\ - \  \eta_{z(t)}^* (y_{z(t)}(0))\ = \ 0,
$$
one can rewrite (\ref{e:dec-fast-4-0-51}) in the form
$$
\frac{1}{T_{z(t)}}\int_0^{T_{z(t)}} (r(u_{z(t)}(\tau),y_{z(t)}(\tau))+\nabla \zeta^*(z(t))^T  h (u_{z(t)}(\tau),y_{z(t)}(\tau))
+ \nabla \eta_{z(t)}^* (y_{z(t)}(\tau))^T  f(u_{z(t)}(\tau),y_{z(t)}(\tau))
$$
\vspace{-.2in}
\begin{equation}\label{e:dec-fast-4-0-52}
 +\  C (\zeta^*(z_0) - \zeta^*(z(t)))-\tilde{a}^*(z_0))d\tau \ = \ 0\ \ \ \forall t\in A.
\end{equation}
The latter along with (\ref{e:dec-fast-4-0-10}) imply the validity of the equality (\ref{e:dec-fast-4-0-9}) for almost all $\tau\in [0,T_{z(t)}] $. This, in turn, implies  the validity of (\ref{e:dec-fast-4-0-9}) for $t\in P_t $, where $P_t$ is such that
$meas\{[0,\infty)\setminus P_t\} = 0$. This completes the proof.
 $\ \Box$

\begin{Corollary}\label{Cor-necessary}
Let Assumption \ref{Ass-Existence-Duals} be satisfied and the equality (\ref{D-SP-new-2}) be valid. Then for an  ACG family satisfying the periodicity condition
(\ref{e:H-tilde-10-3-1001}) to be optimal it is necessary that,  for any $\tau\in P_t$ and any $t\in A$,
\begin{equation}\label{e:H-tilde-10-3-1002}
u_{z(t)}(\tau)\in {\rm Argmin}_{u\in U}\{ r(u,y_{z(t)}(\tau))+\nabla \zeta^* (z(t))^T  h (u,y_{z(t)}(\tau))
+ \nabla \eta^*_{z(t)}(y_{z(t)}(\tau))^T f(u,y_{z(t)}(\tau) ) \} 
\end{equation}
 \end{Corollary}
{\it Proof.}
If the ACG $(u_z(\cdot), y_z(\cdot)) $ is optimal, then (\ref{e:dec-fast-4-0-9}) is satisfied. The latter and
(\ref{e:dec-fast-4-0-10}) imply that
$$
(u_{z(t)}(\tau), y_{z(t)}(\tau))\in {\rm Argmin}_{(u,y)\in U\times Y}\{ r(u,y)+\nabla \zeta^* (z(t))^T  h (u,y)+ \nabla \eta^*_{z(t)} (y)^T  f(u,y )  \}\ \ \ \ \forall \tau\in P_t, \ \ \forall t\in A,
$$
which, in turn, implies (\ref{e:H-tilde-10-3-1002}).
$\ \Box$

Having in mind (\ref{e:H-tilde-10-3-1002}), one may pose a question whether the feedback control
$u^*(y,z)$ defined as a minimizer
\begin{equation}\label{e:F-control-1}
u^*(y,z)\BYDEF {\rm argmin}_{u\in U}\{ r(u,y)+\nabla \zeta^* (z)^T  h (u,y)+\eta^*_{z} (y)^T  f(u,y)\}
\end{equation}
generates an optimal ACG family. 

\begin{Definition}\label{Def-feedback-control} A feedback control function $u^*(y,z)$  is referred to as that generating an optimal ($\alpha$-near optimal) ACG family if a pair $(u_z^*(\cdot),y_z^*(\cdot))$,
where $y_z^*(\tau)\in Y_z $ is a solution of the system
\begin{equation}\label{e-new-101}
 \frac{dy(\tau)}{d\tau}= f(u^*(y(\tau),z), y(\tau))
 \end{equation}
and  $u_z^*(\tau)=u^*(y_z^*(\tau),z) $, generates an optimal ($\alpha$-near optimal) ACG family.
 \end{Definition}

\begin{Proposition}\label{Prop-construction-opt}
Let Assumption \ref{Ass-Existence-Duals} be satisfied and  the equality (\ref{D-SP-new-2}) be valid.
Let $u^*(y,z)$ be as in (\ref{e:F-control-1}). Assume that, for any $z\in Z$,
  the set  of occupational measures $W^*(z) $ that are generated by the pairs $(u_z^*(\cdot),y_z^*(\cdot))$,
 where $y_z^*(\tau) $ is a solution of the system (\ref{e-new-101}) and $u_z^*(\tau)=u^*(y_z^*(\tau),z) $,
 is a singleton, $W^*(z)=\{\mu_z^*\}$.
  Assume also that an optimal ACG family $(  u_z(\cdot),  y_z(\cdot)) $ satisfying the periodicity condition
(\ref{e:H-tilde-10-3-1001}) exists and
the problem
\begin{equation}\label{e:F-control-4}
\min_{u\in U}\{ r(u, y_{ z(t)}(\tau))+\nabla \zeta^* ( z(t))^T  h (u,  y_{  z(t)}(\tau)) + \eta^*_{z(t)} (y_{z(t)}(\tau))^T  f(u,y_{z(t)}(\tau) ) \}
\end{equation}
has a unique minimizer for almost all $\tau\in [0,\infty) $ and  almost all $t\in [0,\infty) $, where
$ z(t)$ is the solution of the system
(\ref{e-opt-OM-1}).
 Then $u^*(y,z)$
generates an optimal
  ACG family.
\end{Proposition}
{\it Proof.}
By Corollary \ref{Cor-necessary},
$  u_{ z(t)}(\tau)$ is a minimizer in the problem (\ref{e:F-control-4}) for all $ \  \tau\in P_t$ and for all $  t\in A $. Due to the fact that the minimizer of this problem is unique,
$$
 u_{ z(t)}(\tau) = u^*( y_{z(t)}(\tau), z(t))
$$
for almost all $\tau\in P_t $ and almost all $t\in A $. Since $\ meas\{[0,\infty)\setminus P_t\}=0  $, $\ y_{z(t)}(\tau) $ is a solution of
(\ref{e-new-101}) and, thus, the occupational measure $\mu_{z(t)}$ generated by $(  u_z(\cdot),  y_z(\cdot)) $ is an element of $W^*(z(t))$ for almost all $t\in A $. The set $W^*(z(t))$ consists of only one element $\mu^*_{z(t)} $ and, hence, $\mu_{z(t)}= \mu^*_{z(t)} $
for almost all $t\in [0,\infty)$ (remind that $\ meas\{[0,\infty)\setminus A\} = 0  $).
 This proves the proposition.
$\ \Box$

Solutions of the  averaged and associated dual problems may not exist (that is, Assumption \ref{Ass-Existence-Duals} may not be satisfied) and, even if they exist, finding these solutions analytically is hardly possible. A way of finding approximate solutions of the dual problems  numerically
was proposed in  \cite{GR-long} for a different class of problems. The approach of \cite{GR-long} is applicable for the class of problems under consideration as well. Below, we briefly describe a way  this approach can be used for construction of feedback controls  that generate near optimal ACG families.

Let $\psi_i(\cdot)\in C^1(\R^k) \ , \ \ i=1,2,...,$ be  a sequence of functions such that any $\zeta(\cdot)\in C^1(\R^k)$ and its gradient are simultaneously approximated by a linear
combination of $\psi_i(\cdot)$ and their gradients. Also, let $\phi_i(\cdot)\in C^1(\R^m) \ , \ \ i=1,2,...,$ be  a sequence of functions such that any $\eta(\cdot)\in C^1(\R^m)$ and its gradient are simultaneously approximated by a linear
combination of $\phi_i(\cdot)$ and their gradients. Examples of such sequences are monomials
 $z_1^{i_1}...z_k^{i_k}$, $i_1,...,i_k =0,1,...$ and, respectively, $y_1^{i_1}...y_m^{i_m}$, $i_1,...,i_m =0,1,...$, with $z_j \ (j=1,...,k)$ and
 $y_l\ (l=1,...,m) $ standing for the components of $z$ and $y$ (see, e.g., \cite{Llavona}).

Let us introduce the following notations:
\begin{equation}\label{e:2.4-M}
W_M(z) \BYDEF  \{\mu \in \mathcal{P}(U \times Y) \ : \
\int_{U\times Y} \nabla\phi_i(y)^{T} f(u,y,z)
\mu(du,dy) =  0, \ \ \ i=1,...,M \} ,
\end{equation}
 \vspace{-.2in}
  \begin{equation}\label{e:graph-w-101}
\Phi_M\BYDEF \{(\mu , z) \ : \ \mu\in W_M(z), \ \ z\in Z \} \ \subset \mathcal{P}(U\times Y)\times Z,
 \end{equation}
 \vspace{-.2in}
 \begin{equation}\label{e:graph-w-M-1}
\tilde H_M (p,z)\BYDEF \min_{\mu\in W_M(z)}\{\tilde r(\mu )+ p^T \tilde h (\mu )\}
 \end{equation}
 (compare with (\ref{e-CSO-7}), (\ref{e:graph-w})
  and (\ref{e:H-tilde})) and let us consider the maximin  problem
 \begin{equation}\label{e:DUAL-AVE-0-approx-MN}
\sup_{\zeta(\cdot)\in \mathcal{Q}_N} \min_{z\in Z}\{\tilde H_M (\nabla \zeta(z),z) + C (\zeta(F(y_0)) - \zeta(z)) \}\BYDEF \tilde{a}_{M,N}(z_0),
\end{equation}
where $sup$ is over the functions from the finite dimensional space $\mathcal{Q}_N\subset  C^1(\R^k)$,
 \begin{equation}\label{e:DUAL-AVE-0-approx-MN-dim} \mathcal{Q}_N\BYDEF \{\zeta(\cdot)\in C^1(\R^k): \zeta(z)= \sum_{i=1}^N \lambda_i \psi_i(z), \ \ \lambda=
(\lambda_i)\in \mathbb{R}^N\}\end{equation}
(compare with (\ref{e:DUAL-AVE})).
The problem (\ref{e:DUAL-AVE-0-approx-MN-dim}) is dual with respect to the semi-infinite LP problem
\begin{equation}\label{e-ave-LP-opt-di-MN}
\min_{p\in \tilde{\mathcal{D}}_{M,N}(z_0)}\int_{F_M}\tilde r(\mu)p(d\mu , dz)= \tilde{a}_{M,N}(z_0),
\end{equation}
where
 \begin{equation}\label{D-SP-new-MN}
\tilde{\mathcal{D}}_{M,N}(z_0)\BYDEF \{ p \in {\cal P} (\Phi_M): \;
 \int_{\Phi_M }(\nabla\psi_i(z)^T \tilde{h}(\mu) + C (\psi_i (z_0) - \psi_i (z) )) p(d\mu,dz) = 0, \ \ \ i=1,...,N\},
  \end{equation}
(compare with  (\ref{D-SP-new}) and (\ref{e-ave-LP-opt-di})). Note that
\begin{equation}\label{e-MN-converg-1}
\tilde{a}_{M,N}(z_0)\leq\tilde{a}^{*}(z_0) \ \ \forall M,N=1,2,... \ ,\ \ \ \ \ \ \ \ \ \ \ \lim_{N\rightarrow\infty} \lim_{M\rightarrow\infty}\tilde{a}_{M,N}(z_0)=\tilde{a}^{*}(z_0)
\end{equation}
(see Proposition 7.1 in \cite{GR-long}).
Denote by $\zeta^{M,N}(z)$ a solution of the problem (\ref{e:DUAL-AVE-0-approx-MN}) (assuming that it exists)
and consider the problem
\begin{equation}\label{e:dec-fast-4-0-1-MN}
\sup_{\eta(\cdot)\in \mathcal{V}_M}\min_{(u,y)\in U\times Y_z}\{r(u,y)+\nabla \zeta^{M,N}(z)^T  h (u,y) + \nabla \eta(y)^T  f(u,y)\}= \tilde H_M (\nabla \zeta^{M,N}(z),z)
\ \ \ \ \  \forall z\in Z,
\end{equation}
where $sup$ is over the functions from the finite dimensional space $\mathcal{V}_M\subset C^1(\R^m) $,
\begin{equation}\label{e:HJB-8-Associated}
\mathcal{V}_M\BYDEF \{\eta(\cdot)\in C^1(\R^m): \eta(y)= \sum_{i=1}^M \omega_i \phi_i(y), \ \ \omega=
(\omega_i)\in \mathbb{R}^M\}.
\end{equation}
The problem (\ref{e:HJB-8-Associated}) is dual with respect to the problem in the right-hand-side of (\ref{e:graph-w-M-1}) taken with $p= \nabla \zeta^{M,N}(z)$. Denote by $\eta_z^{M,N}(y) $ a solution
of the problem (\ref{e:dec-fast-4-0-1-MN}) (assuming that it exists too) and define $u^{M,N}(y,z)$ as a minimizer
\begin{equation}\label{e:F-control-1-MN}
u^{M,N}(y,z)\BYDEF {\rm argmin}_{u\in U}\{ r(u,y)+\nabla \zeta^{M,N} (z)^T  h (u,y) + \nabla \eta^{M,N}_z (y)^T  f (u,y)\},
\end{equation}
For the class of problems considered in  \cite{GR-long}, it has been shown that, under certain conditions,  the feedback control $u^{M,N}(y,z)$ generates a $\nu(M,N)$-near optimal ACG family, with  $\ \lim_{N\rightarrow\infty}\lim_{M\rightarrow\infty}\nu(M,N)=0 $ (see Theorem 9.4 in  \cite{GR-long}).
The argument used in \cite{GR-long} can be extended to the class of problems that we are dealing with.  We, however, do not present this argument in the paper, and we conclude by just indicating a  way to numerically evaluate a \lq\lq degree of near optimality" of the control $u^{M,N}(y,z)$ defined in (\ref{e:F-control-1-MN}) (provided that it generates an ACG family).
 Note first of all that, if solutions of the problems (\ref{e:DUAL-AVE-0-approx-MN}) and (\ref{e:dec-fast-4-0-1-MN}) exist, then they are presentable in the form
\begin{equation}\label{e:DUAL-AVE-0-approx-MN-1}
\zeta^{M,N}(z)=\sum_{i=1}^N \lambda_i^{M,N} \psi_i(z), \ \ \ \  
\eta^{N,M}_z(y)= \sum_{i=1}^M \omega_{z,i}^{N,M} \phi_i(y).
\end{equation}
The coefficients of the expansions in (\ref{e:DUAL-AVE-0-approx-MN-1}) can be numerically found via solving the semi-infinite LP problem (\ref{e-ave-LP-opt-di-MN}) (see the description of the algorithm in Section 4.3 of \cite{GR-long}). Once these coefficients are found, one can define $u^{M,N}(y,z)$ according to (\ref{e:F-control-1-MN}) and verify whether the system
\begin{equation}\label{e-new-101-1}
 \frac{dy(\tau)}{d\tau}= f(u^{M,N}(y(\tau),z), y(\tau)),\ \ \ \ \ \ y(0)=y_0\in Y_z
 \end{equation}
has a  solution $y_z^{M,N}(\tau) $ such that   the pair $(u_z^{M,N}(\cdot),y_z^{M,N}(\cdot))$,
where $u_z^{M,N}(\tau)=u^{M,N}(y^{M,N}(\tau),z) $, generates an occupational measure $\mu^{M,N}_z $.
Finally, one can numerically integrate the system
\begin{equation}\label{e-opt-OM-1-101-101}
z'(t)=\tilde h (\mu^{M,N}_{z(t)}), \ \ \  \  z(0) = z_0 ,
\end{equation}
and (provided that the solution $z^{M,N}(t) $ of the latter exists and is contained in $Z$) evaluate the  integral
\begin{equation}\label{e-opt-OM-1-101-102}
\int_0 ^ { + \infty }e^{-C t}\tilde r(\mu^{M,N}_{z^{M,N}(t)})dt\BYDEF \tilde R^{M,N}(y_0),
\end{equation}
the latter being the value of the objective function of the averaged system (\ref{Vy-ave}) obtained with  $\mu(t)=\mu^{M,N}_{z^{M,N}(t)}$.
By (\ref{D-SP-new-2-inequality}) and (\ref{e-MN-converg-1}),
\begin{equation}\label{e-opt-OM-1-101-103}
\tilde R^{M,N}(y_0)- R^{*}(y_0)\leq \tilde R^{M,N}(y_0) - \frac{\tilde{a}_{M,N}(z_0)}{C},
\end{equation}
and, thus, the right-hand-side in (\ref{e-opt-OM-1-101-103}) gives an estimate of near optimality of the ACG family generated by the control $u^{M,M}(y,z) $.

\section{Construction of  asymptotically  optimal/near optimal controls for the perturbed problem}\label{Sec-near-opt-construction}
In this section, we will describe a way how asymptotically optimal or near optimal controls for the perturbed problem can be constructed on the basis of controls that generate optimal/near optimal ACG families. Everywhere in what follows  it is assumed that all the conditions of Corollary \ref{Corollary-equailty-values} are satisfied and, hence, (\ref{e-opt-val-ineq-1}) is valid (the limit of the optimal value of the perturbed problem exists and is equal to the optimal value of the averaged problem).
Moreover, for simplicity of the exposition, we  restrict ourselves to the consideration of the case when the reduced system is uncontrolled (that is, it is of the form (\ref{e-CSO-2-SC-I})) and we assume that  Assumption \ref{def-1} is valid. Note that the perturbed system (\ref{e-CSO}) takes  the form
\begin{equation}\label{e-CSO-101}
\frac{dy_{\epsilon}(t)}{dt}= \frac{1}{\epsilon}f(y_{\epsilon}(t)) + g(u(t), y_{\epsilon}(t)), \ \ \ \ \ y(0)=y_0
\end{equation}
in this case.
\begin{Definition}\label{Def-asympt-opt}
A control $u_{\epsilon}^*(\cdot) $  is called asymptotically optimal for the perturbed problem  (\ref{Vy-perturbed}) if
\begin{equation}\label{e-opt-convergence-1}
\lim_{\epsilon\rightarrow 0}\int_0 ^ { + \infty }e^{- C t}r(u_{\epsilon}^*(t), y_{\epsilon}^*(t)) dt=
\lim_{\epsilon\rightarrow 0}R^*(\epsilon, y_0),
\end{equation}
where $y_{\epsilon}^*(\cdot) $ is the solution of the system (\ref{e-CSO}) obtained with the control
$u_{\epsilon}^*(\cdot) $.
A control $u_{\epsilon}^*(\cdot) $ is called asymptotically
$\alpha$-near optimal ($\alpha > 0$)  for the perturbed problem  (\ref{Vy-perturbed}) if
\begin{equation}\label{e-opt-convergence-1-1}
\lim_{\epsilon\rightarrow 0}\int_0 ^ { + \infty }e^{- C t}r(u_{\epsilon}^*(t), y_{\epsilon}^*(t)) dt\leq
\lim_{\epsilon\rightarrow 0}R^*(\epsilon, y_0) + \alpha.
\end{equation}
\end{Definition}

Let $u^*(y,z)$ be a feedback control that generates an optimal or near optimal ACG family  and let $\mu^*_z $ be the occupational measure that is generated by the pair
$(u^*(y_z(\tau),z), y_z(\tau)) $, where $y_z(\cdot)$ is a solution of (\ref{e-CSO-2-SC-I}). The latter is $T_z$-periodic (due to Assumption \ref{def-1}) and, hence, for any continuous function $q(u,y)$,
\begin{equation}\label{e-per-measure}
\int_{Y_z}q(u,y)\mu^*_z(du,dy) = \frac{1}{T_z}\int_0^{T_z}q(u^*(y_z(\tau),z),y_z(\tau)) d\tau .
\end{equation}


Let us construct a control $u_{\epsilon}^*(t)$ that will be shown to be asymptotically optimal (near optimal) in the perturbed problem. To this end, partition the interval $[0,\infty)$ by the points
\begin{equation}\label{e:contr-rev-100-1}
t_l = l\Delta(\epsilon), \ l=0,1,...
\end{equation}
where
\begin{equation}\label{e-CSO-29}
\Delta(\epsilon) = \frac{\epsilon}{2L_f}\ln\frac{1}{\epsilon},
\end{equation}
$L_f $ being a Lipschitz constant of $f(\cdot)$.
Note that
 \begin{equation}\label{e:contr-rev-100-2}
\lim_{\epsilon\rightarrow 0}\Delta(\epsilon)=0, \ \ \ \ \  \ \ \ \lim_{\epsilon\rightarrow 0}\frac{\Delta(\epsilon)}{\epsilon}=\infty .
\end{equation}
Define $u_{\epsilon}^*(t) $ on the interval $[0,t_1)  $ by the equation
\begin{equation}\label{e:contr-rev-100-3}
u_{\epsilon}^*(t) = u^*(y_{y_0}(\frac{t}{\epsilon}),z_0) \ \ \ \forall \ t\in [0,t_1),
\end{equation}
where  $y_{y_0}(\tau)$ is the solution
of the  system (\ref{e-CSO-2-SC-I}) considered on the interval $ [0,\frac{t_{1}}{\epsilon}]  $
with the initial condition  $\ y_{y_0}(0)=y_0 $ (recall that $z_0\BYDEF F(y_0)$).

Let us assume that the control $u_{\epsilon}^*(t)$ has been defined on the interval $t\in[0,t_l) $ and let $y_{\epsilon}^*(t) $ be the corresponding solution of the perturbed system (\ref{e-CSO-101}) on this interval.
Extend  the definition of $u_{\epsilon}^*(t) $ to the interval $[0,t_{l+1}] $ by taking
\begin{equation}\label{e:contr-rev-100-4}
u_{\epsilon}^*(t) =  u^*(y_{y_{\epsilon}^*(t_l)}(\frac{t}{\epsilon}),z_{\epsilon}^*(t_l)) \ \ \ \ \forall \ t\in [t_l,t_{l+1}),\
\ \ \ \ \ \ l=1,2,...\ ,
\end{equation}
where
$y_{y_{\epsilon}^*(t_l)}(\tau) $ is the solution of the system (\ref{e-CSO-2-SC-I})  considered on the interval $\ [\frac{t_l}{\epsilon},\frac{t_{l+1}}{\epsilon}]  $
with the initial condition  $\ y_{y_{\epsilon}^*(t_l)}(\frac{t_l}{\epsilon})= y_{\epsilon}^*(t_l) $, and where $z_{\epsilon}^*(t_l)\BYDEF F(y_{\epsilon}^*(t_l))$.

\begin{Proposition}\label{Prop-asym-opt-contr}
Let $u^*(y,z)$ generate an optimal ($\alpha$-near optimal) ACG family.  Assume that
\begin{equation}\label{e-opt-OM-1-2-107}
y^*_{\epsilon}(t)\in Y, \ \ \ \ \ \ F(y^*_{\epsilon}(t))\in Z \ \ \ \ \ \forall t\in [0,\infty),
\end{equation}
where $y^*_{\epsilon}(\cdot) $ is the solution of the system (\ref{e-CSO-101}) obtained with use of the control $u^*_{\epsilon}(\cdot) $ constructed as described above. Assume also that
 the function $\tilde h^*(z) $,
\begin{equation}\label{e-opt-OM-1-2}
\tilde h^*(z)\BYDEF \int_{Y_z}h(u,y)\mu^*_z(du,dy) = \frac{1}{T_z}\int_0^{T_z}h(u^*(y_z(\tau),z),y_z(\tau)), d\tau 
\end{equation}
and the function
$\tilde r^*(z) $,
\begin{equation}\label{e-opt-OM-1-2-11}
\tilde r^*(z)\BYDEF \int_{Y_z}r(u,y)\mu^*_z(du,dy) =\frac{1}{T_z}\int_0^{T_z}r(u^*(y_z(\tau),z),y_z(\tau)) , d\tau 
\end{equation}
satisfy Lipschitz conditions in a neighbourhood of the solution $z^*(\cdot)$ of the system
\begin{equation}\label{e-opt-OM-1-1}
z'(t)=\tilde h^* (z(t)), \ \ \  \  z(0) = z_0 .
\end{equation}
 Then the control $u_{\epsilon}^*(t)$
  is asymptotically optimal ($\alpha$-near optimal) in the perturbed problem (\ref{Vy-perturbed}).
\end{Proposition}
{\it Proof.}
The proof of the proposition is given in Section \ref{Sec-Proofs-I}.
$\ \Box$

The construction of an asymptotically optimal (near optimal) control is simplified if  the function $u^*(y,z) $ is continuous.

 \begin{Proposition}\label{Prop-asym-opt-contr-simple}
Let  $u^*(y,z) $ that generate an optimal ($\alpha$-near optimal) ACG family  be continuous function of $(y,z)$ and let the solution $ y_{\epsilon}^*(\cdot)$ of the system
\begin{equation}\label{e-final-2}
\frac{dy_{\epsilon}^*(t)}{dt}= \frac{1}{\epsilon}f(y_{\epsilon}^*(t)) + g(u^*(y_{\epsilon}^*(t),F(y_{\epsilon}^*(t))), y_{\epsilon}^*(t)), \ \ \ \ \ y_{\epsilon}^*(0)=y_0
\end{equation}
satisfy the inclusions (\ref{e-opt-OM-1-2-107}). Let the functions $\tilde h^*(z), \ \tilde r^*(z) $ be Lipschitz continuous in a neighbourhood of the solution $z^*(\cdot)$ of the system (\ref{e-opt-OM-1-1}). Then the feedback
control $u^*(y,F(y)) $ is asymptotically optimal ($\alpha$-near optimal) in (\ref{Vy-perturbed}). That is,
\begin{equation}\label{e-final-1}
\lim_{\epsilon\rightarrow 0}\int_0 ^ { + \infty }e^{-C t}r(u^*(y_{\epsilon}^*(t),F(y_{\epsilon}^*(t))), y_{\epsilon}^*(t))dt =  \lim_{\epsilon\rightarrow 0}R^*(\epsilon, y_0) \ \ \ (\leq \lim_{\epsilon\rightarrow 0}R^*(\epsilon, y_0) + \alpha).
\end{equation}
\end{Proposition}

{\it Proof.}
The proof  is  in Section \ref{Sec-Proofs-I}.
$\ \Box$

\begin{Remark}
{\rm Note that the assumption that the functions $\tilde h^*(z)$ and $\tilde r^*(z)$ are Lipschitz continuous in a neighbourhood of $z^*(\cdot)$ can be replaced by a weaker assumption that they are piece-wise Lipschitz continuous in a neighbourhood of $z^*(\cdot)$ (see Definition 4.4 in \cite{GR-long}). The proof of the result is, however, much more technical in this case and we do not present it in the paper.}
\end{Remark}
\begin{Remark}\label{Remark-near-opt}
{\rm
Assume that the conditions of Proposition \ref{Prop-asym-opt-contr-simple} are satisfied with $u^*(y,z)=u^{M,N}(y,z)$ , where $u^{M,N}(y,z)$ is defined as a minimizer (\ref{e:F-control-1-MN}). Denote
  by $  R^{M,N}(\epsilon ,y_0)$ the value of the objective function (\ref{Vy-perturbed}) obtained with the use of the feedback control $u^{M,N}(y,F(y))$ in the perturbed system. Due to (\ref{e-opt-val-ineq-1}) and due to (\ref{D-SP-new-2-inequality}), (\ref{e-MN-converg-1}),
  $$
  \lim_{\epsilon\rightarrow 0}R^{M,N}(\epsilon ,y_0) - \lim_{\epsilon\rightarrow 0}R^{*}(\epsilon ,y_0)
 \ = \
  \lim_{\epsilon\rightarrow 0}R^{M,N}(\epsilon ,y_0) - \tilde R^{*}(y_0) \ \leq \ \lim_{\epsilon\rightarrow 0}R^{M,N}(\epsilon ,y_0) - \frac{\tilde{a}_{M,N}(z_0)}{C}
  $$
  (compare with (\ref{e-opt-OM-1-101-103})). Thus,
  $R^{M,N}(\epsilon ,y_0) - \frac{\tilde{a}_{M,N}(z_0)}{C}$
  may serve as a measure of near optimality of the feedback control $u^{M,N}(y,F(y))$ in the perturbed problem (provided that $\epsilon$ is small enough)}.
\end{Remark}

Let us demonstrate the construction of asymptotically near optimal control in the problem of minimization considered on the solutions of the system (\ref{e:ex-1}) (see Example 2 in Section \ref{S-SC1})

{\it Example 2 (continued).}
Consider the problem of optimal control
\begin{equation}\label{e:ex-2}
\epsilon \inf_{u(\cdot)} \int_0^{\infty} e^{-0.1 \epsilon  \tau}[u^2(\tau) + (F(y(\tau)) +2.05)^2]d\tau
\BYDEF R^*(\epsilon, y_0)
\end{equation}
where $inf$ is over all controls and the corresponding solutions of the system (\ref{e:ex-1}) that satisfy the initial condition $y(0) = y_0= (0.8916, \ 3.1370)$ (note that $z_0\BYDEF F(y_0)= -3$).
Differentiating the observable $z(\tau)= F(y(\tau))$ (with $F(y)$ being as in (\ref{e:ex-1-LV-101})) along a solution of the
 system (\ref{e:ex-1}), one can readily come to the conclusion that
its dynamics  is described  by the equation (compare with the general case (\ref{e-CSO-0}))
\begin{equation}\label{e:ex-1-1}
\frac{dz(\tau)}{d\tau}= \epsilon u(\tau)(y_1(\tau)-1), \ \ \ \ \ \ z(0)=z_0.
\end{equation}
Note that,  the solutions of (\ref{e:ex-1}) and (\ref{e:ex-1-1}) are written without the subscript $\epsilon$ (to simplify the notations).
 Having in mind the structure of the objective function $(\ref{e:ex-2})$,
one can expect that the control  minimizing (\ref{e:ex-2}) will make the state trajectory  spiral down from a larger orbit corresponding to $z_0=-3\ $ to a smaller orbit corresponding to $z_0=-2.05\ $ and that the solution of the augmented perturbed system (\ref{e:ex-1}), (\ref{e:ex-1-1}) obtained with this control stays in $Y\times Z$
(where $Y$ and $Z$ are defined in (\ref{e:ex-1-LV-101-1})). This kind of  behaviour  was indeed obtained with  the use of control constructed as described at the end of  Section \ref{Sec-ACG}.

In more detail, the problems (\ref{e:DUAL-AVE-0-approx-MN}), (\ref{e-ave-LP-opt-di-MN}) and (\ref{e:dec-fast-4-0-1-MN}) were formed with the use
of the monomials  $\psi(z)=z^i, \ i=1,..., 10, $ and $\phi(y)=y_1^{i_1}y_2^{i_2}, \ 1\leq i_1+i_2\leq 5 $. That is, $N=10$ (the dimension of the space $\mathcal{Q}_N $ in (\ref{e:DUAL-AVE-0-approx-MN-dim})) and
$M=35 $ (the  dimension of the space $\mathcal{V}_M$ in (\ref{e:HJB-8-Associated})).
 The coefficients  of the expansions (\ref{e:DUAL-AVE-0-approx-MN-1}) for $\zeta^{M,N}(z)$
 and $\eta^{M,N}_z (y) $    were numerically found with the help of a software implementation of the algorithm described in Section 4.3 of \cite{GR-long}. Note that the derivative of the function $\zeta^{M,N}(z)$ was verified to be negative
\begin{equation}\label{e-last-example-0}
\frac{d\zeta^{M,N}(z)}{dz}<0 \ \ \ \ \forall z\in [-3, - 2.05].
\end{equation}
  Having in mind the fact that in the given example
 $$f(u,y) = f(y), \ \ \ \
 \ h(u,y)= u(y_1-1), \ \ \ \  r(u,y)= u^2 + (F(y) +2.05)^2, \ \  \ \ U=[0,1],
 $$ one can readily establish that
 the minimizer in
(\ref{e:F-control-1-MN}) is of the form
\begin{equation}\label{eq-u-1}
u^{M,N}(y,z)= \min\{\ - \frac{1}{2}\frac{d\zeta^{M,N}(z)}{dz}(y_1 - 1),\ 1\} \ \ \ {\rm if} \ \ \ y_1 > 1;  \ \ \ \ \ \ \ \ \ \ \   u(y,z)= 0 \ \ \ {\rm if} \ \ \ y_1\leq 1.
\end{equation}
 Note that, due to (\ref{e-last-example-0}),   $u^{M,N}(y,z)$ is positive if $y_1 > 1 $ and it is zero if $y_1 \leq 1 $. Applying the feedback control
$u^{M,N}(y,F(y))$ in system (\ref{e:ex-1}) (considered with $\epsilon = 0.1 $) until the moment the state trajectory reaches the smaller orbit   characterized by the equality $F(y)=-2.05 $ and applying the \lq\lq zero control" after that moment, one obtains the solution $y(t)=(y_1(t),y_2(t)) $. The state trajectory of this solution and the graph of  the control $u(t)=  u^{M,N}(y(t),F(y(t))) $ are depicted in Figures 2 and 3. The graph of $y_1(t) $  is depicted in figure 4 (the graph of $y_2(t)$ looks similar).
The value of the objective function (\ref{e:ex-2}) thus obtained is  $\approx 1.25 $ \ ($R^{M,N}(\epsilon, y_0)\approx 1.25$ with $\epsilon =0.1 $).
The optimal value of the objective function of the semi-infinite LP problem  (\ref{e-ave-LP-opt-di-MN}) was numerically evaluated to be  $\approx 0.117$ \ ($ \tilde{a}_{M,N}(z_0)\approx 0.117)$ . That is, $\frac{\tilde{a}_{M,N}(z_0)}{C} \approx 1.17$ as $C=0.1$ in this case. Thus, the difference $1.25 - 1.17= 0.08 $ provides us with a measure of near optimality of the found solution (see Remark \ref{Remark-near-opt}).

 \begin{center}
\includegraphics[width=3.20in]{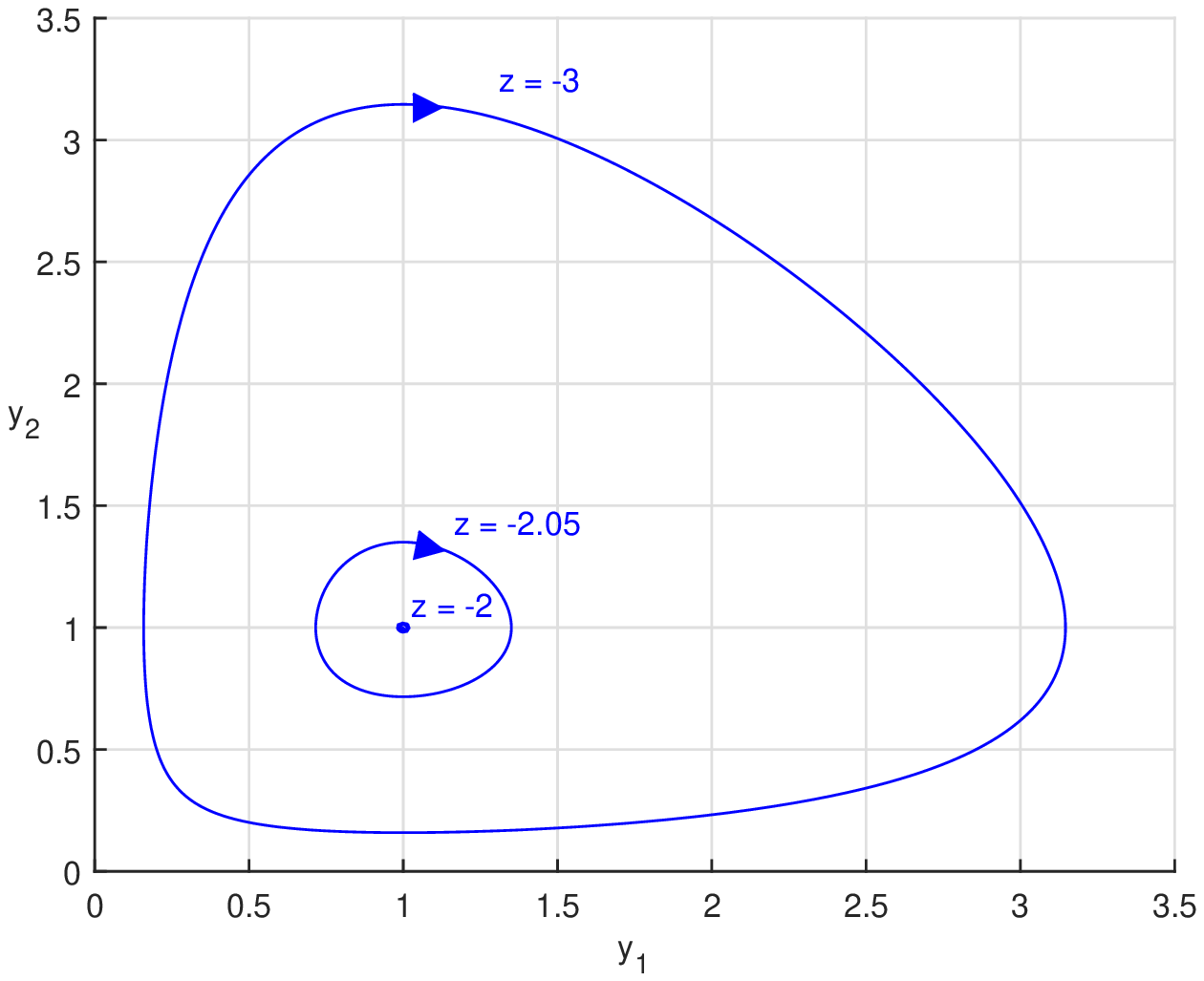}
\end{center}
\vspace{-.1in}
 \ \ \ \ \ \ \ \ \  \ \ \ \ \ \ \ \ \ \ \ \
\ \ \ \ \ \ \ \ \ \ \ \ \ \ \ Fig. 1: \ State trajectories for $\epsilon = 0$
  
  \begin{center}
\includegraphics[width=3.20in]{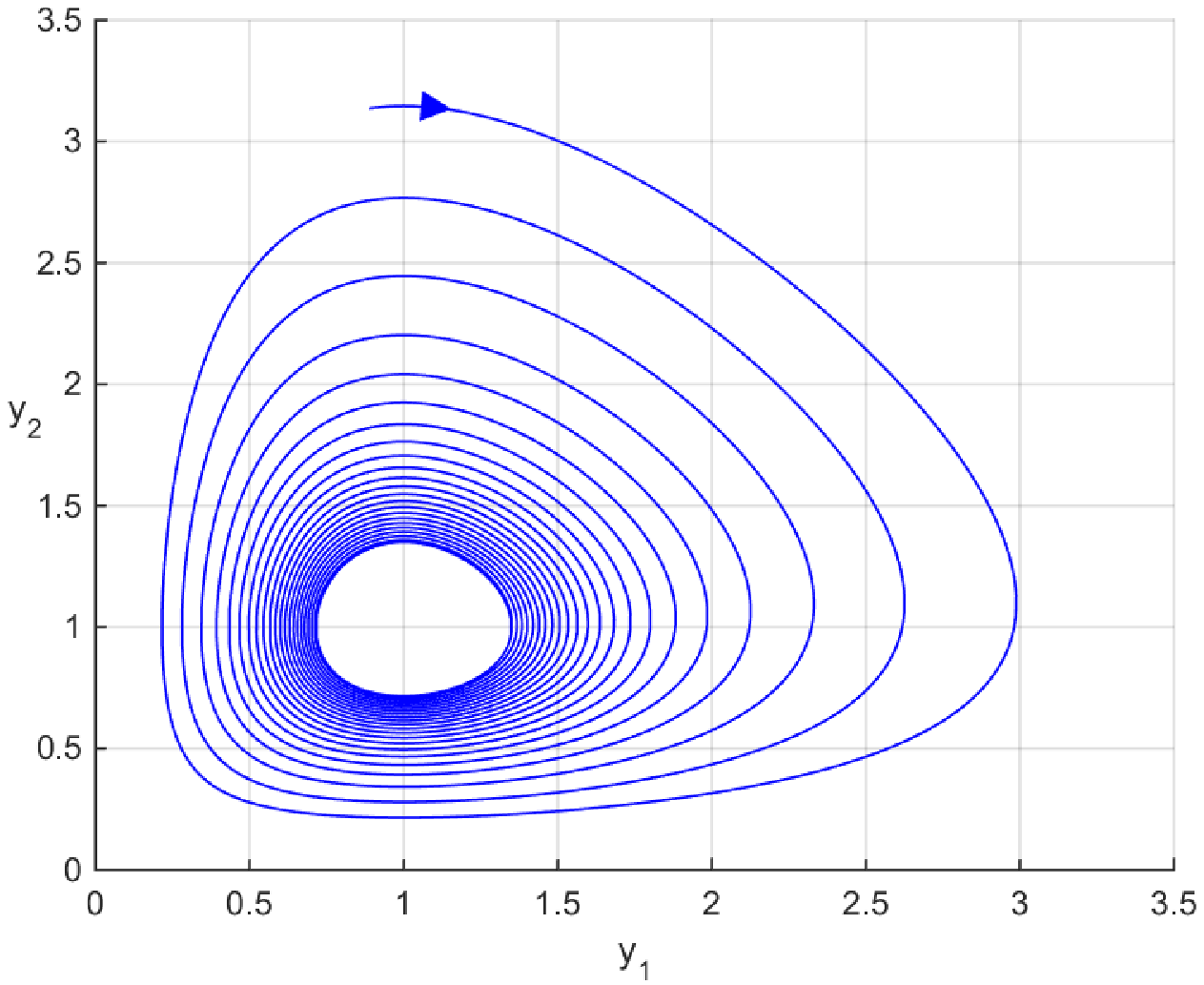}
\end{center}
\vspace{-.1in}
\ \ \ \ \ \ \ \ \ \ \ 
  \ \ \ \ \  \ \ \ \ \ \ \ \ \ \ \ \ \ \ \ \ \ Fig. 2:\ Near optimal state trajectory for $\epsilon =0.1$
  
  \begin{center}
\includegraphics[width=3.20in]{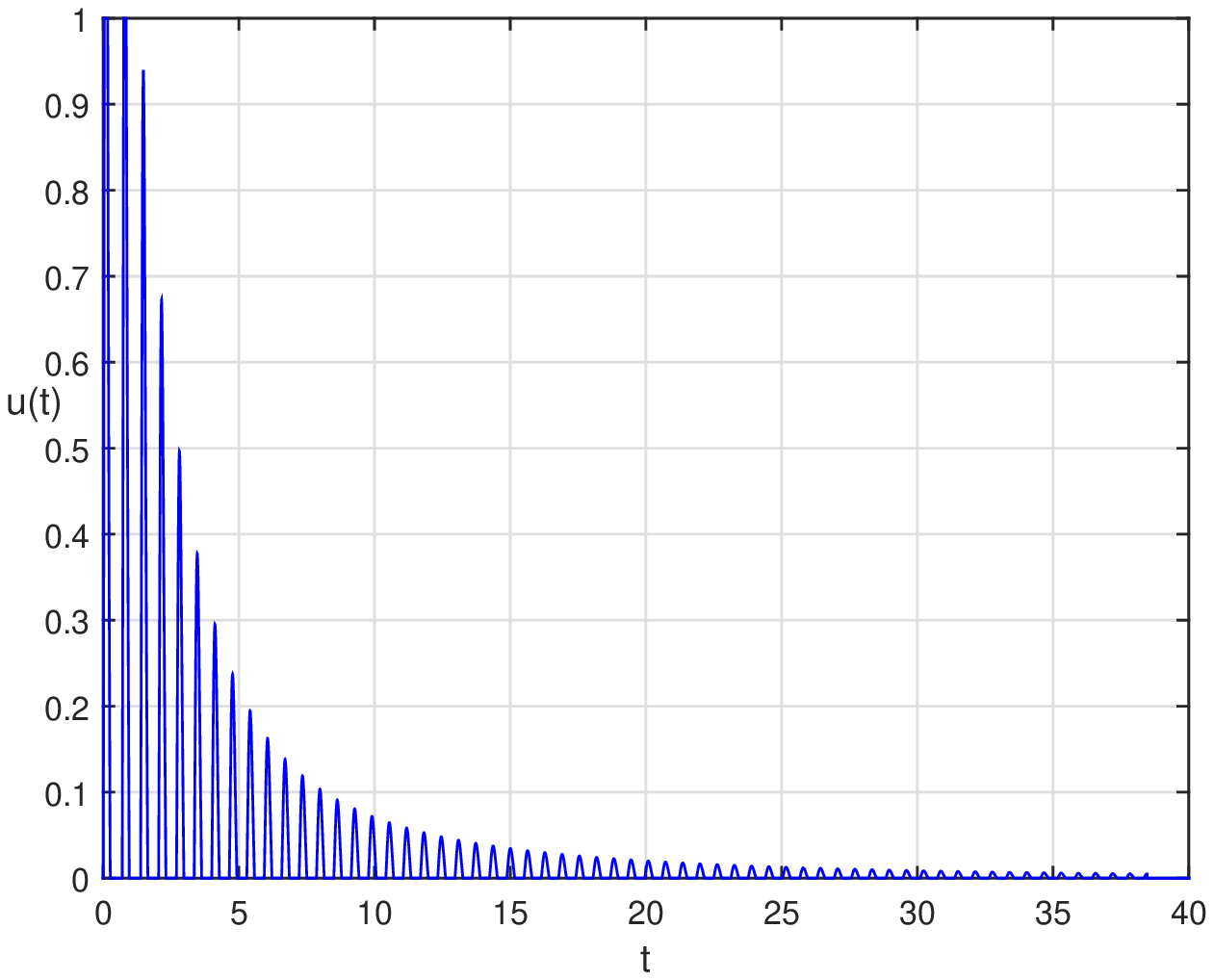}
\end{center}
\vspace{-.1in}
 \ \ \ \ \ \ \ \ \  \ \ \ \ \ \ \ \ \ \ \ \
\ \  \ \ \ \ \ \ \ \ \ \ \ \ \ Fig. 3:\  $u(t)=  u^{M,N}(y(t),F(y(t))) $ 

\begin{center}
\includegraphics[width=3.20in]{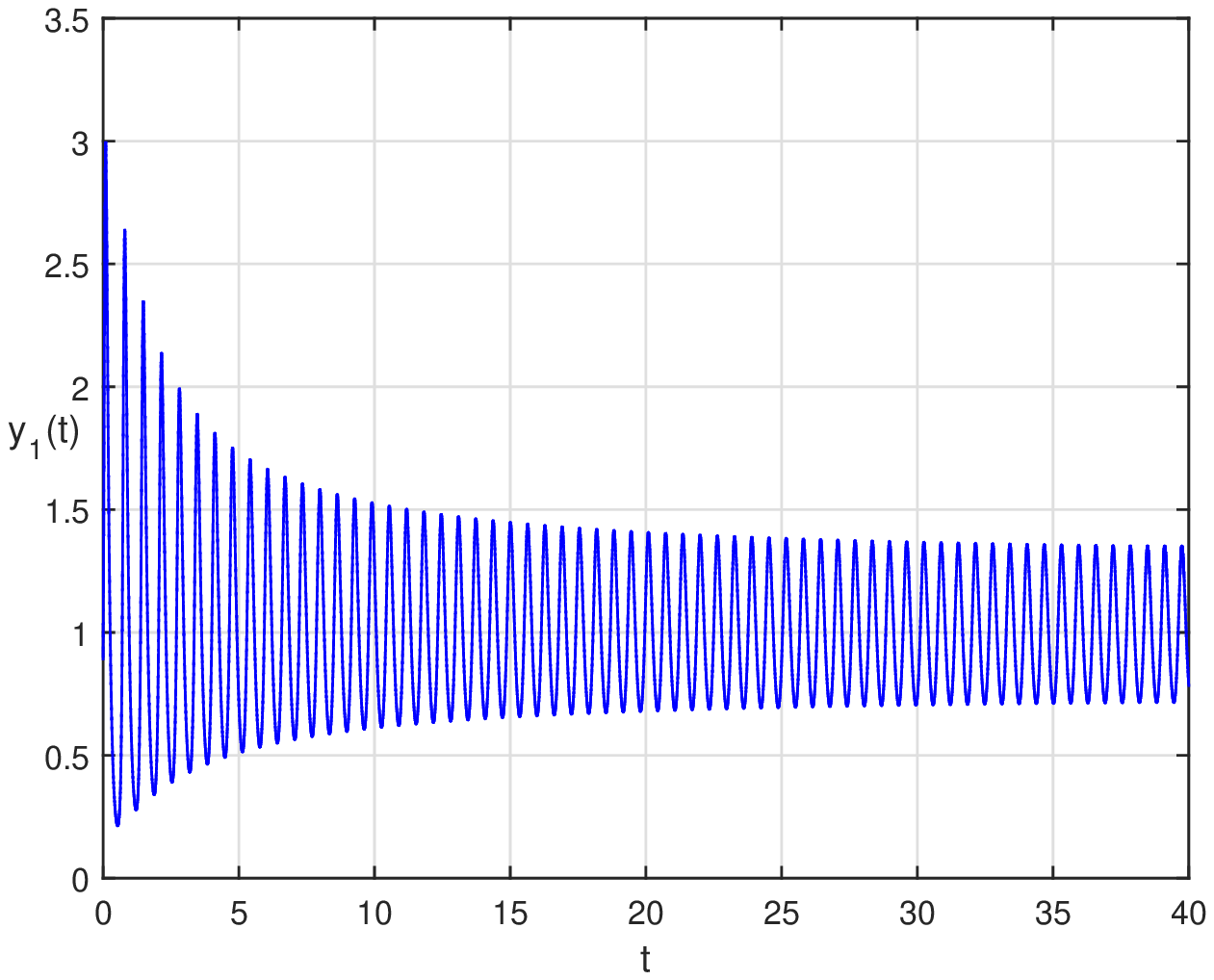}
\end{center}
\vspace{-.1in}
 \ \ \ \ \ \ 
\ \  \ \ \ \ \ \ \ \ \ \ \ \ \ 
  \ \ \ \ \  \ \ \ \ \ \ \ \ \ \ \ \ \ \ \ \ \ \ \ \ \ \ \ \ \ \ \ Fig. 4:\ $y_1(t)$ 

\bigskip

\section{Proofs of Theorems \ref{Prop-Averaging-1} and \ref{Prop-converse-convergence} and a proof of Proposition \ref{Prop-asym-opt-contr}}\label{Sec-Proofs-I}\

{\it Proof of Theorem \ref{Prop-Averaging-1}. }
Let
\begin{equation}\label{e-CSO-21-ext-0}
X(t,u,y) \BYDEF (h(u,y), e^{-Ct}r(u,y))
\end{equation}
and let
\begin{equation}\label{e-CSO-21-ext}
V(t,z)\BYDEF \cup_{\mu\in W(z)} \{\int_{U\times Y}X(t,u,y)\mu(du,dy) \}= \cup_{\mu\in W(z)}\{(\tilde h(\mu), e^{-Ct}\tilde r(\mu)) \}.
\end{equation}
Note that the map $V(t,z)$ is convex and compact valued (due to the fact that
$W(z): Z\leadsto \mathcal{P}(U\times Y) $  is convex and compact valued). Also from (\ref{e-CSO-21-1})
it follows that
 \begin{equation}\label{e-CSO-21-1-ext}
d_H(V(t',z'),V(t'', z''))\leq \hat L (|t'-t''| +||z'-z''||) \ \ \ \ \forall t',t''\geq 0, \ \forall z',z''\in Z,\ \ \ \ \          \hat L= {\rm const}.
\end{equation}
 Let $\theta_{\epsilon}(t)$ stand for  the solution of the equation
\begin{equation}\label{e-theta-1}
\frac{d\theta_{\epsilon} (t)}{dt} = e^{-Ct}r(u_{\epsilon}(t),y_{\epsilon}(t)), \ \ \ \ \ \ \ \ \theta_{\epsilon} (0)=0.
\end{equation}
Then
$x_{\epsilon}(t)\BYDEF (z_{\epsilon}(t),\theta_{\epsilon}(t)) $ is the solution of the equation
(from this point on, we will  write $u(t)$ instead of $u_{\epsilon}(t)$ to simplify the notations)
\begin{equation}\label{e-theta-1-11}
\frac{dx_{\epsilon} (t)}{dt} = X(t,u(t),y_{\epsilon}(t)), \ \ \ \ \ \ \ \ x_{\epsilon} (0)=(z_0,0).
\end{equation}

Due to  Filippov's selection theorem (see, e.g., Theorem 8.2.10 in \cite{aubin1}), to prove the estimates (\ref{e-CSO-22}) and (\ref{e:intro-0-3-8}), it is sufficient to show that
there exists a solution $x(t)= (z(t),\theta(t)) $ of the differential inclusion
\begin{equation}\label{e-CSO-18-1}
\frac{dx(t)}{dt} \in  V(t, z(t)), \ \ \ \ \ x(0)=(z_0, 0)
\end{equation}
such that
\begin{equation}\label{e-theta-2}
\max_{t\in [0,T]}||x_{\epsilon}(t) - x(t)||\leq \beta(\epsilon, T).
\end{equation}
From the validity of (\ref{e-CSO-22}) it will follow that $z(t)\in Z \ \forall t\in [0,T]$ (for $\epsilon$ small enough). Also  the solution $x(t)$ of (\ref{e-theta-2}) is extendable to the interval  $[T,\infty)$ in such a way that $z(t)\in Z \ \forall t\in [T,\infty)$ (due to the viability of the averaged system). Thus, the theorem will be proved if one establishes that a solution of (\ref{e-CSO-18-1}) satisfying (\ref{e-theta-2}) exists.

 Since $x_{\epsilon}(t)$ satisfies (\ref{e-theta-1-11}), one can write down the following relationships
\begin{equation}\label{e-CSO-23}
x_{\epsilon}(t_{l+1}) = x_{\epsilon}(t_{l}) + \int_{t_l}^{t_{l+1}}X(t, u(t), y_{\epsilon}(t))dt=x_{\epsilon}(t_{l}) +\epsilon
 \int_{\frac{t_l}{\epsilon}}^{\frac{t_{l+1}}{\epsilon}}X(\epsilon \tau,u(\epsilon \tau), y_{\epsilon}(\epsilon \tau))d\tau
 , \ \ \ \ \l=0,1,..., N_{\epsilon}-1,
\end{equation}
where $t_l $ are as defined in (\ref{e:contr-rev-100-1}) and $\ N_{\epsilon}\BYDEF \lfloor \frac{T}{\Delta(\epsilon)} \rfloor$ ($\lfloor\cdot\rfloor $ standing for the floor function).
Consider also the difference equation
\begin{equation}\label{e-CSO-24}
\bar x_{\epsilon}(l+1) = \bar x_{\epsilon}(l) + \Delta(\epsilon) \int_{U\times Y}X(t_l, u,y)\bar{\mu}_l(du,dy), \ \ \ \ l=0,1,..., N_{\epsilon}-1, \ \ \ \ \ \ \bar x_{\epsilon}(0)=(z_0, 0),
\end{equation}
where $\bar x_{\epsilon}(l) \BYDEF (\bar z_{\epsilon}(l), \bar \theta_{\epsilon}(l))$, and  $\bar{\mu}_l(du,dy) $ are constructed iteratively in such a way that
\begin{equation}\label{e-CSO-24-1}
\bar{\mu}_l\in W(\bar z_{\epsilon}(l))
\end{equation}
and
\begin{equation}\label{e-CSO-24-2}
||(\frac{\Delta(\epsilon)}{\epsilon})^{-1}\int_{\frac{t_l}{\epsilon}}^\frac{t_{l+1}}{\epsilon}X(t_l, u(\epsilon\tau), y_{\epsilon}(\epsilon\tau))d\tau
- \int_{U\times Y}X(t_l, u, y)\bar{\mu}_l(du,dy)|| \leq \hat L ||z_{\epsilon}(t_{l}) - \bar z_{\epsilon}(l) || + \beta_1(\epsilon)
\end{equation}
for any $\ l=0,1,..., N_{\epsilon}-1$ (with $\hat L$ being a Lipschitz constant from (\ref{e-CSO-21-1-ext}) and with $\ \lim_{\epsilon\rightarrow 0}\beta_1(\epsilon) = 0$; the construction of $\bar{\mu}_l, \ l=0,1,..., N_{\epsilon}-1, $ that satisfy (\ref{e-CSO-24-1}) and (\ref{e-CSO-24-2}) is described at the end of the proof of the theorem).

By subtracting (\ref{e-CSO-24}) from (\ref{e-CSO-23}), one obtains the inequalities
$$
||x_{\epsilon}(t_{l+1}) - \bar x_{\epsilon}(l+1)||\leq ||x_{\epsilon}(t_{l}) - \bar x_{\epsilon}(l)|| + \Delta(\epsilon) ( \hat L ||z_{\epsilon}(t_{l}) - \bar z_{\epsilon}(l) || + \beta_2(\epsilon))
$$
\begin{equation}\label{e-CSO-30}
\leq ||x_{\epsilon}(t_{l}) - \bar x_{\epsilon}(l)|| + \Delta(\epsilon)(\hat L||x_{\epsilon}(t_{l}) - \bar x_{\epsilon}(l) || + \beta_2(\epsilon)), \ \ \ \ \l=0,1,..., N_{\epsilon}-1,
\end{equation}
where $\ \lim_{\epsilon\rightarrow 0}\beta_2(\epsilon) = 0$.
 These imply (see Proposition 5.1 in \cite{Gai1}) that
\begin{equation}\label{e-CSO-31}
||x_{\epsilon}(t_{l}) - \bar x_{\epsilon}(l)|| \leq \kappa(\epsilon , T)
, \ \ \ \ \l=0,1,..., N_{\epsilon}-1, \ \ \ \ \ {\rm where} \ \ \ \ \ \lim_{\epsilon\rightarrow 0}\kappa(\epsilon , T) = 0.
\end{equation}
Define the piecewise linear function $\bar{\bar x}_{\epsilon}(t)=(\bar{\bar z}_{\epsilon}(t), \bar{\bar \theta}_{\epsilon}(t))$ on the interval $[0,T]$ by the equations
$$
\bar{\bar x}_{\epsilon}(t)= \bar x_{\epsilon}(l) + (t-t_l) v_l \ \ \forall t\in [t_l,t_{l+1}), \ \ l=0,1,..., N_{\epsilon}-1,
$$
\vspace{-.2in}
$$
\bar{\bar x}_{\epsilon}(t)= \bar x_{\epsilon}(l) + (t-t_{N_{\epsilon}}) v_{N_{\epsilon}-1} \ \ \forall t\in [t_{N_{\epsilon}},t_{N_{\epsilon}+1}], \ \ \ \
 \ \ \ t_{N_{\epsilon}+1}\BYDEF T,
$$
where
$$
v_l\BYDEF \int_{U\times Y}X(t_l,u,y)\bar{\mu}_l(du,dy)\in V(t_l, \bar z_{\epsilon}(l)).
$$
Note that
\begin{equation}\label{e-CSO-35}
max_{t\in [t_l,t_{l+1}]}|| \bar{\bar x}_{\epsilon}(t)- \bar x_{\epsilon}(l)||\leq\Delta(\epsilon)||v_l|| \leq M\Delta(\epsilon), \ \ \ \ \ {\rm where} \ \ \ \ \
M\BYDEF \max_{(u,y)\in U\times Y}(||h(u,y)||+|r(u,y)|).
\end{equation}
Also, due to (\ref{e-CSO-21-1-ext}), for $t\in (t_l, t_{l+1}) $ (with $l=0,1,...N_{\epsilon}$),
$$
d(\frac{d\bar{\bar x}_{\epsilon}(t)}{dt}, V(t,\bar{\bar z}_{\epsilon}(t)))\leq d(v_l, V(t_l ,\bar z_{\epsilon}(l)))+ d_H(V(t_l,\bar z_{\epsilon}(l)), V(t, \bar{\bar z}_{\epsilon}(t)))
$$
$$
= d_H(V(t_l, \bar z_{\epsilon}(l)), V(t,\bar{\bar z}_{\epsilon}(t)))\leq L_1 \Delta(\epsilon), \ \ \ \ \ L_1={\rm const}.
$$
From Filippov existence theorem (see, e.g. Theorem 10.4.1 in \cite{aubin1}) it now follows that there exists a solution $x(t)$ of the differential inclusion
(\ref{e-CSO-18-1}) such that
\begin{equation}\label{e-CSO-36}
max_{t\in [0,T]}||\bar{\bar x}_{\epsilon}(t)- x(t)||\leq L_T\Delta(\epsilon), \ \ \ \ \ L_T={\rm const}.
\end{equation}
This (along with (\ref{e-CSO-31}) and (\ref{e-CSO-35})) prove the required statement.

To finalize the proof of the theorem, let us now describe  the construction of $\bar{\mu}_l$  ($l=0,1,..., N_{\epsilon}-1 $)  satisfying (\ref{e-CSO-24-1}) and (\ref{e-CSO-24-2}).
For $\ \tau\in [\frac{t_l}{\epsilon},\frac{t_{l+1}}{\epsilon}]  $,\ \  let $y_{y_{\epsilon}(t_l)}(\tau) \ $ be the solution of the system (\ref{e-CSO-2-SC-I})   that is obtained with the use of the control $u(\epsilon \tau) $
and that satisfies  the initial condition $\ y_{y_{\epsilon}(t_l)}(\frac{t_l}{\epsilon})= y_{\epsilon}(t_l) $. That is,
\begin{equation}\label{e-CSO-25}
y_{y_{\epsilon}(t_l)}(\tau)= y_{\epsilon}(t_l) + \int_{\frac{t_l}{\epsilon}}^{\tau}f(u(\epsilon\tau),y_{y_{\epsilon}(t_l)}(\tau'))d\tau' \ \ \ \ \ \ \ \ \ \forall \tau\in [\frac{t_l}{\epsilon},\frac{t_{l+1}}{\epsilon}].
\end{equation}
Due to the fact that $y_{\epsilon}(t) $ is a solution of (\ref{e-CSO}), it satisfies the equation
\begin{equation}\label{e-CSO-26}
y_{\epsilon}(\epsilon\tau)= y_{\epsilon}(t_l) + \int_{\frac{t_l}{\epsilon}}^{\tau}[f(u(\epsilon\tau), y_{\epsilon}(\epsilon\tau'))+ \epsilon g(u(\epsilon\tau'), y_{\epsilon}(\epsilon\tau'))]d\tau'\ \ \ \ \ \ \ \ \ \forall \tau\in [\frac{t_l}{\epsilon},\frac{t_{l+1}}{\epsilon}].
\end{equation}
Hence,
\begin{equation}\label{e-CSO-27}
||y_{\epsilon}(\epsilon\tau)- y_{y_{\epsilon}(t_l)}(\tau)||\leq
L_f \int_{\frac{t_l}{\epsilon}}^{\tau}||y_{\epsilon}(\epsilon\tau')- y_{y_{\epsilon}(t_l)}(\tau')||d\tau'+M_g\Delta(\epsilon)
\ \ \ \ \ \ \ \forall \tau\in [\frac{t_l}{\epsilon},\frac{t_{l+1}}{\epsilon}],
\end{equation}
where $L_f$ is a Lipschitz constant of $f(\cdot)$ and $\ M_g\BYDEF \max_{(u,y)\in U\times Y}||g(u,y)||$. By Gronwall-Bellman lemma, the latter implies
that
\begin{equation}\label{e-CSO-28}
\max_{\tau\in [\frac{t_l}{\epsilon},\frac{t_{l+1}}{\epsilon}]}||y_{\epsilon}(\epsilon\tau)- y_{y_{\epsilon}(t_l)}(\tau)||\leq
M_g\Delta(\epsilon)e^{L_f\frac{\Delta(\epsilon)}{\epsilon}}\leq M_g(\frac{\epsilon}{2L_f}\ln\frac{1}{\epsilon})(\frac{1}{\epsilon^{\frac{1}{2}}})< \epsilon^{\frac{1}{4}},
\end{equation}
the second to last inequality being valid due to (\ref{e-CSO-29}) ($\epsilon$ is assumed to be small enough).
 Note that from (\ref{e-CSO-28}) it follows that
\begin{equation}\label{e-CSO-29-1}
||(\frac{\Delta(\epsilon)}{\epsilon})^{-1}\int_{\frac{t_l}{\epsilon}}^\frac{t_{l+1}}{\epsilon}X(t_l, u(\epsilon\tau), y_{\epsilon}(\epsilon\tau))d\tau
- (\frac{\Delta(\epsilon)}{\epsilon})^{-1}\int_{\frac{t_l}{\epsilon}}^\frac{t_{l+1}}{\epsilon}X(t_l, u(\epsilon\tau), y_{y_{\epsilon}(t_l)}(\tau))d\tau ||\leq
\hat L\epsilon^{\frac{1}{4}}.
\end{equation}
Denote by $\mu_l$ the occupational measure generated by the pair $(u(\epsilon\tau), y_{y_{\epsilon}(t_l)}(\tau))$ on the interval
$[\frac{t_l}{\epsilon},\frac{t_{l+1}}{\epsilon}]$. That is, $\mu_l$ is such that for any continuous function $q(u,y)$
\begin{equation}\label{Proof-1-extra-1}
\int_{U\times Y}q(u,y)\mu_l(du,dy) = (\frac{\Delta(\epsilon)}{\epsilon})^{-1}\int_{\frac{t_l}{\epsilon}}^\frac{t_{l+1}}{\epsilon}q(u(\epsilon\tau), y_{y_{\epsilon}(t_l)}(\tau))d\tau .
\end{equation}
By Proposition \ref{Prop-inv-measure-set-2} (see (\ref{e-CSO-11})), there exists $\tilde{\mu}_l\in W(z_{\epsilon}(t_l))$ such that
\begin{equation}\label{Proof-1-extra-2}
\rho(\mu_l, \tilde{\mu}_l)\leq \beta ((\frac{\Delta(\epsilon)}{\epsilon})^{-1})\BYDEF \beta_3(\epsilon), \ \ \ \ \ \ \lim_{\epsilon\rightarrow 0}\beta_3(\epsilon) = 0,
\end{equation}
the latter implying that
\begin{equation}\label{e-CSO-32}
||(\frac{\Delta(\epsilon)}{\epsilon})^{-1}\int_{\frac{t_l}{\epsilon}}^\frac{t_{l+1}}{\epsilon}X(t_l, u(\epsilon\tau), y_{y_{\epsilon}(t_l)}(\tau))d\tau
- \int_{U\times Y}X(t_l, u, y)\tilde{\mu}_l(du,dy)|| \leq \beta_4(\epsilon), \ \ \ \ \ \ \lim_{\epsilon\rightarrow 0}\beta_4(\epsilon) = 0.
\end{equation}
Also, due to  (\ref{e-CSO-21-1-ext}), there exists $\bar{\mu}_l\in W(\bar z_{\epsilon}(l)) $ such that
\begin{equation}\label{e-CSO-32-1}
||\int_{U\times Y}X(t_l,u, y)\tilde{\mu}_l(du,dy) - \int_{U\times Y}X(t_l, u, y)\bar{\mu}_l(du,dy)||\leq \hat L ||z_{\epsilon}(t_l) - \bar z_{\epsilon}(l)||.
\end{equation}
Summarizing (\ref{e-CSO-29-1}), (\ref{e-CSO-32}) and (\ref{e-CSO-32-1}), one obtains that thus defined $\bar{\mu}_l $ satisfies (\ref{e-CSO-24-2}) (with  $\beta_1(\epsilon)= \hat L\epsilon^{\frac{1}{4}} +\beta_4(\epsilon) $).
$\ \Box$


 {\it Proof of Theorem \ref{Prop-converse-convergence}.}
 Let  $(\mu(\cdot), z(\cdot)) $ satisfy (\ref{e-CSO-17})-(\ref{e-CSO-18}) and let $\theta(t)$ be the solution of the equation
 \begin{equation}\label{e-theta-12}
\frac{d\theta (t)}{dt} = e^{-Ct}\tilde r(\mu(t)), \ \ \ \ \ \ \ \ \theta (0)=0.
\end{equation}
  Then $x(\cdot)= (z(t),\theta (t) )$ is a solution of the differential inclusion (\ref{e-CSO-18-1}) and, due to (\ref{e-CSO-21-1-ext}),
  for any $t\in [t_l,t_{l+1}] $ ($t_l, \  l=0,...,N_{\epsilon}-1$, are as in (\ref{e:contr-rev-100-1})),
 \begin{equation}\label{eq-Proof2-1}
 \frac{dx(t)}{dt} \in  V(t, z(t))\subset V(t_l ,z(t_l)) + \hat L (|t-t_l| +||z(t)-z(t_l)||) \bar B \subset V(t_l, z(t_l))+L_2\Delta(\epsilon)\bar B ,
  \end{equation}
   where $L_2$ is a sufficiently large constant and $\bar B $ is the closed unit ball in $\R^{k+1}$. Consequently,
  $$
\Delta(\epsilon)^{-1} \int_{t_l}^{t_{l+1}}\frac{dx(t)}{dt}dt   \in  V(t_l, z(t_l)) +L_2 \Delta(\epsilon)\bar B, \ \  \ \ l=0,...,N_{\epsilon}-1
$$
 \vspace{-.2in}
 \begin{equation}\label{eq-Proof2-3}
\Rightarrow \ \ \ \ \ dist \left(\Delta(\epsilon)^{-1} \int_{t_l}^{t_{l+1}}\frac{dx(t)}{dt}dt , V(t_l, z(t_l))\right)\leq L_2 \Delta(\epsilon), \ \  \ \ l=0,...,N_{\epsilon}-1.
 \end{equation}
 Let $\bar v_l\in V(t_l, z(t_l)), \ l=0,...,N_{\epsilon}-1$, be defined by the equation
 \begin{equation}\label{eq-Proof2-4}
\bar v_l \BYDEF {\rm argmin} \{ || v - \Delta(\epsilon)^{-1} \int_{t_l}^{t_{l+1}}\frac{dx(t)}{dt}dt ||\ | \ v\in V(t_l , z(t_l)) \}.
\end{equation}
 Note that, by (\ref{eq-Proof2-3}),
 \begin{equation}\label{eq-Proof2-5}
 || \bar v_l - \Delta(\epsilon)^{-1} \int_{t_l}^{t_{l+1}}\frac{dx(t)}{dt}dt ||\leq  L_2\Delta(\epsilon), \ \ \ \ l=0,...,N_{\epsilon}-1 ,
\end{equation}
and note  that from the definition of the maltivalued map $V(\cdot, \cdot) $ (see (\ref{e-CSO-21-ext})) it follows that
there exists $\bar \mu_l\in W(z(t_l)) $ such that
\begin{equation}\label{eq-Proof2-6}
 \bar v_l = \int_{U\times Y }X(t_l, u,y)\bar \mu_l(du,dy), \ \ \ \ l=0,...,N_{\epsilon}-1.
\end{equation}
Consider the difference equation
\begin{equation}\label{e-CSO-24-1-1}
\bar x_{\epsilon}(l+1) = \bar x_{\epsilon}(l) + \Delta(\epsilon) \bar v_l, \ \ \ \ l=0,1,..., N_{\epsilon}-1, \ \ \ \ \ \ \bar x_{\epsilon}(0)=(z_0, 0),
\end{equation}
where $\bar x_{\epsilon}(l) \BYDEF (\bar z_{\epsilon}(l), \bar \theta_{\epsilon}(l))$.
Subtracting the latter from
\begin{equation}\label{e-CSO-24-1-2}
 x(t_{l+1}) =  x(t_{l}) + \Delta(\epsilon) \left(\Delta(\epsilon)^{-1} \int_{t_l}^{t_{l+1}}\frac{dx(t)}{dt}dt\right), \ \ \ \ l=0,1,..., N_{\epsilon}-1, \ \ \ \ \ \ \bar z_{\epsilon}(0)=z_0,
\end{equation}
one can establish the validity of the inequalities
\begin{equation}\label{e-CSO-24-1-3}
|| x(t_{l+1})- \bar x_{\epsilon}(l+1)||\leq ||  x(t_{l})- \bar x_{\epsilon}(l)|| +
  L_2\Delta^2(\epsilon) , \ \ \ \ l=0,1,..., N_{\epsilon}-1.
\end{equation}
These imply (see Proposition 5.1 in \cite{Gai1})
\begin{equation}\label{e-CSO-24-1-4}
||x(t_{l}) - \bar x_{\epsilon}(l)|| \leq \kappa_1(\epsilon , T)
, \ \ \ \ \l=0,1,..., N_{\epsilon}-1, \ \ \ \ \ {\rm where} \ \ \ \ \ \lim_{\epsilon\rightarrow 0}\kappa_1(\epsilon , T) = 0.
\end{equation}
Assume now that the control $u(t)$ has been defined on the interval $[0,t_l], \ l<N_{\epsilon}$, and extend its definition to the interval
 $[0,t_{l+1}]$. Denote by $(y_{\epsilon}(t),x_{\epsilon}(t))  $ the solution of the system (\ref{e-CSO}),(\ref{e-theta-1-11}) obtained
 with the use of the control $u(t)$ on the interval $[0,t_l] $. Note that, due to (\ref{e-CSO-21-1-ext}), there exists
 a measure  $\tilde \mu_l\in W(z_{\epsilon}(t_l))$ such that (see (\ref{eq-Proof2-6}))
$$
 ||\int_{U\times Y }X(t_l, u,y)\tilde \mu_l(du,dy)- \bar v_l|| = ||\int_{U\times Y }X(t_l, u,y)\tilde \mu_l(du,dy) -
 \int_{U\times Y }X(t_l,u,y) \bar \mu_l(du,dy) ||
 $$
 \vspace{-.2in}
  \begin{equation}\label{eq-Proof2-8}
 \ \ \ \ \ \ \ \ \ \ \ \ \ \ \ \ \ \ \ \ \ \ \ \ \ \ \ \ \ \leq \hat L||z_{\epsilon}(t_l)-z(t_l)||.
\end{equation}
  Define the control $u(t) $ on the interval $[t_l , t_{l+1}]  $ in such a way that the pair $(u(\epsilon\tau), y_{y_{\epsilon}(t_l)}(\tau)) $ (where
  $y_{y_{\epsilon}(t_l)}(\tau) $ is the solution of the reduced system (\ref{e-CSO-2-SC-I}) on the interval $[\frac{t_l}{\epsilon} , \frac{t_{l+1}}{\epsilon}]$ which is obtained  with the use of the control $u(\epsilon\tau) $ and with the initial condition
  $y_{y_{\epsilon}(t_l)}(\frac{t_l}{\epsilon})= y_{\epsilon}(t_l)  $)
    generates  an occupational measure $\mu_l $ satisfying (\ref{Proof-1-extra-2}) and (\ref{e-CSO-32}) (the existence of such control follows from
 Assumption \ref{Ass-existence-LOMS}). Use
  the control $u(t)$ to obtain the solution $(y_{\epsilon}(t),x_{\epsilon}(t))  $ of the system (\ref{e-CSO}),(\ref{e-theta-1-11}) on the interval $[t_l, t_{l+1}] $. Note that, similarly to the proof of Theorem \ref{Prop-Averaging-1},
 it can be verified that the estimate (\ref{e-CSO-29-1}) is valid, which along with (\ref{e-CSO-32}) and (\ref{eq-Proof2-8}) imply that
 \begin{equation}\label{eq-Proof2-8-1}
||(\frac{\Delta(\epsilon)}{\epsilon})^{-1}\int_{\frac{t_l}{\epsilon}}^\frac{t_{l+1}}{\epsilon}X(t_l, u(\epsilon\tau), y_{\epsilon}(\epsilon\tau))d\tau
- \bar v_l||\leq \hat L||z_{\epsilon}(t_l)-z(t_l)|| + \beta_1(\epsilon),
\end{equation}
where $\lim_{\epsilon\rightarrow 0}\beta_1(\epsilon) = 0 $.

  By continuing this process, one can construct the control $u(t)$ on the interval
   $[0,t_{N_{\epsilon}}]$. On the interval $[t_{N_{\epsilon}},T]$, the control $u(t)$ can be taken to be equal to an arbitrary
 $u\in U$. Let us show that the solution $(y_{\epsilon}(t),x_{\epsilon}(t))  $ of the system (\ref{e-CSO}),(\ref{e-theta-1-11}) obtained with this control satisfies  (\ref{e-CSO-22}) and (\ref{e:intro-0-3-8}).

Having in mind the validity of (\ref{eq-Proof2-8-1}) and arguing as in the proof of Theorem \ref{Prop-Averaging-1}, one can establish the validity of the estimates (\ref{e-CSO-31}), where $\{\bar x_{\epsilon}(l)\} $ is the solution of (\ref{e-CSO-24-1-1}). These  and (\ref{e-CSO-24-1-4}) imply that
\begin{equation}\label{e-CSO-24-1-4-1}
||x_{\epsilon}(t_{l}) - x(t_{l})|| \leq  \kappa_2(\epsilon , T)
, \ \ \ \ \l=0,1,..., N_{\epsilon}-1, \ \ \ \ \ {\rm where} \ \ \ \ \ \lim_{\epsilon\rightarrow 0}\kappa_2(\epsilon , T) = 0.
\end{equation}
Since
$$
\max_{t\in [t_l,t_{l+1}]}||x_{\epsilon}(t)-x_{\epsilon}(t_l)||\leq M\Delta(\epsilon), \ \ \ \ \ \ \
\max_{t\in [t_l,t_{l+1}]}||x(t)-x(t_l)||\leq M\Delta(\epsilon),
$$
the estimates (\ref{e-CSO-24-1-4-1}) imply the validity of (\ref{e-theta-2}) (with $\beta (\epsilon , T) \BYDEF  \kappa_2(\epsilon , T) + 2 M\Delta(\epsilon) $), which in turn implies
(\ref{e-CSO-22}) and (\ref{e:intro-0-3-8}) . Due to (\ref{e-CSO-22}) and (\ref{e-CSO-28})  and due to (\ref{e-CSO-20}), the inclusions (\ref{e-CSO-20-101}) are valid for $t\in [0,T] $. Since
the augmented perturbed system (\ref{e-CSO}), (\ref{e-CSO-0-2}) is viable in $Y\times Z$, the solution of the latter can be extended in such a way that the inclusions (\ref{e-CSO-20-101}) will be satisfied for all $t\in [0,\infty) $. This completes the proof of the theorem.
 $\ \Box$

{\it Proof of Proposition \ref{Prop-asym-opt-contr}}. Let $y_{\epsilon}^*(t)$ be the solution of  the system
(\ref{e-CSO-101}) obtained with the control $u_{\epsilon}^*(t) $ and let $x_{\epsilon}^*(t)\BYDEF ( z_{\epsilon}^*(t), \theta_{\epsilon}^*(t) ) $ be defined by the equation
\begin{equation}\label{e-theta-1-11-1}
\frac{dx_{\epsilon}^* (t)}{dt} = X(t,u_{\epsilon}^*(t),y_{\epsilon}^*(t)), \ \ \ \ \ \ \ \ x_{\epsilon}^* (0)=(z_0,0),
\end{equation}
where $X(\cdot)$ is as in (\ref{e-CSO-21-ext-0}). Note that a part of (\ref{e-theta-1-11-1}) is the equation
\begin{equation}\label{e-theta-1-11-1-11}
\frac{d\theta_{\epsilon}^* (t)}{dt} = e^{-ct}r(u_{\epsilon}^*(t),y_{\epsilon}^*(t)), \ \ \ \ \ \ \ \ \theta_{\epsilon}^* (0)=0.
\end{equation}
 From (\ref{e-theta-1-11-1}) it follows that
$$
x_{\epsilon}^*(t_{l+1}) = x_{\epsilon}^*(t_{l}) + \int_{t_l}^{t_{l+1}}X(t, u^*_{\epsilon}(t), y_{\epsilon}^*(t))dt
= x_{\epsilon}^*(t_{l}) + \int_{t_l}^{t_{l+1}}X(t , u^*(y_{y_{\epsilon}^*(t_l)}(\frac{t}{\epsilon}),z_{\epsilon}^*(t_l)), y_{\epsilon}^*(t))dt
$$
\vspace{-.2in}
\begin{equation}\label{e-CSO-23-11}
=x_{\epsilon}^*(t_{l}) +\epsilon
 \int_{\frac{t_l}{\epsilon}}^{\frac{t_{l+1}}{\epsilon}}X(\epsilon\tau , u^*(y_{y_{\epsilon}^*(t_l)}(\tau),z_{\epsilon}^*(t_l)), y_{\epsilon}^*(\epsilon \tau))d\tau .
\end{equation}
Similarly to (\ref{e-CSO-28}), one can verify that
\begin{equation}\label{e-CSO-28-11}
\max_{\tau\in [\frac{t_l}{\epsilon},\frac{t_{l+1}}{\epsilon}]}||y_{\epsilon}^*(\epsilon\tau)- y_{y_{\epsilon}^*(t_l)}(\tau)||< \epsilon^{\frac{1}{4}}.
\end{equation}
Consequently,
\begin{equation}\label{e-CSO-28-12}
||\frac{\epsilon}{\Delta(\epsilon)}\int_{\frac{t_l}{\epsilon}}^{\frac{t_{l+1}}{\epsilon}}X(t_l, u^*(y_{y_{\epsilon}^*(t_l)}(\tau),z_{\epsilon}^*(t_l)), y_{\epsilon}^*(\epsilon \tau))d\tau
-
\frac{\epsilon}{\Delta(\epsilon)} \int_{\frac{t_l}{\epsilon}}^{\frac{t_{l+1}}{\epsilon}}X(t_l, u^*(y_{y_{\epsilon}^*(t_l)}(\tau),z_{\epsilon}^*(t_l)), y_{y_{\epsilon}^*(t_l)}(\tau))d\tau||
\end{equation}
\vspace{-.3in}
$$
\ \ \ \ \ \ \ \ \ \ \ \ \ \ \ \ \ \ \ \ \ \ \ \ \ \ \ \ \ \ \leq L_X \epsilon^{\frac{1}{4}},
$$
where $L_X$ is a Lipschitz constant of $X(t,u,y)$ in $y$. Also, due to periodicity of $y_{y_{\epsilon}^*(t_l)}(\tau) $, one has
(see (\ref{e-opt-OM-1-2}))
\begin{equation}\label{e-CSO-28-13}
||\frac{\epsilon}{\Delta(\epsilon)} \int_{\frac{t_l}{\epsilon}}^{\frac{t_{l+1}}{\epsilon}}X(t_l, u^*(y_{y_{\epsilon}^*(t_l)}(\tau),z_{\epsilon}^*(t_l)), y_{y_{\epsilon}^*(t_l)}(\tau))d\tau
- \tilde X(t_l, z_{\epsilon}^*(t_l)) ||\leq
\bar L \frac{\epsilon  }{\Delta(\epsilon)},
\end{equation}
where $\bar L $ is a positive constant and where
\begin{equation}\label{e-opt-OM-1-2-101}
\tilde X(t, z)\BYDEF \frac{1}{T_z}\int_0^{T_z}X(t,  u^*(y_z(\tau),z),y_z(\tau)) d\tau .
\end{equation}
Note that, by (\ref{e-opt-OM-1-2}) and (\ref{e-opt-OM-1-2-11}),
\begin{equation}\label{e-opt-OM-1-2-102}
\tilde X(t, z) = (\tilde h^*(z), e^{-Ct}\tilde r^*(z)) .
\end{equation}
Let $\theta^* (t)$ be the solution of the equation
\begin{equation}\label{e-opt-OM-1-2-103}
\frac{d\theta (t)}{dt} = e^{-Ct}\tilde r^*(z^*(t)), \ \ \ \ \ \ \ \ \theta (0)=0
\end{equation}
and let $x(t)\BYDEF (z^*(t),\theta^*(t))$, where $z^*(t)$  is the solution of (\ref{e-opt-OM-1-1}).
Then
\begin{equation}\label{e-theta-1-11-11}
\frac{dx (t)}{dt} = \tilde X(t,z^*(t)), \ \ \ \ \ \ \ \ x(0)=(z_0,0).
\end{equation}
Consider the difference equation
\begin{equation}\label{e-CSO-24-11}
 x_{\epsilon}(l+1) =  x_{\epsilon}(l) + \Delta(\epsilon)\tilde X (t_l, z_{\epsilon}(l)),  \ \ \ \ x_{\epsilon}(0)=(z_0,0),
\end{equation}
where $x_{\epsilon}(l)\BYDEF (z_{\epsilon}(l),\theta_{\epsilon}(l)) $.
From the assumed Lipschitz continuity of $\tilde h(z) $ and $\tilde r(z) $ it follows that $\tilde X(t,z) $ is Lipschitz
continuous, which, as can be readily verified, implies the validity of the estimates
\begin{equation}\label{e-CSO-24-12}
 ||x(t_l) - x_{\epsilon}(l)|| \leq \sigma_1(\epsilon , T), \ \ \ \ l=0,1,...,  \lfloor \frac{T}{\Delta(\epsilon)} \rfloor, \ \ \ \ \ \ \ {\rm where}
 \ \ \ \ \ \ \ \lim_{\epsilon\rightarrow 0}\sigma_1(\epsilon , T) = 0,
\end{equation}
where $x(t)$ is the solution of (\ref{e-theta-1-11-11}).
By subtracting (\ref{e-CSO-24-11}) from (\ref{e-CSO-23-11}) and taking into account (\ref{e-CSO-28-12}) and (\ref{e-CSO-28-13}),
one can obtain
$$
||x_{\epsilon}^*(t_{l+1}) - x_{\epsilon}(l+1)||\leq ||x_{\epsilon}^*(t_{l}) - x_{\epsilon}(l)||
+ \Delta(\epsilon)\ || \tilde X (t_l, z_{\epsilon}^*(t_{l})) - \tilde X (t_l,, z_{\epsilon}(l))|| +\bar L_1 \Delta (\epsilon)\ (\epsilon^{\frac{1}{4}} +  \frac{\epsilon  }{\Delta(\epsilon)}+ \Delta(\epsilon))
$$
\vspace{-.2in}
$$
\leq ||x_{\epsilon}^*(t_{l}) - x_{\epsilon}(l)||
+ \Delta(\epsilon)\ L_{\tilde X}||z_{\epsilon}^*(t_{l}) - z_{\epsilon}(l)|| +\bar L_1 \Delta (\epsilon)\ (\epsilon^{\frac{1}{4}} +  \frac{\epsilon  }{\Delta(\epsilon)}+ \Delta(\epsilon)),
$$
\vspace{-.2in}
$$
\leq ||x_{\epsilon}^*(t_{l}) - x_{\epsilon}(l)||
+ \Delta(\epsilon)\ L_{\tilde X}||x_{\epsilon}^*(t_{l}) - x_{\epsilon}(l)|| +\bar L_1 \Delta (\epsilon)\ (\epsilon^{\frac{1}{4}} +  \frac{\epsilon  }{\Delta(\epsilon)}+ \Delta(\epsilon)),
$$
where $L_{\tilde X} $ is a Lipschitz constant of $\tilde X(t,\cdot)$ and $\bar L_2$ is sufficiently large. The latter  implies that
  \begin{equation}\label{e-extra-11}
 ||x_{\epsilon}^*(t_{l}) - x_{\epsilon}(l)|| \leq \sigma_2(\epsilon , T), \ \ \ \ l=0,1,..., \lfloor \frac{T}{\Delta(\epsilon)} \rfloor, \ \ \ \ \ \ \ {\rm where}
 \ \ \ \ \ \ \ \lim_{\epsilon\rightarrow 0}\sigma_2(\epsilon , T) = 0,
\end{equation}
    (see Proposition 5.1 in \cite{Gai1}). This along with (\ref{e-CSO-24-12}) implies that
      $$
 ||x_{\epsilon}^*(t_{l}) - x(t_l)|| \leq \sigma_1(\epsilon , T) + \sigma_2(\epsilon , T), \ \ \ \ l=0,1,..., \lfloor \frac{T}{\Delta(\epsilon)} \rfloor,
$$
which, in turn, leads to the estimate
\begin{equation}\label{e-CSO-24-15}
 \max_{t\in [0,T]}||x_{\epsilon}^*(t) - x(t)|| \leq \sigma_3(\epsilon , T), \ \ \ \ \ \ \ {\rm where}
 \ \ \ \ \ \ \ \lim_{\epsilon\rightarrow 0}\sigma_3(\epsilon , T) = 0.
\end{equation}
From (\ref{e-CSO-24-15}) it follows (see (\ref{e-theta-1-11-1-11}) and (\ref{e-opt-OM-1-2-103})) that
\begin{equation}\label{e-CSO-24-15-1}
 \left|\int_0^T e^{-ct}r(u_{\epsilon}^*(t) ,y_{\epsilon}^*(t))dt  - \int_0^T e^{-ct}\tilde r^*(z^*(t))dt\right|
 \leq \sigma_3(\epsilon , T).
\end{equation}
Using an argument similar to that used in Corollary \ref{Corollary-inequality-values}, one can verify that
(\ref{e-CSO-24-15-1}) implies.
\begin{equation}\label{e-CSO-24-15-2}
 \left|\int_0^{\infty} e^{-ct}r(u_{\epsilon}^*(t) ,y_{\epsilon}^*(t))dt  - \int_0^{\infty} e^{-ct}\tilde r^*(z^*(t))dt\right|
 \leq \sigma_4(\epsilon ), \ \ \ \ \ \ \ {\rm where}
 \ \ \ \ \ \ \ \lim_{\epsilon\rightarrow 0}\sigma_4(\epsilon ) = 0.
\end{equation}
By  definition
\begin{equation}\label{e-CSO-24-15-3}
\int_0^{\infty} e^{-ct}\tilde r^*(z^*(t))dt = \tilde R^*(y_0)
\end{equation}
if $u^*(y,z) $ generates an optimal ACG family and
\begin{equation}\label{e-CSO-24-15-3-101}
\int_0^{\infty} e^{-ct}\tilde r^*(z^*(t))dt \leq \tilde R^*(y_0) + \alpha
\end{equation}
if $u^*(y,z) $ generates an $\alpha $-near optimal ACG family (see Definitions  \ref{Def-ACG-opt} and \ref{Def-feedback-control}).
These along with (\ref{e-CSO-24-15-2})  prove the proposition. 
$\ \Box$

{\it Proof of Proposition  \ref{Prop-asym-opt-contr-simple}}.
 Let $ x_{\epsilon}^* (t)\BYDEF  (y_{\epsilon}^* (t),  \theta_{\epsilon}^* (t) )$ be the solution of the equation
\begin{equation}\label{e-theta-1-11-1-100}
\frac{d  x_{\epsilon}^*  (t)}{dt} = X(t,u^*( y_{\epsilon}^* (t) , F( y_{\epsilon}^* (t))), y_{\epsilon}^* (t)), \ \ \ \ \ \ \ \ x_{\epsilon}^*  (0)=(z_0,0),
\end{equation}
where $X(\cdot)$ is as in (\ref{e-CSO-21-ext-0}). Note that a part of (\ref{e-theta-1-11-1-100}) is the equation
\begin{equation}\label{e-theta-1-11-1-11-11}
\frac{d \theta_{\epsilon}^*  (t)}{dt} = e^{-ct}r(u^*( y_{\epsilon}^* (t) , F( y_{\epsilon}^* (t))),y_{\epsilon}^* (t)), \ \ \ \ \ \ \ \ \theta_{\epsilon}^*  (0)=0.
\end{equation}
From (\ref{e-theta-1-11-1-100}) it follows that
$$
 x_{\epsilon}^* (t_{l+1}) =  x_{\epsilon}^* (t_{l}) + \int_{t_l}^{t_{l+1}}X(t,u^*( y_{\epsilon}^* (t), F( y_{\epsilon}^* (t))), y_{\epsilon}^* (t))dt
$$
\vspace{-.3in}
\begin{equation}\label{e-CSO-23-11-1}
= x_{\epsilon}^* (t_{l}) +\epsilon
 \int_{\frac{t_l}{\epsilon}}^{\frac{t_{l+1}}{\epsilon}}X(\epsilon\tau , u^*( y_{\epsilon}^* (\epsilon\tau), F( y_{\epsilon}^* (\epsilon\tau))), y_{\epsilon}^* (\epsilon \tau))d\tau .
\end{equation}
Let $y_{y_{\epsilon}^* (t_l)}(\tau) $ stand for the solution of the  system (\ref{e-CSO-2-SC-I}) considered on the interval $\ [\frac{t_l}{\epsilon},\frac{t_{l+1}}{\epsilon}]  $
with the initial condition  $\ y_{y_{\epsilon}^* (t_l)}(\frac{t_l}{\epsilon})= y_{\epsilon}^* (t_l) $. Then, due to  (\ref{e-CSO-28-11}) and due to continuity
of the function $X(t,u^*(y,F(y)),y) $,
$$
||\frac{\epsilon}{\Delta(\epsilon)}\int_{\frac{t_l}{\epsilon}}^{\frac{t_{l+1}}{\epsilon}}X(\epsilon\tau , u^*( y_{\epsilon}^*(\epsilon\tau), F( y_{\epsilon}^*(\epsilon\tau))), y_{\epsilon}^*(\epsilon \tau))d\tau
$$
\vspace{-.2in}
\begin{equation}\label{e-CSO-23-11-2}
-\
\frac{\epsilon}{\Delta(\epsilon)} \int_{\frac{t_l}{\epsilon}}^{\frac{t_{l+1}}{\epsilon}}X(t_l, u^*( y_{y_{\epsilon}^*(t_l)}(\tau), F( y_{y_{\epsilon}^*(t_l)}(\tau))), y_{y_{\epsilon}^*(t_l)}(\tau))d\tau||\leq \kappa(\epsilon), \ \ \ \ \lim_{\epsilon\rightarrow 0}\kappa(\epsilon)=0
\end{equation}
(note that $\ F( y_{y_{\epsilon}^*(t_l)}(\tau))= F( y_{y_{\epsilon}^*(t_l)}(\frac{t_l}{\epsilon}))\BYDEF z_{\epsilon}^*(t_l) \ \ \forall \tau\in [\frac{t_l}{\epsilon}, \frac{t_{l+1}}{\epsilon}]$). Since the estimate similar to (\ref{e-CSO-28-13}) is valid, one can repeat the argument used in the proof of Proposition \ref{Prop-asym-opt-contr} to establish that the estimates (\ref{e-extra-11}) are valid, where $x_{\epsilon}(l), \ l=0,1,..., $ are defined by (\ref{e-CSO-24-11}). Continuing further as in the proof of Proposition \ref{Prop-asym-opt-contr}, one can establish the validity of (\ref{e-CSO-24-15}) and, subsequently, the validity of (\ref{e-CSO-24-15-2}). The latter along with (\ref{e-CSO-24-15-3}) and (\ref{e-CSO-24-15-3-101}) prove the proposition.
$\ \Box$

\end{document}